\documentclass{article}

\usepackage{arxiv}

\usepackage{doi}

\usepackage{amsmath,amssymb,ifthen,url,graphicx,color,array,theorem}

\usepackage{tikz}
\usepackage{tikz-3dplot}

\usepackage{afterpage}
\usepackage{pgf,pgfplots,pgfplotstable}
\usepackage{pdflscape,rotating}
\usepackage{subcaption}

\usetikzlibrary{shapes.geometric}
\usetikzlibrary{positioning,calc}
\usetikzlibrary{arrows.meta}
\usetikzlibrary{fit, backgrounds}

\usetikzlibrary{intersections,patterns,pgfplots.fillbetween}

\makeatletter
\tikzoption{canvas is xy plane at z}[]{%
  \def\tikz@plane@origin{\pgfpointxyz{0}{0}{#1}}%
  \def\tikz@plane@x{\pgfpointxyz{1}{0}{#1}}%
  \def\tikz@plane@y{\pgfpointxyz{0}{1}{#1}}%
  \tikz@canvas@is@plane
}
\makeatother

\usepackage{natbib}
\usepackage{booktabs,multicol,threeparttable}
\usepackage{multirow}

\newcolumntype{L}[1]{>{\raggedright\let\newline\\\arraybackslash\hspace{0pt}}m{#1}}
\newcolumntype{C}[1]{>{\centering\let\newline\\\arraybackslash\hspace{0pt}}m{#1}}
\newcolumntype{R}[1]{>{\raggedleft\let\newline\\\arraybackslash\hspace{0pt}}m{#1}}

\usepackage{color}
\definecolor{strcolor}{rgb}{0.6, 0.2, 0.6}
\definecolor{commentcolor}{rgb}{0.3125, 0.5, 0.3125}
\definecolor{keycol}{rgb}{0, 0, 1}

\usepackage{xcolor}
\usepackage[linesnumbered,ruled,vlined,noend]{algorithm2e}
\usepackage{sankey}

\SetCommentSty{mycommfont}

\SetKwInput{KwInput}{Input}                
\SetKwInput{KwOutput}{Output}              

\SetInd{0.5em}{0.5em}

\newtheorem{definition}{Definition}

\newcommand{\Distance}[1]
{
	e_{#1}
}

\newcommand{\PlacementNeq}[1]
{
	r_{#1}^{-}
}

\newcommand{\PlacementPos}[1]
{
	r_{#1}^{+}
}

\newcommand{\Orientation}[1]
{
    q_{#1}
}

\newcommand{\Overlap}[2]
{
	m_{#1}^{#2}
}

\newcommand{\Support}[2]
{
	s_{#1}^{#2}
}

\newcommand{\Start}[1]
{
	p_{#1}^{\mathrm{S}}
}

\newcommand{\End}[1]
{
	p_{#1}^{\mathrm{E}}
}

\newcommand{\Diff}[1]
{
	e_{#1}^{\mathrm{S}}
}

\newcommand{\Area}[2]
{
	a_{#1}^{\mathrm{S},#2}
}

\newcommand{\Position}[1]
{
	u_{#1}
}

\newcommand{\Succeed}[1]
{
	u^{\mathrm{succ}}_{#1}
}

\newcommand{\setOfItems}{I}

\newcommand{\setOfOrientations}{O}
\newcommand{\setOfDimensions}{D}

\newcommand{\setOfHeavyItems}{\mathcal{H}}
\newcommand{\setOfLightItems}{\mathcal{L}}

\newcommand{\basicBC}{B\&C\textsuperscript{B}}
\newcommand{\completeBC}{B\&C\textsuperscript{C}}

\newcommand{\allConstraints}{\textit{all constraints}}
\newcommand{\noFragility}{\textit{no fragility}}
\newcommand{\noSupport}{\textit{no support}}
\newcommand{\noLifo}{\textit{no LIFO}}
\newcommand{\loadingOnly}{\textit{loading only}}

\newcommand{\lifoNoSeq}{\textit{LIFO without sequence}}

\title{A branch-and-cut algorithm for vehicle routing problems with three-dimensional loading constraints}

\newif\ifuniqueAffiliation

\ifuniqueAffiliation 
\author{ \href{https://orcid.org/0000-0000-0000-0000}{\includegraphics[scale=0.06]{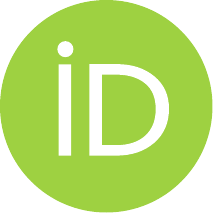}\hspace{1mm}David S.~Hippocampus}\thanks{Use footnote for providing further
		information about author (webpage, alternative
		address)---\emph{not} for acknowledging funding agencies.} \\
	Department of Computer Science\\
	Cranberry-Lemon University\\
	Pittsburgh, PA 15213 \\
	\texttt{hippo@cs.cranberry-lemon.edu} \\
	\And
	\href{https://orcid.org/0000-0000-0000-0000}{\includegraphics[scale=0.06]{orcid.pdf}\hspace{1mm}Elias D.~Striatum} \\
	Department of Electrical Engineering\\
	Mount-Sheikh University\\
	Santa Narimana, Levand \\
	\texttt{stariate@ee.mount-sheikh.edu} \\
}
\else
\usepackage{authblk}

\newbox{\orcid}\sbox{\orcid}{\includegraphics[scale=0.06]{orcid.pdf}} 
\author[1,2]{%
	\href{https://orcid.org/0000-0003-1650-8936}{\usebox{\orcid}\hspace{1mm}Felix Tamke\thanks{\texttt{felix.tamke@tu-dresden.de}}}%
}
\author[1]{%
	\href{https://orcid.org/0000-0002-5869-0425}{\usebox{\orcid}\hspace{1mm}Florian Linß\thanks{\texttt{florian.linss@tu-dresden.de}}}%
}
\author[1,2]{%
	\href{https://orcid.org/0000-0002-9359-8835}{\usebox{\orcid}\hspace{1mm}Leopold Kuttner\thanks{\texttt{leopold.kuttner@tu-dresden.de}}}%
}
\author[1]{%
	\href{https://orcid.org/0000-0003-4711-2184}{\usebox{\orcid}\hspace{1mm}Udo Buscher\thanks{\texttt{udo.buscher@tu-dresden.de}}}%
}
\affil[1]{Faculty of Business and Economics, TU Dresden}
\affil[2]{Clover Optimization}
\fi

\renewcommand{\shorttitle}{A branch-and-cut algorithm for 3L-CVRPs}

\hypersetup{
pdftitle={A branch-and-cut algorithm for vehicle routing problems with three-dimensional loading constraints},
pdfsubject={cs.DM, math.OC},
pdfauthor={Felix Tamke, Florian Linß, Leopold Kuttner, Udo Buscher},
pdfkeywords={vehicle routing; container loading; branch-and-cut; constraint programming},
}

\begin{document}
\maketitle

\begin{abstract}
	This paper presents a new branch-and-cut algorithm based on infeasible path elimination for the three-dimensional loading capacitated vehicle routing problem (3L-CVRP) with different loading problem variants. We show that a previously infeasible route can become feasible by adding a new customer if support constraints are enabled in the loading subproblem and call this the incremental feasibility property. Consequently, different infeasible path definitions apply to different 3L-CVRP variants and we introduce several variant-depending lifting steps to strengthen infeasible path inequalities. The loading subproblem is solved exactly using a flexible constraint programming model to determine the feasibility or infeasibility of a route. An extreme point-based packing heuristic is implemented to reduce time-consuming calls to the exact loading algorithm. Furthermore, we integrate a start solution procedure and periodically combine memoized feasible routes in a set-partitioning-based heuristic to generate new upper bounds. A comprehensive computational study, employing well-known benchmark instances, showcases the significant performance improvements achieved through the algorithmic enhancements. Consequently, we not only prove the optimality of many best-known heuristic solutions for the first time but also introduce new optimal and best solutions for a large number of instances.
\end{abstract}


\section{Introduction}

The classic capacitated vehicle routing problem (CVRP) typically provides only approximate solutions for real-world transportation problems; this is because the loading capacity of a vehicle is approximated one-dimensionally by volume, weight, or pallet space, instead of considering the actual dimensions of the loaded goods being considered in the solving of a multi-dimensional packing problem. The former approach is suitable for bulk cargo and homogeneous goods, while the latter is particularly relevant for transporting heterogeneous goods with defined dimensions, such as pallets, crates, or packages and including shipping freight and parcels, food and beverages for retail, furniture, and pharmaceutical logistics. 
Here, relying on a one-dimensional approximation often results in infeasible or inferior solutions, as the goods may not fit on the vehicle or cargo space may end up being wasted, leading to delays in delivery execution, underutilized vehicles, and generally increased costs.

Two extensions of the CVRP---called two-dimensional and three-dimensional loading capacitated vehicle routing problem (2L-CVRP and 3L-CVRP, respectively)---have been proposed to integrate load and route planning. Both extensions combine two $\mathcal{NP}$-hard problems and are therefore computationally very challenging, but the additional effort is usually worthwhile. \citet{Cote2017} compare an integrated 2L-CVRP and three non-integrated 2L-CVRP approaches, which tackle the planning problems separately. Their study reveals that the routing and loading decisions are closely connected; the non-integrated approaches produce solutions that are, on average, approximately 7\% more costly than those generated by the integrated approach, and up to 26\% more costly in the worst case.

While heuristic and exact solution methods for the 2L-CVRP have been thoroughly studied, the literature on the 3L-CVRP focuses primarily on heuristic solution algorithms, and exact approaches have been studied to a lesser extent: to the best of our knowledge, there is currently no exact algorithm that can optimally solve even the smallest instances of the benchmark set introduced in \citet{Gendreau2006} given their loading constraint definitions. However, exact approaches are useful in determining the structures and properties of optimization problems, and they can provide baselines for evaluating heuristics.  In addition, exact approaches are the best choice both for showcasing the advantages of integrating load optimization into route planning and for highlighting the challenges of solving the integrated problem. To help fill this research gap, this paper makes the following contributions: 

\begin{itemize}
    \item We describe the \textit{incremental feasibility} property for three-dimensional loading problems. Given a minimum support area constraint, this property describes that a previously infeasible loading can become feasible by adding items to the container.
    \item We present a comprehensive branch-and-cut algorithm for multiple 3L-CVRP variants. As part of this, we introduce new inequalities and lifting procedures along with a flexible constraint programming model and an efficient packing metaheuristic to solve the container loading subproblem.
    \item We conduct a computational study for different variants of the container loading subproblem using small to mid-sized instances from the standard 3L-CVRP benchmark set presented in \citet{Gendreau2006}. For the first time, we present proven optimal solutions for many instances as well as provide new best-known solutions in some cases with significant savings compared to heuristic results.
    \item We publish a repository including detailed solutions, a simple web app to visualize solutions and solver statistics, and the source code of our algorithm to facilitate validation and encourage further research.
\end{itemize}

The remainder of this paper is structured as follows. Section~\ref{sec:Problem} describes the integrated problem and the \textit{incremental feasibility} property. Section~\ref{sec:Literature} presents a brief overview of the relevant literature. Section~\ref{sec:VRPModel} provides a two-index vehicle flow formulation for the 3L-CVRP based on infeasible paths. Section~\ref{sec:LoadingProblem} presents an exact constraint programming (CP) formulation and a heuristic approach for the container loading subproblem. In Section~\ref{sec:VI}, several model improvements, including the lifting of infeasible path inequalities, are introduced. Section~\ref{sec:BCAlgorithm}, describes the branch-and-cut algorithm. Section~\ref{sec:ComputationalStudies} presents the computational study. Finally, Section~\ref{sec:Conclusion} concludes the paper with a summary of its main findings.  It should be noted that 'packing' and 'loading' are used interchangeably in this paper. 

\section{The vehicle routing problem with three-dimensional loading constraints}\label{sec:Problem}

\subsection{Problem description}
The definition of the 3L-CVRP by \citet{Gendreau2006} closely follows the classical model of the CVRP. Specifically, a fleet of homogeneous vehicles with given cargo space dimensions and a maximum weight capacity is located at a single depot. Customers request a non-negative number of rectangular goods, called 'items', with specific dimensions and weights. The distance between all customers and the depot and between any two customers is known, and costs are proportional to the distance. The goal is to determine a minimum-cost collection of vehicle routes, with each originating and terminating at the depot, such that every customer is visited exactly once and each loading is feasible subject to the vehicles' weight limit and given loading constraints. The latter is particularly important for the 3L-CVRP and is part of the container loading subproblem.


\subsection{Container loading subproblem} \label{sec:containerLoadingDef}
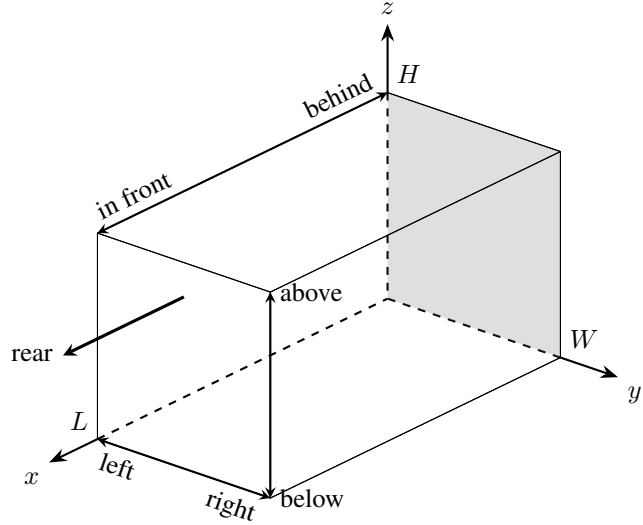
\begin{figure}
    \centering  
\newcommand{\tdplotsetcoordcart}[4]{%
\coordinate (#1) at (#2,#3,#4);
\coordinate (#1xy) at (#2,#3,0);
\coordinate (#1xz) at (#2,0,#4);
\coordinate (#1yz) at (0,#3,#4);
\coordinate (#1x) at (#2,0,0);
\coordinate (#1y) at (0,#3,0);
\coordinate (#1z) at (0,0,#4);
}

\tdplotsetmaincoords{66}{130}
\begin{tikzpicture}[scale=1,tdplot_main_coords,line join=round]
\coordinate (O) at (0,0,0);
\draw[thick,-{Stealth}] (6,0,0) -- (7,0,0) node[anchor=north east]{$x$};
\draw[thick,-{Stealth}] (0,3,0) -- (0,4,0) node[anchor=north west]{$y$};
\draw[thick,-{Stealth}] (0,0,3) -- (0,0,4) node[anchor=south]{$z$};

\tdplotsetcoordcart{Container}{6}{3}{3}


  \draw[stealth-stealth, thick, color=black]
    ($ (Containerz) + (0,0,0) $) --  ($(Containerxz) + (0,0,0)$)
    node[pos=0, anchor=south east, sloped]{behind}
    node[pos=1, anchor=south west, sloped]{in front};
  \draw[stealth-stealth, thick, color=black]
    ($ (Containerx) + (0,0,0) $) --  ($(Containerxy) + (0,0,0)$)
    node[pos=0, anchor=north west, sloped]{left}
    node[pos=1, anchor=north east, sloped]{right};
  \draw[stealth-stealth, thick, color=black]
    ($ (Container) + (0,0,0) $) --  ($(Containerxy) + (0,0,0)$)
    node[pos=0, anchor=west]{above}
    node[pos=1, anchor=west]{below};

\draw[dashed, thick, color=black] (Containerx) -- (O);
\draw[dashed, thick, color=black] (Containery) -- (O);
\draw[dashed, thick, color=black] (Containerz) -- (O);

\node[anchor=south west] at (Containerz) {$H$};
\node[anchor=south west] at (Containery) {$W$};
\node[anchor=south east] at (Containerx) {$L$};

\draw[draw=none, fill=gray, opacity=0.25] (Containeryz) -- (Containery) -- (O) -- (Containerz) -- cycle;


\coordinate (EndOfContainer) at ($(Containerxz)!.5!(Containerxy)$);
\draw[-stealth,color=black, very thick] ($(EndOfContainer) + (0,0,1)$) -- ($(EndOfContainer) + (2.5,0,1)$) node[anchor=east]{rear};




\draw[solid, color=black] (Containerx) -- (Containerxz) -- (Containerz) -- (Containeryz) -- (Containery) -- (Containerxy) -- cycle;

\draw[solid, color=black] (Container) -- (Containerxy) -- (Containerx) -- (Containerxz) -- cycle;
\draw[solid, color=black] (Container) -- (Containerxy) -- (Containery) -- (Containeryz) -- cycle;
\draw[solid, color=black] (Container) -- (Containerxz) -- (Containerz) -- (Containeryz) -- cycle;
\end{tikzpicture}
    \caption{Definition of the coordinate system, the relative directions, and the access direction for unloading.}\label{fig:CoordDefintion}
\end{figure}

The problem of finding a feasible loading of multiple small items in three-dimensional space in a large container is defined as the container loading problem (CLP) 
, cf. \cite{Bortfeldt2013a}. It should be noted that while CLPs usually consider an objective function for optimization, we are only interested in proving feasibility or infeasibility; accordingly, we consider the constraints used in the benchmark instances of \citet{Gendreau2006} and in most heuristic approaches. 

The definitions in this paragraph describe the considered loading restrictions. All of a customer's items have the same group ID, which indicates their position in the unloading sequence and the route; higher group IDs are unloaded first. Relative directions of items are defined as seen from the rear of the container as shown in Figure~\ref{fig:CoordDefintion}. An item $i$ is below, to the left of, or behind another item $j$ if the top face of item $i$ has a z-coordinate lower than or equal to the bottom face of item $j$, if the right face of item $i$ has a y-coordinate lower than or equal to the left face of item $j$, or if the front face of item $i$ has an x-coordinate lower than or equal to the back face of item $j$. 'To the right of', 'above', and 'in front of' are defined analogously. 
Given the definitions above, the following constraints apply to the container loading subproblem of the 3L-CVRP:
\begin{itemize}
    \item No-overlap: Items must not overlap, nor exceed container dimensions.
    \item Rotation: Item rotations around the z-axis are allowed, swapping length and width, resulting in two possible orientations.
    \item Fragility: Items can be fragile or non-fragile. Fragile items can be stacked onto any other item, while non-fragile items must not touch fragile items from above.
    \item Support: Items require a percentage of their base area to be supported by other items or the base of the container.
    \item LIFO: Unloading must respect the last-in-first-out (LIFO) policy. When a customer is visited, each of its items $i$ must be accessible through a sequence of straight movements parallel to the x-axis. Thus, no item $j$ with a lower group ID may be placed between $i$ and the rear of the vehicle or above item $i$ if they intersect in the xy-plane.
\end{itemize}
\cite{Gendreau2006} introduce five different loading variants for the 3L-CVRP to evaluate the influence of each loading constraint: the \allConstraints{} variant considers the full set of loading constraints defined above; the \noFragility{}, \noLifo{}, and \noSupport{} variants each relax their corresponding constraint; and the \loadingOnly{} variant considers only the no-overlap and rotations constraints. Note that without the support constraints, as in the \noSupport{} and \loadingOnly{} variants, zero percent of the bottom face of an item must be supported, allowing it to hover. If the support constraints must be respected, a special property for CLPs occurs: 

\begin{definition}
    Given a container loading problem with support constraints, \textit{incremental feasibility} describes the property that adding an item to an infeasible loading can generate a feasible loading by providing additional support area.
\end{definition}

\begin{figure}
    \centering  
    \begin{subfigure}[b]{0.48\textwidth}
        {
\newcommand{\tdplotsetcoordcart}[4]{%
\coordinate (#1) at (#2,#3,#4);
\coordinate (#1xy) at (#2,#3,0);
\coordinate (#1xz) at (#2,0,#4);
\coordinate (#1yz) at (0,#3,#4);
\coordinate (#1x) at (#2,0,0);
\coordinate (#1y) at (0,#3,0);
\coordinate (#1z) at (0,0,#4);
}

\newcommand{\drawBox}[8]{
\tdplotsetcoordcart{Box}{#1+#4}{#2+#5}{#3+#6}
\draw[solid, color=black, fill=#7, fill opacity=#8] (Box) -- ($(Box)+(0,0,-#6)$) -- ($(Box)+(0,-#5,-#6)$) -- ($(Box)+(0,-#5,0)$) -- cycle;
\draw[solid, color=black, fill=#7, fill opacity=#8] (Box) -- ($(Box)+(0,0,-#6)$) -- ($(Box)+(-#4,0,-#6)$) -- ($(Box)+(-#4,0,0)$) -- cycle;
\draw[solid, color=black, fill=#7, fill opacity=#8] (Box) -- ($(Box)+(0,-#5,0)$) -- ($(Box)+(-#4,-#5,0)$) -- ($(Box)+(-#4,0,0)$) -- cycle;
}

\tdplotsetmaincoords{55}{130} 

\noindent
\begin{tikzpicture}[scale=1,tdplot_main_coords,line join=round]
\coordinate (O) at (0,0,0);

\draw[thick, color = black, -{Stealth}] (4,0,0) -- (4.5,0,0) node[anchor=north east]{$x$};
\draw[thick, color = black, -{Stealth}] (0,2.5,0) -- (0,3,0) node[anchor=north west]{$y$};
\draw[thick, color = black, -{Stealth}] (0,0,1.5) -- (0,0,2) node[anchor=south]{$z$};

\tdplotsetcoordcart{Container}{4}{2.5}{1.5}

\draw[thick, color=black] (Containerx) -- (O);
\draw[thick, color=black] (Containery) -- (O);
\draw[thick, color=black] (Containerz) -- (O);

\draw[solid, thick, color=black] (Containerx) -- (Containerxz) -- (Containerz);
\draw[solid, thick, color=black] (Containery) -- (Containeryz) -- (Containerz);
\draw[solid, thick, color=black] (Containerx) -- (Containerxy) -- (Containery);

\drawBox{0}{0}{0}{1.5}{2}{0.5}{black!80}{0.8}
\drawBox{0}{0}{0.5}{3}{1}{0.5}{gray!40}{0.5}



\end{tikzpicture}}
        \caption{Infeasible loading given support constraints}
    \end{subfigure}
    \hfill
    \begin{subfigure}[b]{0.48\textwidth}
        {
\newcommand{\tdplotsetcoordcart}[4]{%
\coordinate (#1) at (#2,#3,#4);
\coordinate (#1xy) at (#2,#3,0);
\coordinate (#1xz) at (#2,0,#4);
\coordinate (#1yz) at (0,#3,#4);
\coordinate (#1x) at (#2,0,0);
\coordinate (#1y) at (0,#3,0);
\coordinate (#1z) at (0,0,#4);
}

\newcommand{\drawBox}[8]{
\tdplotsetcoordcart{Box}{#1+#4}{#2+#5}{#3+#6}
\draw[solid, color=black, fill=#7, fill opacity=#8] (Box) -- ($(Box)+(0,0,-#6)$) -- ($(Box)+(0,-#5,-#6)$) -- ($(Box)+(0,-#5,0)$) -- cycle;
\draw[solid, color=black, fill=#7, fill opacity=#8] (Box) -- ($(Box)+(0,0,-#6)$) -- ($(Box)+(-#4,0,-#6)$) -- ($(Box)+(-#4,0,0)$) -- cycle;
\draw[solid, color=black, fill=#7, fill opacity=#8] (Box) -- ($(Box)+(0,-#5,0)$) -- ($(Box)+(-#4,-#5,0)$) -- ($(Box)+(-#4,0,0)$) -- cycle;
}

\tdplotsetmaincoords{55}{130} 

\noindent
\begin{tikzpicture}[scale=1,tdplot_main_coords,line join=round]
\coordinate (O) at (0,0,0);

\draw[thick, color = black, -{Stealth}] (4,0,0) -- (4.5,0,0) node[anchor=north east]{$x$};
\draw[thick, color = black, -{Stealth}] (0,2.5,0) -- (0,3,0) node[anchor=north west]{$y$};
\draw[thick, color = black, -{Stealth}] (0,0,1.5) -- (0,0,2) node[anchor=south]{$z$};

\tdplotsetcoordcart{Container}{4}{2.5}{1.5}

\draw[thick, color=black] (Containerx) -- (O);
\draw[thick, color=black] (Containery) -- (O);
\draw[thick, color=black] (Containerz) -- (O);

\draw[solid, thick, color=black] (Containerx) -- (Containerxz) -- (Containerz);
\draw[solid, thick, color=black] (Containery) -- (Containeryz) -- (Containerz);
\draw[solid, thick, color=black] (Containerx) -- (Containerxy) -- (Containery);

\drawBox{0}{0}{0}{1.5}{2}{0.5}{black!80}{0.8}
\drawBox{1.5}{0}{0}{2}{1.5}{0.5}{gray}{0.8}
\drawBox{0}{0}{0.5}{3}{1}{0.5}{gray!40}{0.5}



\end{tikzpicture}}
        \caption{Feasible loading given support constraints}
    \end{subfigure}
    \caption{Visualization of the incremental feasibility property in the 3L-CVRP.} \label{fig:StabilityRequirement}
\end{figure}
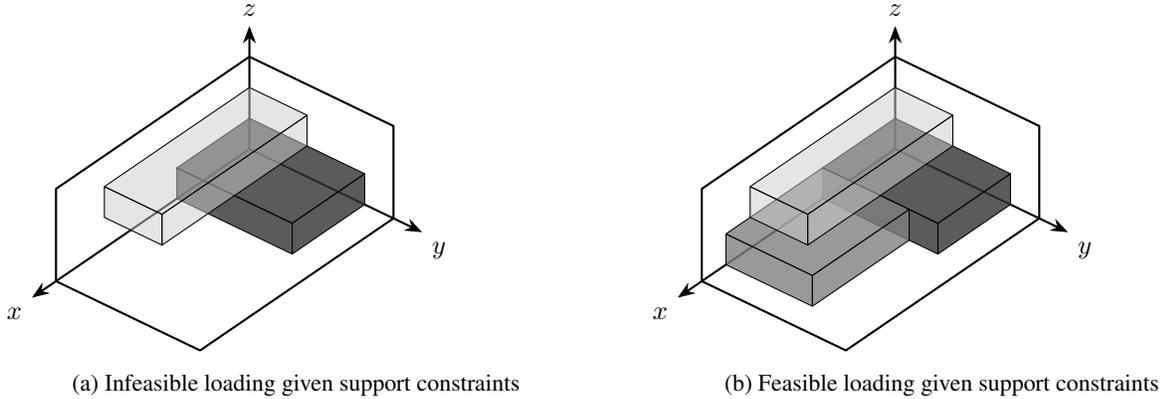

Figure~\ref{fig:StabilityRequirement} shows an example of this property. If we assume that the support area must be 100\%, the loading on the left side is infeasible, as the bottom face of the upper item is not fully supported. Moreover, the two items by themselves can never achieve a feasible loading if they must be stacked. However, the loading becomes feasible if a suitable third item is added. Therefore, only this particular infeasible set of items---not all supersets---can be excluded from further consideration in algorithms. This property represents a key distinguishing characteristic from one- and two-dimensional packing problems, where adding an item cannot affect the infeasibility of a packing.

\section{Literature review}\label{sec:Literature}

This section reviews the main literature related to the 3L-CVRP. After presenting an overview of the 3L-CVRP and its variants that have been studied thus far, we focus on the algorithms used to solve the 3L-CVRP.  

\subsection{Literature on existing 3L-CVRPs}


The first study of the 3L-CVRP is published by \cite{Gendreau2006}, who introduce a standard CVRP with a limited fleet size and incorporate a CLP approach that includes the constraints of: no-overlap, LIFO, fragility, item rotation, and vertical support, as detailed in Section~\ref{sec:Problem}. 
While \cite{Pollaris2015} provide a comprehensive review of the variants existing as of 2015, we focus here on more recent extensions. For one notable variant, \cite{Bortfeldt2020} describe the option of split deliveries; the authors distinguish between forced splitting, where the complete demand cannot be loaded into one vehicle, and optional splitting, where the demand can be split to reduce costs. 

The generalization of vehicle routing problems with loading constraints reveals interesting additional complexities in some well-studied VRP variants that are not present in the CVRP, such as reloading efforts with and without pickup tasks or backhauls. \cite{Mannel2016} study the pickup and delivery problem with three-dimensional loading constraints; similar to the 3L-CVRP, no reloading is allowed here, and so after an item is loaded, it should not be moved before unloading. In contrast, \cite{Reil2018} address the option of reloading in the context of the VRP with backhauls and three-dimensional loading constraints. In addition to rear loading, side loading is also considered.  
 
\cite{Bortfeldt2013a} conduct a comprehensive review of packing problems and their associated constraints, which are relevant in the context of 3L-CVRPs. \cite{Krebs2021} address advanced loading constraints for the 3L-CVRP, such as weight distribution to consider axle weights, robust stability, and reachability constraints.

\subsection{Literature on algorithms for the 3L-CVRP}

Starting with the algorithm of \cite{Gendreau2006}, the literature mainly sees the 3L-CVRP solved heuristically. Most heuristics decompose the integrated problem into its subproblems; the outer routing heuristic aims to improve the routes, while the inner packing heuristic checks the loading feasibility by solving the CLP for these routes.

Common approaches for solving the routing problem are tabu search \citep{Gendreau2006,Tarantilis2009,Bortfeldt2012,Zhu2012,Tao2015}, ant colony optimization \citep{Fuellerer2010}, greedy randomized adaptive search procedure \citep{Lacomme2013}, adaptive variable neighborhood search \citep{Wei2014}, and column-generation-based heuristic algorithms \citep{Mahvash2017,Rajaei2022}. The 'packing first, routing second' heuristic from \cite{Bortfeldt2013} uses a memoization technique to determine loading feasibility by overestimating the required load meters, thereby reducing the number of CLP subproblems that need to be solved. This method operates by strategically executing load meter approximations for both single and multiple customer requests, which are then stored in a tree structure.

The frequently-invoked loading subproblem is primarily addressed with constructive heuristics. Notably, the bottom-left-fill algorithm and the touching area algorithm of \cite{Gendreau2006} are commonly employed for this evaluation. Efficiently limiting the evaluation of possible placement points involves such concepts as the corner-point approach \citep{Martello2000}, the normal-pattern approach \citep{Christofides1977}, and the extreme-point approach \citep{Crainic2008}. The most effective heuristic algorithm for the 3L-CVRP is presented by \cite{Zhang2015}; their algorithm employs an evolutionary local search for routing and an open-space-based heuristic using extreme points for efficient item placement during loading evaluation.

Addressing the challenge of results verification, \cite{Krebs2023a} have recently published a solution validator and visualizer. This is especially crucial, as existing approaches in the literature often present only summarized results, making verification difficult.

With respect to exact solution approaches, the 2L-CVRP benefits from a variety of methods, including state-of-the-art algorithms such as the branch-and-cut algorithm by \cite{Zhang2022} and the branch-and-price-and-cut algorithm by \cite{Zhang2022a}. In contrast, the number of exact formulations available for the 3L-CVRP is relatively limited. \cite{Junqueira2013} propose a first mixed-integer programming (MIP) formulation, solving the integrated problem of the 3L-CVRP. To check vertical stability, the load-bearing strength of the items is taken into account instead of the minimum percentage of the support area; in addition, no rotations are considered. The authors prove optimality for up to 9 customers and 32 items.
\cite{Vega2019} develop a nonlinear MIP formulation for the integrated problem and investigate multiple objectives such as fleet balancing and weight distribution; their analysis utilizes the instances proposed by \cite{Junqueira2013}. \cite{Hokama2016} extend their branch-and-cut approach for the 2L-CVRP so that some three-dimensional loading constraints can also be considered. The routing subproblem is formulated as an MIP, and cuts are added to exclude routes with infeasible loading. The loading checks are performed with a CP model, which considers only the LIFO policy and no rotation, fragility, or vertical stability. The branch-and-cut approach optimally solved the relaxed instances of \cite{Gendreau2006} for up to 30 customers and 32 items. Most recently, \cite{Kuecuek2022} formulate a CP model for the routing problem. The loading is solved using both an evolutionary heuristic and a CP model, but only the no-overlap constraints and rotation are considered. The algorithm is evaluated on the simplified instances of \cite{Moura2009} and can prove optimality for instances of up to 25 customers and 115 items. Notably, none of the existing exact approaches fully account for all the constraints specified by \cite{Gendreau2006}.

\section{Model formulation}\label{sec:VRPModel}
\subsection{Two-index vehicle flow formulation}
We define $G=(N,A)$ as a complete directed graph, where $N$ is the set of nodes $\{0,1,2,\dots,n\}$, with node $0$ representing the depot. $C=\{1,2,\dots,n\}$ denotes the set of customers. The arc set is defined as $A = \{(i,j)|i,j \in N\}$, with $c_{ij}$ denoting the travel cost of arc $(i,j) \in A$. There are $K$ homogeneous vehicles located at the depot; each vehicle has a weight capacity of $Q$ and a cuboid-shaped loading space defined by length $L$, width $W$, and height $H$, with $V = L\cdot W\cdot H$ is the loading volume. Each customer $i$ requires a set of $m_i$ three-dimensional items $I_{ik}$ with $k \in \{1,\dots,m_i\}$, having width $w_{ik}$, length $l_{ik}$, and height $h_{ik}$, total weight $q_i$, and the total required loading volume $v_i = \sum_{k=1}^{m_i} l_{ik} w_{ik}  h_{ik}$.

Given a node $i \in N$, $\delta^-\left(i\right)$ and $\delta^+\left(i\right)$ denote the sets of outgoing and incoming arcs, respectively. Moreover, $A(S) = \{(i, j ) \in A: i , j \in S\}$ is the set of all arcs connecting nodes in an arbitrary set of customers $S \subseteq C$, and $q\left(S\right)$ and $v\left(S\right)$ are the total weight and volume demanded by all $i \in S$, respectively. In addition, $P = \left\{v_1,\dots,v_p\right\}$ represents a sequence of nodes, and $A(P) = \left\{\left(v_i, v_{i+1}\right) |i=1,\dots,p-1\right\}$ denotes the set of all arcs included in $P$. We define three special cases of paths: (i) if $v_1 \neq 0, v_p \neq 0$, $P$ corresponds to a regular path $P^\mathrm{r}$; (ii) if $v_1 \neq 0, v_p=0$, $P$ corresponds to a tail path $P^0$; and (iii) if $v_1=v_p=0$, $P$ corresponds to a complete route $R$. The depot node 0 can only be at the beginning or end of $P$, and $P$ is always a simple path except for its special case $R$. 

Let $x_{ij}$ be a binary variable for each arc $(i,j)\in A$, where $x_{ij}=1$ if $(i,j)$ is used by one of the vehicles and $x_{ij}=0$ otherwise. Thus, we can formulate the two-index vehicle flow formulation for the 3L-CVRP as follows:

\begin{align}
\allowdisplaybreaks
\min &\sum_{\left( i,j \right) \in A} c_{ij} x_{ij}\label{eq:MIP:ObjFunction} \\
\text{s.t.} &\sum_{\left( i, j \right) \in \delta^+\left(i\right)} x_{ij} = 1 \quad \forall~ i \in C, \label{eq:MIP:InDegree}\\
    &\sum_{\left( i, j \right) \in \delta^-\left(i\right)} x_{ij} = 1 \quad \forall~ i \in C, \label{eq:MIP:OutDegree}\\
    &\sum_{\left( i, j \right) \in \delta^+\left( 0 \right)} x_{ij} \leq K, \label{eq:MIP:StartingArcss} \\
    &\sum_{\left( i,j \right) \in A\left( S \right)} x_{ij} \leq |S| - r\left( S \right) \quad \forall ~ S \subseteq C, S \neq \emptyset, \label{eq:MIP:SubtourElim}\\
    &\sum_{\left( i,j \right) \in A\left(P\right)} x_{ij} \leq |A(P)| - 1 \quad \forall ~ \text{infeasible paths } P.  \label{eq:MIP:InfPath}
\end{align}

The objective function~\eqref{eq:MIP:ObjFunction} aims to minimize the distance traveled. Constraints~\eqref{eq:MIP:InDegree} and \eqref{eq:MIP:OutDegree} ensure that every customer is visited exactly once. Constraint~\ref{eq:MIP:StartingArcss} limits the number of routes to the maximum number of available vehicles. Constraints~\eqref{eq:MIP:SubtourElim} are generalized subtour elimination constraints (GSEC), as they ensure that routes are connected to the depot and respect the one-dimensional capacities of volume and weight; here, $r\left(S\right) = \max\left( \left \lceil \frac{q\left( S \right)}{Q} \right \rceil, \left \lceil \frac{v\left( S \right)}{V} \right \rceil \right)$\ imposes a lower bound on the number of vehicles needed to serve the customers in $S$. Finally, constraints~\eqref{eq:MIP:InfPath} define the infeasible path constraints that eliminate all infeasible paths with respect to the loading constraints.

\subsection{Loading problem-dependent definition of infeasible paths} \label{sec:InfeasiblePathDef}
The definition of an infeasible path depends on the loading variant applied in the container loading subproblem (see Section~\ref{sec:containerLoadingDef}).
Usually, an infeasible route $R=\left\{0, v_2,\dots,v_{p-1}, 0 \right\}$ can be excluded by prohibiting only path $P^\mathrm{r}= \left\{v_2,\dots,v_{p-1}\right\}$---i.e., the sequence of customers in $R$ \citep[cf.][]{Ascheuer2000,ropkes-2007,Zhang2022a}. This is because in the considered routing problems, an infeasible route always remains infeasible when a customer is added; consequently, constraints \eqref{eq:MIP:InfPath} exclude $P^\mathrm{r}$ from being part of any route. However, with the \textit{incremental feasibility} property, this only applies to the \noSupport{} and \loadingOnly{} variants: if support constraints must be respected, we could add another customer to $P^\mathrm{r}$, perhaps leading to a feasible path. To maintain the infeasible sequence of customers, a new customer can only be added at the beginning or end. A customer added at the beginning of the route is unloaded first; therefore, its items can never support later unloaded items and simultaneously comply with the LIFO policy, as the later unloaded items would have to be above the first items. In contrast, a customer added at the end is unloaded last, and its items could always support items unloaded earlier and satisfy the LIFO policy. As a result, $P$ is defined as a tail path $P^0=\left\{v_2,\dots,v_{p-1},0\right\}$, to prohibit infeasible customer sequences without excluding feasible routes given support and LIFO constraints. For the \noLifo{} variant, the depot node must also be added at the beginning, thus prohibiting only the complete route $R$. To distinguish the different versions of constraints \eqref{eq:MIP:InfPath} in the following sections, we call them \textit{infeasible regular path} (IRP) inequality with $P = P^\mathrm{r}$, \textit{infeasible tail path} (ITP) inequality with $P = P^0$, and \textit{infeasible route} (IR) inequality with $P = R$. We can state the following transitive relation concerning the strengths of the three path elimination inequalities: IRP $\succ$ ITP $\succ$ IR. It should be noted that weaker constraints could replace stronger ones without eliminating feasible solutions. Information regarding the infeasible paths of all variants is summarized in Table~\ref{tab:InfPathConfigs}. 

\begin{table}[h]
    \centering
    \footnotesize
    \caption{Definitions of infeasible paths for an infeasible route regarding different loading variants.}\label{tab:InfPathConfigs}
    \begin{tabular}{lll}
    \toprule
    Path $P$ & Variant(s) & Name of constraints \eqref{eq:MIP:InfPath} with $P$ \\
    \midrule
    $P^\mathrm{r} = P$ with $v_1\neq0, v_p\neq0$ & \noSupport{}, \loadingOnly{} & infeasible regular path (IRP) inequality\\
    $P^0 = P$ with $v_1 \neq0, v_p=0$ & \allConstraints{}, \noFragility{} &  infeasible tail path (ITP) inequality\\
    $R = P$ with $v_1=v_p=0$ & \noLifo{} & infeasible route (IR) inequality \\
    \bottomrule
    \end{tabular}
\end{table}

\section{Container loading problem}\label{sec:LoadingProblem}

This section describes the solution methods for the CLP. First, we present the CP model formulation, which can be applied to several variants of the loading problem. In addition, the extended packing heuristic based on extreme points is described so as to quickly find a feasible solution.

\subsection{Constraint programming model}


The CLP can be described by the CP model presented in this section. Variable and parameter descriptions are provided in Table \ref{tab:VariablesCP}. The flexible model can be used for several common container loading variants, as can be seen in Table \ref{tab:PackingProblemVariants}.

\begin{table}[h]
    \centering
    \caption{Definition of variables, parameters, and sets for the CP model.} \label{tab:VariablesCP}
    \begin{tabular}{lll}
    \toprule
    Variable & Domain & Description\\
    \midrule
        $\Start{id}$ & integer & start of item $i$ in dimension $d$ \\  
        $\End{id}$ & integer & end of item $i$ in dimension $d$ \\
        $\Distance{id}$ & integer & spatial extent of item $i$ in dimension $d$ \\
        $\Orientation{io}$ & boolean & item $i$ is placed in orientation $o$ \\
        $\PlacementNeq{ijd}$ & boolean & item $i$ is placed in negative direction of dimension $d$ relative to item $j$ \\
        $\PlacementPos{ijd}$ & boolean & item $i$ is placed in positive direction of dimension $d$ relative to item $j$ \\
        $\Overlap{ij}{\mathrm{xy}}$ & boolean & item $i$ and item $j$ overlap in $\mathrm{xy}$-plane \\
        $\Support{ij}{\mathrm{xy}}$ & boolean & item $i$ is supported by item $j$ in $\mathrm{xy}$-plane\\        
        $\Diff{ijd}$ & integer & relevant overlapping range of items $i$ and $j$ in dimension $d$\\
        $\Area{ij}{\mathrm{xy}}$ & integer & size of the overlapping area of items $i$ and $j$ in $\mathrm{xy}$-plane \\
        $\Position{f}$ & integer & position of customer $f$ in route \\
        $\Succeed{fh}$ & boolean & customer $f$ is visited after customer $h$\\[9pt]
        Parameter \\ \midrule
        $A_i$ & integer & area of item $i$ \\
        $\alpha$ & $[0, 1]$ & support area requirement as a percentage of an item's base area \\
        $l_{io}$, $w_{io}$ & integer & length and width of item $i$ in orientation $o$ \\
        $h_i$ & integer & height of item $i$ \\
        $L, W, H$ & integer & length, width, and height of the container \\
        Set \\ \midrule
        $\setOfItems$ & & set of items \\
        $\setOfOrientations$ & & set of orientations \\
        $\setOfDimensions$ & & set of dimensions \\
        \bottomrule
    \end{tabular}
\end{table}

\subsubsection{Item positions}
The start and end positions $\Start{id}$ and $\End{id}$ of item $i$ in dimension $d \in D$ are defined on a set of placement points $\mathcal{P}_d$, i.e., $\Start{id} \in P_d, \End{id} \in P_d$. Different sets of placement points are valid depending on the problem variant (cf. Table~\ref{tab:PackingProblemVariants}). We distinguish between three different sets of placement points. The simplest considers every coordinate as a placement point according to unit discretization $\mathcal{U}_i$ \citep{Kurpel2020}. Placement points can be reduced by taking item dimensions into account with regular normal patterns $\mathcal{B}_i$ \citep{Boschetti2002}. An even more efficient set of placement points can be determined by minimal meet-in-the-middle patterns $\mathcal{M}_i^{\mathrm{1,2}}$. We use minimal meet-in-the-middle patterns with preprocessing steps 1 and 2 and the placement point union as the threshold selection function \citep{Cote2018}. 
The start and end points of item $i$ in dimension $d$ are connected in the way shown in equation \eqref{eq:CP:Intervals}.
\allowdisplaybreaks
\begin{align}
    \Start{id} + \Distance{id} = \End{id}, \quad \forall ~ i \in \setOfItems, d \in \setOfDimensions. \label{eq:CP:Intervals}
\end{align}

\subsubsection{No-overlap of items}
Relative directions $\PlacementNeq{ijd}$ and $\PlacementPos{ijd}$ of items $i$ and $j$ in dimension $d$ are determined in constraints \eqref{eq:CP:RelativeDirectionsA}--\eqref{eq:CP:RelativeDirectionsSymmetryHint}.
\begin{align}
    \PlacementNeq{ijd} \Rightarrow \End{jd} \leq \Start{id}, \quad \forall ~ i, j \in \setOfItems, i\neq j, d \in \setOfDimensions, \label{eq:CP:RelativeDirectionsA} \\
    \neg \PlacementNeq{ijd} \Rightarrow \Start{id} < \End{jd}, \quad \forall ~ i, j \in \setOfItems, i\neq j, d \in \setOfDimensions, \label{eq:CP:RelativeDirectionsB} \\
    \PlacementPos{ijd} \Rightarrow \End{id} \leq \Start{jd}, \quad \forall ~ i, j \in \setOfItems, i\neq j, d \in \setOfDimensions, \label{eq:CP:RelativeDirectionsC} \\
    \neg \PlacementPos{ijd} \Rightarrow \Start{jd} < \End{id}, \quad \forall ~ i, j \in \setOfItems, i\neq j, d \in \setOfDimensions, \label{eq:CP:RelativeDirectionsD} \\
    \PlacementNeq{ijd} = \PlacementPos{jid}, \quad \forall ~ i, j \in \setOfItems, i\neq j, d \in \setOfDimensions. \label{eq:CP:RelativeDirectionsSymmetryHint}
\end{align}
The relative directions of items for the different dimensions are: $\PlacementNeq{ijd}$ is to the left (y-axis), behind (x-axis), or below (z-axis); $\PlacementPos{ijd}$ is to the right (y-axis), in front (x-axis), or above (z-axis). Using the relative direction variables, constraints \eqref{eq:CP:RelativeDirectionsBoolOr} guarantee that two items do not overlap in at least one dimension.
\begin{align}
        \bigvee\limits_{d \in \setOfDimensions} \left( \PlacementNeq{ijd} \lor \PlacementPos{ijd} \right) \quad \forall ~ i, j \in \setOfItems, i\neq j. \label{eq:CP:RelativeDirectionsBoolOr}
\end{align}

\subsubsection{Item orientation}
Constraints \eqref{eq:CP:OrientationLength}--\eqref{eq:CP:SelectOrientation} set the distance $\Distance{id}$ of an item $i$ in the three dimensions. Since rotation is only allowed in the $\mathrm{xy}$-plane, the distances along the x-axis and the y-axis depend on the chosen orientation $\Orientation{io}$, while the height is fixed. Exactly one orientation must be selected.
\begin{align}
    \Orientation{io} \Rightarrow \Distance{i\mathrm{x}} = l_{io}, \quad \forall ~ i \in \setOfItems, o \in \setOfOrientations, \label{eq:CP:OrientationLength} \\
    \Orientation{io} \Rightarrow \Distance{i\mathrm{y}} = w_{io}, \quad \forall ~ i \in \setOfItems, o \in \setOfOrientations, \label{eq:CP:OrientationWidth} \\
    \Distance{i\mathrm{z}} = h_i, \quad \forall ~ i \in \setOfItems, \label{eq:CP:OrientationHeight} \\
    \sum_{o \in \setOfOrientations} \Orientation{io} = 1, \quad \forall ~ i \in \setOfItems. \label{eq:CP:SelectOrientation}
\end{align}
\subsubsection{Item intersection and support in $\mathrm{xy}$-plane}
Constraints \eqref{eq:CP:IntersectionBoolOr}--\eqref{eq:CP:IntersectionD} determine whether two items overlap in the $\mathrm{xy}$-plane. Variables $\Overlap{ij}{\mathrm{xy}}$ must be true if items $i$ and $j$ are neither to the left of, to the right of, behind, nor in front of each other.
\begin{align}
    \Overlap{ij}{\mathrm{xy}} \lor \PlacementNeq{ij\mathrm{x}} \lor \PlacementPos{ij\mathrm{x}} \lor \PlacementNeq{ij\mathrm{y}} \lor \PlacementPos{ij\mathrm{y}}, \quad \forall i,j \in \setOfItems, j > i, \label{eq:CP:IntersectionBoolOr} \\
    \PlacementNeq{ij\mathrm{x}} \Rightarrow \neg \Overlap{ij}{\mathrm{xy}}, \quad \forall i,j \in \setOfItems, j > i, \label{eq:CP:IntersectionA} \\
    \PlacementPos{ij\mathrm{x}} \Rightarrow \neg \Overlap{ij}{\mathrm{xy}}, \quad \forall i,j \in \setOfItems, j > i, \label{eq:CP:IntersectionB} \\
    \PlacementNeq{ij\mathrm{y}} \Rightarrow \neg \Overlap{ij}{\mathrm{xy}}, \quad \forall i,j \in \setOfItems, j > i, \label{eq:CP:IntersectionC} \\
    \PlacementPos{ij\mathrm{y}} \Rightarrow \neg \Overlap{ij}{\mathrm{xy}}, \quad \forall i,j \in \setOfItems, j > i. \label{eq:CP:IntersectionD}
\end{align}
In addition, constraints~\eqref{eq:CP:SupportA}--\eqref{eq:CP:SupportD} ensure that the vertical stability of item $i$ can only be supported by another item $j$ if they overlap and the top face of $j$ touches the bottom face of $i$. Moreover, an item cannot support itself.
\begin{align}
   \Support{ij}{\mathrm{xy}} \lor \left( \End{j\mathrm{z}} \neq \Start{i\mathrm{z}} \right) \lor \neg \Overlap{ij}{\mathrm{xy}}, \quad \forall ~ i,j \in \setOfItems, j > i, \label{eq:CP:SupportA} \\
   \neg \Overlap{ij}{\mathrm{xy}} \Rightarrow \neg \Support{ij}{\mathrm{xy}}, \quad \forall ~ i,j \in \setOfItems, j > i, \label{eq:CP:SupportB} \\
   \Support{ij}{\mathrm{xy}} \lor \left( \End{j\mathrm{z}} \neq \Start{i\mathrm{z}} \right) \lor \neg \Overlap{ji}{\mathrm{xy}}, \quad \forall ~ i,j \in \setOfItems, i > j, \label{eq:CP:SupportC} \\
   \neg \Overlap{ji}{\mathrm{xy}} \Rightarrow \neg \Support{ij}{\mathrm{xy}}, \quad \forall ~ i,j \in \setOfItems, i > j, \label{eq:CP:SupportD}\\
   \Support{ii}{\mathrm{xy}} = 0,  \quad \forall ~ i \in \setOfItems. \label{eq:CP:SupportE}
\end{align}
\subsubsection{Item fragility}
Item $j$ cannot support the vertical stability of item $i$ if $j$ is fragile and $i$ is non-fragile. Thus, both items cannot simultaneously overlap in the $\mathrm{xy}$-plane while $j$ is directly below $i$ (see constraints~\eqref{eq:CP:SupportA}~and~\eqref{eq:CP:SupportC}).
\begin{align}
    \Support{ij}{\mathrm{xy}} = 0, \quad \forall ~ i,j \in \setOfItems, j \text{ is fragile}, i \text{ is non-fragile}. \label{eq:CP:Fragility}
\end{align}

\subsubsection{Support area in the $\mathrm{xy}$-plane}
Constraints \eqref{eq:CP:SupportAreaDiffX}--\eqref{eq:CP:SupportAreaEnforceDiffY} determine the overlap length of two items $i$ and $j$ along the x-axis, $\Diff{ij\mathrm{x}}$, and the y-axis, $\Diff{ij\mathrm{y}}$. Both lengths are only relevant for the vertical support area if $i$ and $j$ overlap in the $\mathrm{xy}$-plane. Hence, $\End{id} - \Start{jd}$ and $\End{jd} - \Start{id}$ always evaluate to positive values for $d \in \{x,y\}$ in \eqref{eq:CP:SupportAreaDiffX} and \eqref{eq:CP:SupportAreaDiffY}.
\begin{align}
    &\begin{aligned}        
    \Overlap{ij}{\mathrm{xy}} \Rightarrow \Diff{ij\mathrm{x}} = \min\left\{ \End{i\mathrm{x}} - \Start{j\mathrm{x}}, \End{j\mathrm{x}} - \Start{i\mathrm{x}}, \Distance{i\mathrm{x}}, \Distance{j\mathrm{x}}  \right\}, \\ \forall ~ i,j \in \setOfItems, j > i, \label{eq:CP:SupportAreaDiffX}
    \end{aligned} \\
    &\neg \Overlap{ij}{\mathrm{xy}} \Rightarrow \Diff{ij\mathrm{x}} = 0, \ \forall ~ i,j \in \setOfItems, j > i, \label{eq:CP:SupportAreaEnforceDiffX} \\
    &\begin{aligned}
    \Overlap{ij}{\mathrm{xy}} \Rightarrow \Diff{ij}{\mathrm{y}} = \min\left\{ \End{i}{\mathrm{y}} - \Start{j}{\mathrm{y}}, \End{j}{\mathrm{y}} - \Start{i}{}, \Distance{i}{\mathrm{y}}, \Distance{j}{\mathrm{y}}  \right\}, \\ \forall ~ i,j \in \setOfItems, j > i, \label{eq:CP:SupportAreaDiffY}
    \end{aligned} \\
    &\neg \Overlap{ij}{\mathrm{xy}} \Rightarrow \Diff{ij\mathrm{y}} = 0, \ \forall ~ i,j \in \setOfItems, j > i. \label{eq:CP:SupportAreaEnforceDiffY}
\end{align}

The non-negative overlap lengths are used to determine the overlapping area in the $\mathrm{xy}$-plane in constraints \eqref{eq:CP:SupportArea}. Finally, constraints \eqref{eq:CP:SupportAreaRequirement} ensure the minimum support area requirement $\alpha A_i$ for item $i$, where $\alpha \in [0, 1]$. The overlapping area of $i$ and $j$ can only be used for the vertical stability of $i$ if $j$ supports $i$. However, if $i$ is placed on the container floor, the vertical stability requirement is automatically satisfied.
\begin{align}
    & \Area{ij}{\mathrm{xy}} = \Diff{ij\mathrm{x}} \cdot \Diff{ij\mathrm{y}}, \quad \forall ~ i,j \in \setOfItems, j > i, \label{eq:CP:SupportArea} \\
    & \begin{aligned}
    \Start{i\mathrm{z}} > 0 \Rightarrow & \sum_{j \in \setOfItems, j>i} \Area{ij}{\mathrm{xy}} \cdot \Support{ij}{\mathrm{xy}} \\ &+ \sum_{j \in \setOfItems, i>j} \Area{ji}{\mathrm{xy}} \cdot \Support{ij}{\mathrm{xy}} \geq \alpha \cdot A_i, \quad \forall ~ i \in \setOfItems \label{eq:CP:SupportAreaRequirement}.
    \end{aligned}
\end{align}

\subsubsection{LIFO unloading policy}
Constraints \eqref{eq:CP:LifoAllDifferent}--\eqref{eq:CP:LifoB} guarantee that the LIFO unloading policy is respected for a given sequence of customers and also for a set of customers $C^{\mathrm{r}} = \{1,\hdots, \bar{n}\}$ visited in a route. The first three constraints address the customers' positions in a route and the resulting sequence; thus, they determine the group IDs. If customer $f$ succeeds customer $h$ in the route---i.e., $\Succeed{fh}$ is true---customer $f$ is visited after customer $h$ and $\Position{f}$ must be larger than $\Position{h}$; if $\Succeed{fh}$ is false, customer $f$ is visited before customer $h$. For a given customer sequence, the values for $\Succeed{fh}$ and $\Position{f}$ are fixed and not part of the decision to be made. The LIFO unloading policy is realized in constraints \eqref{eq:CP:LifoA} and \eqref{eq:CP:LifoB}, with $c_i$ denoting the customer associated with item $i$. Both constraints are necessary, since $\Succeed{{c_i}, {c_j}}$ is only defined for $c_i > c_j$. If the vehicle serves the customer associated with item $i$ after the customer of item $j$---i.e., $\Succeed{{c_i}, {c_j}} = 1$ in \eqref{eq:CP:LifoA} and $\Succeed{{c_j}, {c_i}} = 0$ in \eqref{eq:CP:LifoB}---item $j$ must be unloaded before item $i$. Therefore, if item $i$ is not placed to the left or right of item $j$, item $i$ must be placed behind or below item $j$, in order to ensure accessibility to item $j$ from the rear when the vehicle arrives at customer $c_j$.
\begin{align}
    \text{alldifferent}(\Position{1}, \hdots, \Position{\bar{n}}), \label{eq:CP:LifoAllDifferent} \\
    \Succeed{fh} \Rightarrow \Position{f} > \Position{h}, \ \forall ~ f,h \in {C}^{\mathrm{r}}, h > f, \label{eq:CP:LifoSuccessor} \\
    \neg \Succeed{fh} \Rightarrow \Position{f} < \Position{h}, \ \forall ~ f,h \in {C}^{\mathrm{r}}, h > f, \label{eq:CP:LifoPredecessor} \\
    \neg \PlacementNeq{ij\mathrm{y}} \land \neg \PlacementPos{ij\mathrm{y}} \land \Succeed{c_i, c_j} \Rightarrow \PlacementNeq{ij\mathrm{x}} \lor \PlacementNeq{ij\mathrm{z}}, \ \forall i,j \in \setOfItems, c_i < c_j, \label{eq:CP:LifoA} \\
    \neg \PlacementNeq{ij\mathrm{y}} \land \neg \PlacementPos{ij\mathrm{y}} \land \neg \Succeed{c_j, c_i} \Rightarrow \PlacementNeq{ij\mathrm{x}} \lor \PlacementNeq{ij\mathrm{z}}, \ \forall i,j \in \setOfItems, c_j < c_i. \label{eq:CP:LifoB}
\end{align}

\subsubsection{Model configurations for different loading problem variants}

The structure of the CP model allows for adaptation to different loading problem variants, such that constraints can be activated or deactivated. Table~\ref{tab:PackingProblemVariants} provides an overview of relevant model formulations and usages of placement points for each dimension. In addition to the known loading problem variants, it is also possible to solve the relaxations considering only the LIFO constraint with and without a given sequence, which are relevant for evaluating integer solutions in the branch-and-cut algorithm.

\begin{table}[]
    \centering
    \footnotesize
    \caption{Different model variants with the corresponding model formulation and valid placement points.} \label{tab:PackingProblemVariants}
    \centering
    \begin{threeparttable}
    \begin{tabular}{lcccccc@{\qquad}c@{\qquad}ccc@{\qquad}c}%
      \toprule
        \\[-3pt]
        & & & &  &  &  &  & \multicolumn{3}{c}{Placement points $\mathcal{P}$} \\
         \cmidrule(rl{3pt}){9-11} \\[-9pt]
         Variant & \rotatebox{45}{\rlap{No-overlap}} & \rotatebox{45}{\rlap{Rotation}} & \rotatebox{45}{\rlap{Support}} & \rotatebox{45}{\rlap{Fragility}} & \rotatebox{45}{\rlap{LIFO}} & \rotatebox{45}{\rlap{Sequence}} & Model & X & Y & Z &  \\ \midrule 
         \allConstraints{} & $\bullet$ & $\bullet$ & $\bullet$ & $\bullet$ & $\bullet$ & $\bullet$ &  \eqref{eq:CP:Intervals}--\eqref{eq:CP:LifoB} & $\mathcal{U}_i$ & $\mathcal{U}_i$ & $\mathcal{B}_i$ \\
         \noFragility{} & $\bullet$ & $\bullet$ & $\bullet$ & & $\bullet$ & $\bullet$ & \eqref{eq:CP:Intervals}--\eqref{eq:CP:SupportE}, \eqref{eq:CP:SupportAreaDiffX}--\eqref{eq:CP:LifoB} & $\mathcal{U}_i$ & $\mathcal{U}_i$ & $\mathcal{B}_i$ \\
         \noLifo{} & $\bullet$ & $\bullet$ & $\bullet$ & $\bullet$ & & & \eqref{eq:CP:Intervals}--\eqref{eq:CP:SupportAreaRequirement} & $\mathcal{U}_i$ & $\mathcal{U}_i$ & $\mathcal{B}_i$ \\
         \noSupport{} & $\bullet$ & $\bullet$ & & $\bullet$ & $\bullet$ & $\bullet$ & \eqref{eq:CP:Intervals}--\eqref{eq:CP:Fragility}, \eqref{eq:CP:LifoAllDifferent}--\eqref{eq:CP:LifoB} & $\mathcal{B}_i$ & $\mathcal{B}_i$ & $\mathcal{U}_i$ \\
         \loadingOnly{} & $\bullet$ & $\bullet$ & & & & & \eqref{eq:CP:Intervals}--\eqref{eq:CP:IntersectionD} & $\mathcal{M}_i^{\mathrm{1,2}}$ & $\mathcal{M}_i^{\mathrm{1,2}}$ & $\mathcal{M}_i^{\mathrm{1,2}}$ \\[6pt]
         Relaxations \\ \midrule
         LIFO without sequence & $\bullet$ & $\bullet$ & & & $\bullet$ & & \eqref{eq:CP:Intervals}--\eqref{eq:CP:SelectOrientation}, \eqref{eq:CP:LifoAllDifferent}--\eqref{eq:CP:LifoB} & $\mathcal{B}_i$ & $\mathcal{B}_i$ & $\mathcal{B}_i$ \\
         LIFO with sequence & $\bullet$ & $\bullet$ & & & $\bullet$ & $\bullet$ & \eqref{eq:CP:Intervals}--\eqref{eq:CP:SelectOrientation}, \eqref{eq:CP:LifoA}--\eqref{eq:CP:LifoB} & $\mathcal{B}_i$ & $\mathcal{B}_i$ & $\mathcal{B}_i$ \\
         \bottomrule
    \end{tabular}
    \begin{tablenotes}[para,flushleft]
        \item[] The placement points are unit discretization ($\mathcal{U}_i$), regular normal patterns ($\mathcal{B}_i$),  and minimal meet-in-the-middle patterns ($\mathcal{M}_i^{\mathrm{1,2}}$).
    \end{tablenotes}
    \end{threeparttable}
\end{table}

\subsection{Augmented extreme point-based packing heuristic}

To efficiently determine feasible solutions to the CLP, a constructive packing heuristic that is similar to \cite{Crainic2008} (cf. Algorithm \ref{alg:ExtremePointHeuristic}) is used within an iterated greedy metaheuristic. 

\begin{algorithm}[h]
	\DontPrintSemicolon
	
	\KwInput{items $I$, dimensions of bin $B$, item placement sequence $S$}
	\KwOutput{packing of items $I$ into bin $B$}
	
	Initialize normal patterns on container floor \tcp*{According to \citet{Beasley1985a}.}
	\ForEach{\upshape item $i$ $\subset I$ according to $S$}
	{
		Select placement for $i$ \tcp*{Best fit lexicographic x, y, z, and short side in x-direction.}
		Perform normal pattern augmentation\;
		Update extreme points \tcp*{According to \citet{Crainic2008}.}
	}
	
	\Return $I$
	\caption{Augmented extreme point-based heuristic}
	\label{alg:ExtremePointHeuristic}
\end{algorithm}

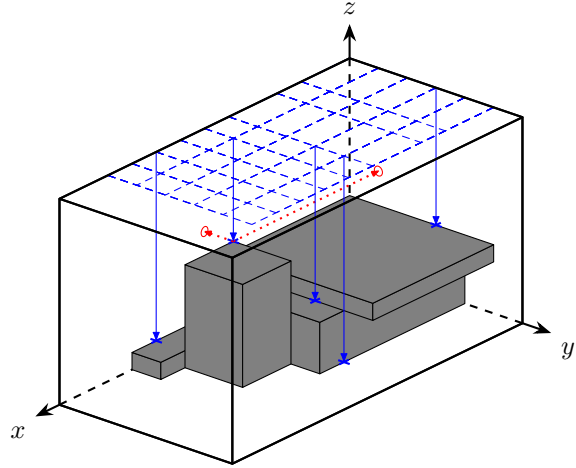
\begin{figure}
    \centering  
    {
\newcommand{\tdplotsetcoordcart}[4]{%
\coordinate (#1) at (#2,#3,#4);
\coordinate (#1xy) at (#2,#3,0);
\coordinate (#1xz) at (#2,0,#4);
\coordinate (#1yz) at (0,#3,#4);
\coordinate (#1x) at (#2,0,0);
\coordinate (#1y) at (0,#3,0);
\coordinate (#1z) at (0,0,#4);
}

\tikzset{
   placementpointEP/.pic={
	   \begin{scope}[canvas is xy plane at z=0]
	     \draw [red] (0,0) circle (2pt);
   		\end{scope}
   },
   placementpointEPX/.pic={
	   \begin{scope}[canvas is xz plane at y=0]
	     \draw [red] (0,0) circle (2pt);
   		\end{scope}
   },
   placementpointEPY/.pic={
	   \begin{scope}[canvas is yz plane at x=0]
	     \draw [red] (0,0) circle (2pt);
	     \draw [red] (0,0) circle (0.1pt);
   		\end{scope}
   },
   placementpointEPZ/.pic={
	   \begin{scope}[canvas is xy plane at z=0]
	     \draw [red] (0,0) circle (2pt);
   		\end{scope}
   },
   placementpointNP/.pic={
     \draw [thick, draw=blue] (0,-0.1,0) -- (0,0.1,0);     
     \draw [thick, draw=blue] (-0.1,0,0) -- (0.1,0,0);
   }
}

\tdplotsetmaincoords{66}{130}
\begin{tikzpicture}[scale=1,tdplot_main_coords,line join=round]
\coordinate (O) at (0,0,0);
\draw[thick,-{Stealth}] (6,0,0) -- (6.5,0,0) node[anchor=north east]{$x$};
\draw[thick,-{Stealth}] (0,3,0) -- (0,3.5,0) node[anchor=north west]{$y$};
\draw[thick,-{Stealth}] (0,0,3) -- (0,0,3.5) node[anchor=south]{$z$};

\tdplotsetcoordcart{Container}{6}{3}{3}


\draw[dashed, thick, color=black] (Containerx) -- (O);
\draw[dashed, thick, color=black] (Containery) -- (O);
\draw[dashed, thick, color=black] (Containerz) -- (O);


\tdplotsetcoordcart{Box1}{3}{2}{0.75}
\draw[solid, color=black, fill=gray, join=round] (Box1) -- (Box1xy) -- (Box1x) -- (Box1xz) -- cycle;
\draw[solid, color=black, fill=gray] (Box1) -- (Box1xy) -- (Box1y) -- (Box1yz) -- cycle;
\draw[solid, color=black, fill=gray] (Box1) -- (Box1xz) -- (Box1z) -- (Box1yz) -- cycle;

\tdplotsetcoordcart{Box2}{2.5}{2.5}{1}
\draw[solid, color=black, fill=gray] (Box2) -- ($(Box2xy)+(0,0,0.75)$) -- ($(Box2x)+(0,0,0.75)$) -- (Box2xz) -- cycle;
\draw[solid, color=black, fill=gray] (Box2) -- ($(Box2xy)+(0,0,0.75)$) -- ($(Box2y)+(0,0,0.75)$) -- (Box2yz) -- cycle;
\draw[solid, color=black, fill=gray] (Box2) -- (Box2xz) -- (Box2z) -- (Box2yz) -- cycle;



\tdplotsetcoordcart{Box3}{4.5}{0.5}{0.25}
\draw[solid, color=black, fill=gray] (Box3) -- ($(Box3xy)+(0,0,0)$) -- ($(Box3x)+(0,0,0)$) -- (Box3xz) -- cycle;
\draw[solid, color=black, fill=gray] (Box3) -- ($(Box3xy)+(0,0,0)$) -- ($(Box3y)+(3,0,0)$) -- ($(Box3yz)+(3,0,0)$) -- cycle;
\draw[solid, color=black, fill=gray] (Box3) -- (Box3xz) -- ($(Box3z)+(3,0,0)$) -- ($(Box3yz)+(3,0,0)$) -- cycle;

\tdplotsetcoordcart{Box4}{4}{1.5}{1.5}
\draw[solid, color=black, fill=gray] (Box4) -- ($(Box4xy)+(0,0,0)$) -- ($(Box4x)+(0,0.5,0)$) -- ($(Box4xz)+(0,0.5,0)$) -- cycle;
\draw[solid, color=black, fill=gray] (Box4) -- ($(Box4xy)+(0,0,0)$) -- ($(Box4y)+(3,0,0)$) -- ($(Box4yz)+(3,0,0)$) -- cycle;
\draw[solid, color=black, fill=gray] (Box4) -- ($(Box4xz)+(0,0.5,0)$) -- ($(Box4z)+(3,0.5,0)$) -- ($(Box4yz)+(3,0,0)$) -- cycle;

\pic at (2.5,2,0) {placementpointNP};

\pic at ($(Box4)+(-1,-1,0)$) {placementpointNP};

\pic at (0,1.5,1) {placementpointNP};

\pic at (4,0,0.25) {placementpointNP};
\pic at (2.5,1.5,0.75) {placementpointNP};

\draw [thick, dotted, red, -{Latex[round,length=4pt]}] ($(Box4)+(-1,-1,0)$) -- (0,0.5,1.5);
\pic at (0,0.5,1.5) {placementpointEPY};
\draw [thick, dotted, red, -{Latex[round,length=4pt]}] ($(Box4)+(-1,-1,0)$) -- (3,0,1.5);
\pic at (3,0,1.5) {placementpointEPX};


\pgfmathsetmacro{\zOffset}{-2}
\pgfmathsetmacro{\zOffsetNew}{3}

%
%
%

\foreach \x in {0,1,1.5,2.5,3,3.5,4,4.5,5}
	\foreach \y in {0,0.5,1,1.5,2,2.5} 
   	{
		\draw [thin, densely dashed, blue] (\x,0,\zOffsetNew) -- (\x,\y,\zOffsetNew);
		\draw [thin, densely dashed, blue] (0,\y,\zOffsetNew) -- (\x,\y,\zOffsetNew);
	} 


\draw [blue, -{Latex[round,length=4pt]}, opacity=1.0] (4,0,\zOffsetNew) -- (4,0,0.25);
\draw [blue, -{Latex[round,length=4pt]}, opacity=1.0] (0,1.5,\zOffsetNew) -- (0,1.5,1);
\draw [blue, -{Latex[round,length=4pt]}, opacity=1.0] (3,0.5,\zOffsetNew) -- (3,0.5,1.5);
\draw [blue, -{Latex[round,length=4pt]}, opacity=1.0] (2.5,1.5,\zOffsetNew) -- (2.5,1.5,0.75);
\draw [-{Latex[round,length=4pt]}, blue] (2.5,2,\zOffsetNew) -- (2.5,2,0);

\draw[solid, thick, color=black] (Containerx) -- (Containerxz) -- (Containerz) -- (Containeryz) -- (Containery) -- (Containerxy) -- cycle;

\draw[solid, thick, color=black] (Container) -- (Containerxy) -- (Containerx) -- (Containerxz) -- cycle;
\draw[solid, thick, color=black] (Container) -- (Containerxy) -- (Containery) -- (Containeryz) -- cycle;
\draw[solid, thick, color=black] (Container) -- (Containerxz) -- (Containerz) -- (Containeryz) -- cycle;
\end{tikzpicture}}
    \caption{Augmentation of placement points by extreme points and normal patterns. Normal patterns (dashed, blue lines and crosses, cf. \cite{Beasley1985a}) and extreme points (dotted, red lines and circles, cf. \cite{Crainic2008}). Five useful normal pattern points that are not in the set of extreme points are shown. Two extreme points on the container walls are not in the dynamically generated set of normal pattern points; they have no supporting item directly below, but they can be useful if partial support is sufficient.}  \label{fig:NormalPatterns}
\end{figure}

In the constructive phase, items are sorted according to the rule proposed by \citet{Zhang2015}. They are placed one after the other on available placement points, which are implicitly used to represent residual space.
Placement points are generated by projecting three corner points of a newly placed item in the directions of the container walls; each of the three corners is projected in two directions, resulting in a total of six so-called extreme points \citep{Crainic2008}.
Because extreme points do not generate many useful placement points when only partial support is required, additional placement points are dynamically generated (cf. Figure \ref{fig:NormalPatterns}).
Initially, normal patterns are generated on the container floor \citep{Herz1972,Christofides1977,Beasley1985a}.
New items are placed with their bottom-left back corner on all available placement points for all allowed orientations. The placement and orientation are chosen to result in a feasible left-most downward position---i.e., the point with the lexicographically smallest x, y, and z coordinates---with ties broken in favor of the orientation that results in the shorter side of the item being placed in the x-direction. After an item is placed, new placement points are generated onto the top face of the placed item according to extreme point projections and normal pattern projections. If an item cannot be placed, it is skipped and no new placement points are generated.
A feasible packing has been found if all items can be placed and the algorithm aborts.


If not all items can be placed, there are two possibilities. If the number of items is below a certain threshold $T^{\mathrm{enum}}$, all possible $T^{\mathrm{enum}}!$ sequences are enumerated and solved with the extreme point heuristic (cf. \citet{Zhang2015}). 
Otherwise, a simple iterated greedy local search is executed (see Algorithm \ref{alg:Metaheuristic}). 

\begin{algorithm}
\DontPrintSemicolon
  
  \KwInput{items ${I}$, dimensions of bin ${B}$, neighborhoods ${N}$, moves ${M}$}
  \KwOutput{feasibility of packing of items ${I}$ into bin ${B}$}
  \KwData{solution $s$, best solution $s^*$, item placement sequence ${S}$}

    \tcp{Constructive phase.}
    Determine initial placement sequence $S_0$ of ${I}$ \tcp*{According to \citet{Zhang2015}.}
    $s^*$ = ExtremePointHeuristic(${I}$, ${B}$, ${S}_0$) \;
    \tcp{Iterated greedy with multiple neighborhoods defined on the placement sequence.}
    \While{\upshape abort criterion is not satisfied}
    {
    	\tcp{Variable neighborhood descent with pipe neighborhood change step.}
    	\ForEach{\upshape neighborhood $n \in {N}(s^*)$}
    	{
    		\tcp{Restricted local search with best improvement; abort if feasible.}
    		\ForEach{\upshape move $m \in {M}_n$}
    		{
    			$s$ = ExtremePointHeuristic(${I}$, ${B}$, ${S}_m$)\;
    			\lIf{\upshape $s$ is feasible}
    			{
    				\Return true
    			}
    			\ElseIf{\upshape number of successfully placed items in $s > s^*$  \textbf{or} \\ number of successfully placed items in $s \geq s^*$ and $v^{\mathrm{LIFO}}(s) \leq v^{\mathrm{LIFO}}(s^*)$}
    			{
    				$s^*$ = $s$
    			}
    		}
    	}
    	Perturb $s^*$ \tcp*{Locally restricted random swap move of placement sequence.}
    }
	\lIf{\upshape $s^*$ is feasible}
	{
		\Return true
	}
	\lElse
	{
		\Return false
	}
	
\caption{Augmented extreme point-based metaheuristic}
\label{alg:Metaheuristic}
\end{algorithm}

The intensification phase employs a variable neighborhood descent with a pipe neighborhood change procedure \citep{Hansen2017}. Swap and insertion neighborhoods are defined based on the placement sequence of items. The constructive extreme point-based heuristic is applied to evaluate a move using the fixed placement sequence determined by the move, rather than relying on a sorting rule. The size of the neighborhood is restricted based on a distance metric of the move sequence, in comparison to a strict LIFO placement sequence.
In contrast to the constructive part, evaluation is aborted when an item cannot be placed. A new incomplete solution is accepted if it either successfully packs strictly more items than the previously best-known solution or if it packs at least as many items, with ties broken in favor of a lower LIFO distance. The latter is the minimum absolute difference of the item placement sequence in reference to the start or end indices of the item's group in a strict LIFO sequence. Given items $i \in I$, an arbitrary position $k$ in the sequence, a tolerance threshold $T$, and the start and end indices $I_{g_i}^{\mathrm{start}}$ and $I_{g_i}^{\mathrm{end}}$ of the item's group $g_i$ in strict LIFO sequence, we define the LIFO distance as
\begin{align}
    d^{\mathrm{LIFO}} (i, k, T) =
    \begin{cases}
        I_{g_i}^{\mathrm{start}} - k, & \quad \text{if } k < I_{g_i}^{\mathrm{start}} - T, \\
        k - I_{g_i}^{\mathrm{end}}, & \quad \text{if } k > I_{g_i}^{\mathrm{end}} + T, \\
        0, & \quad \text{else.}
    \end{cases}
\end{align}
With an item's sequence index $s_i$ and the LIFO distance tolerance threshold $T^{\mathrm{LIFO, tol}}$, the total LIFO distance is calculated as
\begin{align}
    v^{\mathrm{LIFO}} = \sum_{i \in \mathcal{I}} d^{\mathrm{LIFO}} (i, s_i, T^{\mathrm{LIFO, tol}}).
\end{align}

To reduce the complexity of the local search, the neighborhoods are restricted according to a move evaluation distance such that moves are only evaluated if the LIFO distance $d^{\mathrm{LIFO}} (i, s_j, 0)$ of two items $i$ and $j$ is within a threshold $T^{\mathrm{LIFO, eval}}$---i.e., if $d^{\mathrm{LIFO}} (i, s_j, 0) \leq T^{\mathrm{LIFO, eval}}$. 
The diversification mechanism perturbs the current solution through locally restricted random swap moves: a random element $i$ is chosen and swapped with another random element $j$ in its vicinity such that the distance of the indices does not exceed the perturbation swap distance $T^{\mathrm{pert}}$---i.e., $d^{\mathrm{LIFO}} (i, s_j, 0) \leq T^{\mathrm{pert}}$. This limits diversification strength and is favorable for LIFO variants.
The metaheuristic continues until a feasible packing has been determined or a specified number of perturbations has been reached.

\section{Model strengthening}\label{sec:VI}
We apply three different model-strengthening methods to improve the computational performance. First, the two-index vehicle flow formulation presented in Section~\ref{sec:VRPModel} is strengthened by valid inequalities (VI) added during model generation. Second, cutting planes that are valid for the CVRP are added during the branch-and-cut algorithm at fractional solutions. Third, lifting methods are presented to strengthen infeasible path constraints \eqref{eq:MIP:InfPath} found at integer solutions.


\subsection{Valid inequalities} \label{sec:valid_inequalities}
The lower bound $K^{\min}$ for the number of required vehicles is determined by solving the one-dimensional bin packing problem using the items' volume and weight as resources. Constraints~\eqref{eq:VI:LBvehicles} enforce that at least $K^{\min}$ outgoing arcs from the depot are assigned.
\begin{equation}
    \sum_{\left( i, j \right) \in \delta^+\left( 0 \right)} x_{ij} \geq K^{\min}.
    \label{eq:VI:LBvehicles}
\end{equation}
In addition, inequalities~\eqref{eq:VI:Subtour} are applied to eliminate subtours with two customers.
\begin{equation}
    x_{ij} + x_{ji} \leq 1 ~ \forall~ i,j \in C, i < j.
    \label{eq:VI:Subtour}
\end{equation}

\subsection{Fractional solutions}
\label{sec:FractionalCuts}
Solving the CLP is computationally very expensive; to avoid this, we consider only the one-dimensional approximation to separate cuts for fractional solutions. In general, all cutting planes valid for the CVRP are valid for the 3L-CVRP if only the one-dimensional approximation of the 3L-CVRP is considered, using the weight and total volume as customer demand. In our branch-and-cut algorithm, we apply rounded capacity, multistar, framed capacity, and generalized large multistar inequalities \citep{Lysgaard2004}.
In addition, since the CVRP is a generalization of the asymmetric traveling salesman problem (ATSP), cutting planes valid for the ATSP are valid for the 3L-CVRP. We apply $D_k^+$ and $D_k^-$ cuts \citep{Grotschel1985} as well as CAT-cuts \citep{Balas1989}.
\subsection{Integer solutions}
All types of infeasible path constraints described in Section~\ref{sec:InfeasiblePathDef} are quite weak regardless of the path definition because they only cut off a small part of the solution space.
The following section details different options for lifting constraints \eqref{eq:MIP:InfPath}.

\subsubsection{Two-path inequalities}
All infeasible path inequalities can be lifted if there is no feasible path connecting all customer nodes in ${S}={C}(P)$, with ${C}(P)$ denoting the set of customers visited in path $P$. Thus, inequalities \eqref{eq:VI:TwoPathConstraint} can be applied to ensure at least two routes are used to serve all customers in $S$:
\begin{equation}
    \sum_{\left( i,j \right) \in {A}\left( {S} \right)} x_{ij} \leq |{S}| - 2 \label{eq:VI:TwoPathConstraint}.
\end{equation}
\cite{Kohl1999} introduced these as 2-path (2P) inequalities for the vehicle routing problem with time windows. \cite{Zhang2022a} have similarly proposed a similar type of inequality as infeasible set inequalities for the 2L-CVRP.
It is obvious that 2P $\succ$ IRP because they prohibit not just one but all paths in set $S$. To identify a violated 2P inequality, we must check whether there is no route with a feasible loading that visits all customers. However, it is crucial not to exclude any feasible path that connects a superset of $S$. Hence, the support constraint must be omitted, due to the \textit{incremental feasibility} property. We perform the feasibility check using a relaxed loading variant that only considers (i) the \lifoNoSeq{} variant (if the original loading variant must respect the LIFO policy) or (ii) the \loadingOnly{} variant (if the original variant does not include the LIFO policy). 

It should be noted that 2P inequalities can be further strengthened if the actual number of routes needed to serve set ${S}$ is determined by a three-dimensional bin-packing problem with additional loading constraints instead of by the CLP feasibility problem. However, this is computationally very expensive, and therefore it is not implemented.

\subsubsection{Lifting of infeasible regular path inequalities}
An IRP inequality can be lifted by a tournament inequality, originally introduced for the ATSP with time windows by \cite{Ascheuer2000}. To apply regular tournament (RT) inequalities to the 3L-CVRP, we define the transitive closure of $P=P^\mathrm{r}$ as $\left[ P \right] = \left\{ \left(v_i, v_j \right) \in A~|~ 1 \leq i < j \leq p \right\} $. Then, an RT inequality is defined as in \eqref{eq:VI:TournamentConstraints}:
\begin{equation}
    \sum_{\left( i,j \right) \in \left[ P \right]} x_{ij} \leq |A(P)| - 1.  
    \label{eq:VI:TournamentConstraints}
\end{equation}
In addition, since we consider a directed graph $G$, we can check the reversed path $\bar{P} = \left\{ v_p, \dots, v_1 \right\}$ \citep[cf.][]{ropkes-2007}. 
If $\bar{P}$ is infeasible as well, we can add the following inequality:
\begin{equation}
\sum_{\left( i,j \right) \in A\left(P\right)} x_{ij} + \sum_{\left( i,j \right) \in A(\bar{P})} x_{ij} \leq |A(P)| - 1. \label{eq:VI:UndirectedInfeasiblePathConstraint}
\end{equation}
We call inequality \eqref{eq:VI:UndirectedInfeasiblePathConstraint} the undirected infeasible path (UIP) inequality.

\subsubsection{Lifting of infeasible tail path inequalities}
First, an ITP inequality with $P = P^0$ can be lifted to an IRP inequality and consequently to an RT inequality if path $\tilde{P} =P^0 \setminus \{0\}$ remains infeasible with relaxed support constraints, since the \textit{incremental feasibility} property does not have to be respected.

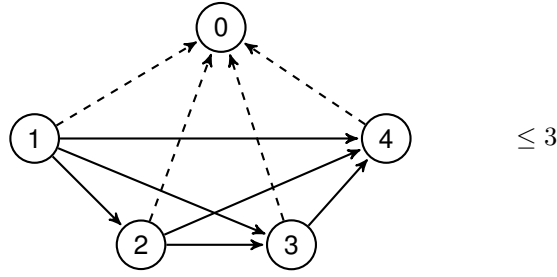
\begin{figure}
    \centering
    {\begin{tikzpicture}[->,>=stealth',shorten >=1pt,auto,node distance=2cm,
                    thick,main node/.style={circle,draw,font=\sffamily}]

  \node[main node] (1) {1};
  \node[main node] (0) [above right = 1cm and 2cm of 1]{0};
  \node[main node] (2) [below right of=1] {2};
  \node[main node] (3) [right of=2] {3};
  \node[main node] (4) [right = 4cm of 1] {4};
  \node (LE) [right of = 4] {$\leq 3$};
  
    \draw (1) -- (2);
    \draw (1) -- (3);
    \draw (1) -- (4);
    \draw (2) -- (3);
    \draw (2) -- (4);
    \draw (3) -- (4);
    \draw [dashed] (1) -- (0);
    \draw [dashed] (2) -- (0);
    \draw [dashed] (3) -- (0);
    \draw [dashed] (4) -- (0);
    
\end{tikzpicture}}
    \caption{Tournament tail inequalities for infeasible path $P = \left\{ 1,2,3,4,0\right\}$ and depot node 0. Solid lines have a multiplier of 1, while dashed lines have a multiplier of 0.5.} \label{fig:TailTournamentConstraint}
\end{figure}

If $\tilde{P}$ is not infeasible with relaxed support constraints, another lifting procedure for $P$ related to RT inequalities can be applied. The original variant is invalid because the depot---the last node in $P$---can have more than one incoming arc. If, for example, path $P = \left\{ 1,2,3,4,0\right\}$ is infeasible, a corresponding RT inequality would also exclude the combination of paths $P_1 = \left\{ 1,2,3,0\right\}$ and $P_2 = \left\{4,0\right\}$. To achieve a good lifting through coefficients as large as possible and to ensure no feasible paths are excluded, all arcs entering the depot are multiplied by a factor of 0.5. We call this the tail tournament (TT) inequality; it is defined in \eqref{eq:VI:TournamentTailConstraints}, and Figure~\ref{fig:TailTournamentConstraint} offers an example.
\begin{equation}
    \sum_{\left( i,j \right) \in \left[ \tilde{P} \right]} x_{ij} + 0.5 \sum_{i \in \tilde{P}} x_{i0} \leq |A(P)| - 1. 
    \label{eq:VI:TournamentTailConstraints}
\end{equation}

Additionally, similar to the UIP inequality, the reversed path of $P$, $\bar{P} = \left\{v_{p-1},\dots, 0 \right\}$ can also be checked in the case of an infeasible path $P$. The undirected infeasible tail path (UITP) inequality, defined in~\eqref{eq:VI:UndirectedInfeasibleTailPathConstraints}, can be applied if $P$ is infeasible in both directions (see Figure~\ref{fig:UndirectedITPC} for an example).
\begin{equation}
    \sum_{\left( i,j \right) \in A\left( \tilde{P} \right)} \left( x_{ij} + x_{ji} \right) + 0.5 \sum_{i \in \tilde{P}} x_{i0} \leq |A(P)| - 1.  
    \label{eq:VI:UndirectedInfeasibleTailPathConstraints}
\end{equation}
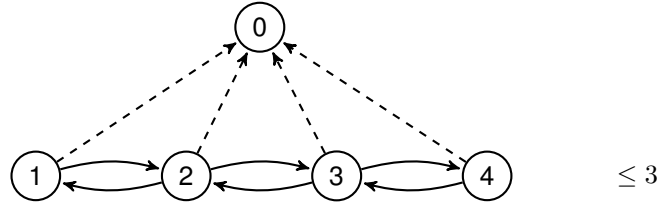
\begin{figure}
    \centering  
    {\begin{tikzpicture}[->,>=stealth',shorten >=1pt,auto,node distance=2cm,
                    thick,main node/.style={circle,draw,font=\sffamily}]

  \node[main node] (1) {1};
  \node[main node] (0) [above right = 1.5cm and 2.5cm of 1]{0};
  \node[main node] (2) [right of=1] {2};
  \node[main node] (3) [right of=2] {3};
  \node[main node] (4) [right of=3] {4};
  \node (LE) [right of = 4] {$\leq 3$};
  
    \draw [->] (1) to [bend left=15] (2);
    \draw [->] (2) to [bend left=15] (3);
    \draw [->] (3) to [bend left=15] (4);
    \draw [->] (4) to [bend left=15] (3);
    \draw [->] (3) to [bend left=15] (2);
    \draw [->] (2) to [bend left=15] (1);
    \draw [dashed] (1) -- (0);
    \draw [dashed] (2) -- (0);
    \draw [dashed] (3) -- (0);
    \draw [dashed] (4) -- (0);
    
\end{tikzpicture}}
    \caption{Undirected infeasible tail path inequalities for infeasible path $P = \left\{ 1,2,3,4,0\right\}$ and depot node 0. Solid lines have a multiplier of 1, while dashed lines have a multiplier of 0.5.} \label{fig:UndirectedITPC}
\end{figure}

\subsubsection{Lifting of infeasible route inequalities}
\begin{figure}
    \centering  
    {\begin{tikzpicture}[->,>=stealth',shorten >=1pt,auto,node distance=2cm,
                    thick,main node/.style={circle,draw,font=\sffamily}]

  \node[main node] (1) {1};
  \node[main node] (2) [below of=1] {2};
  \node[main node] (3) [below of=2] {3};
  \node[main node] (4) [right = 4cm of 1] {4};
  \node[main node] (5) [below of=4] {5};
  \node[main node] (6) [below of=5] {6};
  \node (LE) [right of = 5] {$\leq 1$};
  
    \draw [->] (1) to [bend left=15] (2);
    \draw [->] (2) to [bend left=15] (1);
    \draw [->] (2) to [bend left=15] (3);
    \draw [->] (3) to [bend left=15] (2);
    \draw [->] (1) to [bend left=30] (3);
    \draw [->] (3) to [bend left=30] (1);
    \draw [->, dashed] (1) to [bend left=15] (4);
    \draw [->, dashed] (4) to [bend left=15] (1);
    \draw [->, dashed] (1) to [bend left=15] (5);
    \draw [->, dashed] (5) to [bend left=15] (1);
    \draw [->, dashed] (1) to [bend left=15] (6);
    \draw [->, dashed] (6) to [bend left=15] (1);
    \draw [->, dashed] (2) to [bend left=15] (4);
    \draw [->, dashed] (4) to [bend left=15] (2);
    \draw [->, dashed] (2) to [bend left=15] (5);
    \draw [->, dashed] (5) to [bend left=15] (2);
    \draw [->, dashed] (2) to [bend left=15] (6);
    \draw [->, dashed] (6) to [bend left=15] (2);
    \draw [->, dashed] (3) to [bend left=15] (4);
    \draw [->, dashed] (4) to [bend left=15] (3);
    \draw [->, dashed] (3) to [bend left=15] (5);
    \draw [->, dashed] (5) to [bend left=15] (3);
    \draw [->, dashed] (3) to [bend left=15] (6);
    \draw [->, dashed] (6) to [bend left=15] (3);
    
\end{tikzpicture}}
    \caption{Two-path tail inequalities for sets $S=\{1,2,3\}$ and $C\setminus S=\{4,5,6\}$. Solid lines have a multiplier of 1, while dashed lines have a multiplier of -1.}    \label{fig:2PTC}
\end{figure}
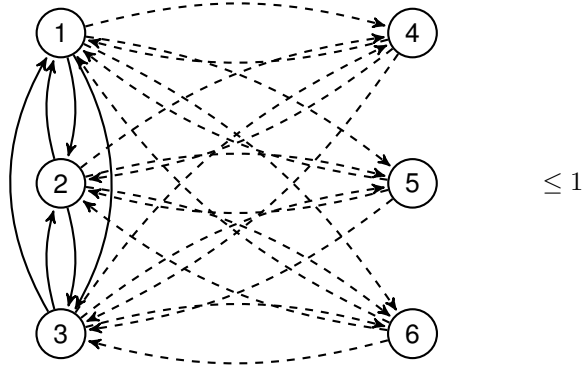

IR inequalities are only present in the \noLifo{} variant; here, the customer sequence is irrelevant to the CLP, as the LIFO policy is not considered. As a result, not just the currently checked route $R$ but all routes between customers in $S=C(R)$ can be prohibited if the load check is infeasible, similar to the 2P inequalities. However, in contrast to the 2P inequalities, the support constraint must be respected in the \noLifo{} variant, and thus the \textit{incremental feasibility} property applies---making the 2P inequalities invalid, as they forbid any path in $S$ to be a partial sequence of any other path with additional customers $C\setminus S$. We introduce the 2-path tail (2PT) inequalities in \eqref{eq:VI:TwoPathTailConstraints} to exclude all routes that consist exclusively of the customers in $S$ (see Figure~\ref{fig:2PTC} for an example).
\begin{equation}
    \sum_{(i, j) \in A(S)} x_{ij} - \sum_{i \in S, j \in C \setminus S} (x_{ij}+x_{ji}) \leq |S| - 2. 
    \label{eq:VI:TwoPathTailConstraints}
\end{equation}

\subsection{Lifting hierarchy for different loading variants}
\begin{figure}
    \centering  
    {\begin{tikzpicture}[->,>=stealth',shorten >=1pt,auto,node distance=2cm and 10cm,
                    thick,main node/.style={rectangle,draw,font=\sffamily}]                   

  \tikzset{
  default/.style={
    draw,
    ->,
  },
  noLifo/.style={
    draw,
    dashdotdotted
  },
  noSupport/.style={
    draw,
    dotted
  },
  noFragility/.style={
    draw,
  },
  loadingOnly/.style={
    draw,
    dashdotted
  },
}
  \node[main node, fill=white] (IRPC) {IRP}; 
  \node[main node] (RTC) [left of=IRPC]{RT};
  \node[main node] (UIPC) [below left=1cm and 1cm of IRPC] {UIP};
  \node[main node] (2PC) [left of=RTC] {2P};
  \node[main node] (TTC) [right= 3cm of IRPC]{TT}; 
  \node[main node, fill=white] (ITPC) [right of=TTC] {ITP};  
  \node[main node] (UITPC) [below left=1cm and 1cm of ITPC]{UITP};
  
  \node[main node] (2PC_loadingOnly) [below= 2.5cm of 2PC] {2P};
  \node[main node, fill=white] (RTC_loadingOnly) [right of=2PC_loadingOnly]{RT};  
  \node[main node, fill=white] (IRPC_loadingOnly) [right of=RTC_loadingOnly] {IRP};

  \node[main node, fill=white] (IRC) [below= 2.5cm of ITPC] {IR};
  \node[main node] (2PTC) [left of=IRC] {2PT};
  \node[main node] (2PC_NoLIFO) [left of=2PTC] {2P}; 

  \draw [default] (ITPC) to (TTC);
  \draw [default,dashed] (ITPC) to node[below right, align=center,font=\fontsize{8pt}{10pt}\selectfont] {rev. path \\ infeasible} (UITPC);
  \draw [default,dashed] (TTC) to node[below, align=center,font=\fontsize{8pt}{10pt}\selectfont] {support constraints \\ disabled} (IRPC);
  \draw [default] (IRPC) to (RTC);
  \draw [default,dashed] (IRPC) to node[below right, align=center,font=\fontsize{8pt}{10pt}\selectfont] {rev. path \\ infeasible} (UIPC);
  \draw [default,dashed] (RTC) to node[below, align=center,font=\fontsize{8pt}{10pt}\selectfont] {no fixed\\LIFO sequence} (2PC);

  \draw [default] (IRPC_loadingOnly) to (RTC_loadingOnly);
  \draw [default] (RTC_loadingOnly) to (2PC_loadingOnly);
  
  \draw [default] (IRC) to (2PTC);
  \draw [default,dashed] (2PTC) to node[below, align=center,font=\fontsize{8pt}{10pt}\selectfont] {no fixed\\LIFO sequence} (2PC_NoLIFO); 

  \node[label={[font=\fontsize{8pt}{10pt}\selectfont]left:\textit{all constraints, no fragility}}] at ($(ITPC.east |- UITPC) + (0.25,-0.6)$) {};
  \node[label={[font=\fontsize{8pt}{10pt}\selectfont]left:\textit{no support}}] at ($(IRPC.east |- UIPC) + (0.25,-0.6)$) {};
  \node[label={[font=\fontsize{8pt}{10pt}\selectfont]left:\textit{no LIFO}}] at ($(IRC.east) + (0.25,-0.6)$) {};
  \node[label={[font=\fontsize{8pt}{10pt}\selectfont]left:\textit{loading only}}] at ($(IRPC_loadingOnly.east) + (0.25,-0.6)$) {};
  
    \begin{scope}[on background layer]
        \fill [draw=black, dashed, rounded corners, fill=gray!40, opacity=0.5] ($(ITPC.east |- ITPC.north) + (0.2,0.4)$) -- ($(2PC.west |-  2PC.north) + (-0.4,0.4)$) -- ($(2PC.west |-  UITPC.north) + (-0.4,0.2)$) -- ($(UITPC.west |- UITPC.north) + (-0.2,0.2)$) -- ($(UITPC.west |- UIPC.south) + (-0.2,-0.5)$) -- ($(ITPC.east |- UITPC.south) + (0.2,-0.5)$) -- cycle;
    \end{scope}

    \begin{scope}[on background layer]
        \fill [draw=black, dashed, rounded corners, fill=gray!80, opacity=0.5] ($(IRPC.east |- IRPC.north) + (0.2,0.2)$) -- ($(2PC.west |-  2PC.north) + (-0.2,0.2)$) -- ($(2PC.west |-  UIPC.north) + (-0.2,0.4)$) -- ($(UIPC.west |- UIPC.north) + (-0.2,0.4)$) -- ($(UIPC.west |- UIPC.south) + (-0.2,-0.5)$) -- ($(IRPC.east |- UIPC.south) + (0.2,-0.5)$) -- cycle;       
    \end{scope}

    \begin{scope}[on background layer]
        \fill [draw=black, dashed, rounded corners, fill=gray!40, opacity=0.5] ($(IRPC_loadingOnly.east |- IRPC_loadingOnly.north) + (0.2,0.2)$) -- ($(2PC_loadingOnly.west |-  IRPC_loadingOnly.north) + (-0.2,0.2)$) -- ($(2PC_loadingOnly.west |-  IRPC_loadingOnly.south) + (-0.2,-0.5)$) -- ($(IRPC_loadingOnly.east |-  IRPC_loadingOnly.south) + (0.2,-0.5)$) -- cycle;
    \end{scope}

    \begin{scope}[on background layer]
        \fill [draw=black, dashed, rounded corners, fill=gray!40, opacity=0.5] ($(IRC.east |- IRC.north) + (0.2,0.2)$) -- ($(2PC_NoLIFO.west |-  IRC.north) + (-0.2,0.2)$) -- ($(2PC_NoLIFO.west |-  IRC.south) + (-0.2,-0.5)$) -- ($(IRC.east |-  IRC.south) + (0.2,-0.5)$) -- cycle;
    \end{scope}
\end{tikzpicture}}
    \caption{Classification of lifted inequalities for the different loading variants. Transitive relationships concerning the strength of the inequalities: never-used inequalities are marked in white, and the dominating inequality is connected via a solid line; in contrast, lifting represented by dashed lines requires additional conditions.}  \label{fig:Lifting}
\end{figure}
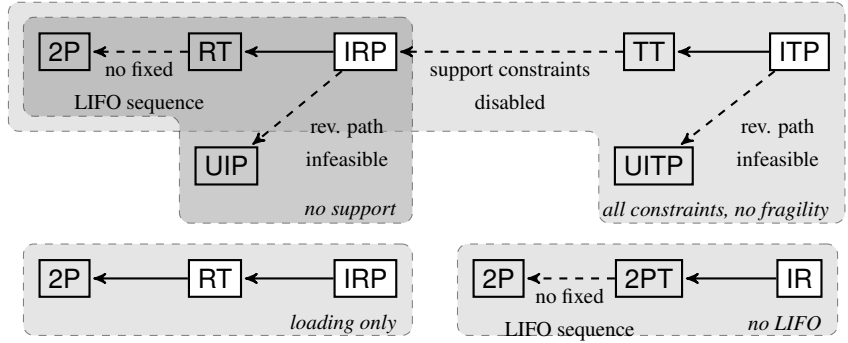

Figure~\ref{fig:Lifting} shows a hierarchy of the inequalities that are used to prohibit infeasible paths in integer solutions for each variant, starting from the infeasible path of the corresponding loading variant as described in Table~\ref{tab:InfPathConfigs}. Not all lifting procedures are necessary for all loading variants, and some inequalities are never used because they are dominated by others without any additional conditions and can be lifted immediately: RT $\succ$ IRP, TT $\succ$ ITP, 2PT $\succ$ IR. Please note that only 2P inequalities are included in the \loadingOnly{} variant, as the customer sequence is irrelevant for the subproblem.


\section{Branch-and-cut algorithm}\label{sec:BCAlgorithm}

This section presents the framework of the branch-and-cut algorithm, first describing, the general structure and then the individual components. It should be noted that the components must be adapted to the respective loading problem variant. Here, we concentrate on the most complex variant, the \allConstraints{} variant.

\subsection{General procedure}

The general structure of the branch-and-cut routine is described in Algorithm~\ref{alg:BCRoutine}. To start, we apply a preprocessing procedure to eliminate infeasible arcs (see Section~\ref{sec:Preprocessing}). Subsequently, an initial solution with at most $K$ vehicles is generated with a constructive heuristic (see Section~\ref{sec:ConstructiveHeuristic}). During the branch-and-bound search, three procedures are used. First, to ensure feasibility, newly obtained integer solutions must be checked (see Section~\ref{sec:IntegerSolutionCallback}). These feasibility checks are computationally very expensive; therefore, we have integrated a memory management to avoid redundant feasibility checks and to speed up the algorithm. Feasible routes are stored in set $\mathcal{R}$, while infeasible customer combinations are stored in set $\mathcal{C}^{\text{Inf}}$. In both sets, the corresponding elements of the respective loading variants and relaxations are stored separately. Infeasible customer combinations are only used for checks without an active support constraint. It should be noted that infeasible paths must not be stored, as they can no longer appear in new solutions once they are prohibited. Second, in addition to integer solutions, we separate and add cuts for fractional solutions to improve the lower bound (LB) (see Section~\ref{sec:FractionalSolutionCallback}). Third, we apply a set-partitioning-based heuristic using set $\mathcal{R}$ to more quickly obtain better upper bounds (UB) (see Section~\ref{sec:SPHeuristic}). The branch-and-cut algorithm's standard termination criterion is the convergence of the UB and LB, the achievement of a proven optimal solution, or the reaching of a time limit.


\begin{algorithm}
\DontPrintSemicolon
  \KwOutput{best solution $s$}
    Set of feasible routes $\mathcal{R} \gets \emptyset$, set of infeasible customer combinations $\mathcal{C}^{\text{Inf}} \gets \emptyset$\;
    Preprocessing \tcp*{Section~\ref{sec:Preprocessing}}
    $s \gets$ ConstructiveHeuristic() \tcp*{Section~\ref{sec:ConstructiveHeuristic}}
    \While {termination criteria not met}
    {
        $n \gets$  current node in branch-and-bound tree\;
        \eIf{$n$ has new integer solution $s'$ with $f(s') < f(s)$}
        {
            $solutionFeasible \gets$ IntegerSolutionCallback($s', \mathcal{R}, \mathcal{C}^{\text{Inf}}$) \tcp*{Section~\ref{sec:IntegerSolutionCallback}}
            \If{solutionFeasible}
            {
                $s \gets s'$ 
            }
        }
        {
            \If{$n$ has new fractional solution $\tilde{s}$ and call condition is satisfied}
            {
                FractionalSolutionCallback($\tilde{s}$) \tcp*{Section~\ref{sec:FractionalSolutionCallback}}
            }
            \If{Call condition for upper bound procedure is satisfied}
            {   
                $s' \gets$  SPHeuristic$\left(\mathcal{R}, s\right)$ \tcp*{Section~\ref{sec:SPHeuristic}}
                \lIf{$s' \neq \emptyset$}{$s \gets s'$}
            }
        }
    }
    \Return $s$
\caption{General branch-and-cut algorithm}
\label{alg:BCRoutine}
\end{algorithm}


\subsection{Preprocessing}\label{sec:Preprocessing}
To speed up the solution process, infeasible arcs $(i,j)$ with $i,j \in C$---i.e., paths with two customer nodes---are identified during preprocessing. The identification has two steps. In the first step, the loading check of the respective variant with relaxed support constraints is performed; if the check is infeasible, arc $(i,j)$ can be excluded (i.e., $x_{ij}$ is fixed to zero).
The second step is only relevant for loading variants that include the support constraint and when $(i,j)$ is feasible without active support constraint. We check arc $(i,j)$ with the complete constraint set of the respective variant. If this is feasible, we add route $R = \left\{0,i,j,0\right\}$ to set $\mathcal{R}$; otherwise, we try to lift each infeasible path $P=\{i,j,0\}$ by extending arc $(i,j)$ to $P^*=\{i,j,k\}$ with $k \in C \setminus \{i,j\}$ and evaluate $P^*$ with relaxed support constraints. If all extended paths $P^*$ are infeasible, arc $(i,j)$ is infeasible and can be excluded; if not, we add an ITP inequality with $P=\{i,j,0\}$. Finally, if arcs $(i,j)$ and $(j,i)$ are infeasible, we add set $\left\{i,j\right\}$ to set $\mathcal{C}^{\text{Inf}}$.

\subsection{Checking integer solutions}
\label{sec:IntegerSolutionCallback}
Algorithm~\ref{alg:CallbackIntegerSol} describes the callback for checking the feasibility of all routes $R \in \mathcal{R}_{s'}$ in a new integer solution $s'$. Each route $R \in \mathcal{R}_{s'}$ is validated with respect to connectivity to the depot and exceeding the volume or weight capacity. If an infeasibility is detected, we add a GSEC for set $S=C(R)$ with the LB $r(S)$; otherwise, a check follows as to whether the route has already been stored in $\mathcal{R}$ or $S$ as a superset of an infeasible customer combination ${T} \in \mathcal{C}^{\text{Inf}}$. If the preliminary checks do not yield a positive result, the augmented extreme point-based packing heuristic is used to quickly check feasibility. Finally, if necessary, the CP model is called to definitively determine feasibility or infeasibility. This step is explained in more detail below.

\begin{algorithm}[h]
\DontPrintSemicolon
  \KwInput{solution $s'$, set of feasible routes $\mathcal{R}$, set of infeasible customer combinations $\mathcal{C}^{\text{Inf}}$}
    Determine set of routes $\mathcal{R}_{s'}$ in $s'$\;
    $cutAdded \gets$ false\;
    \ForEach{\upshape Route $R \in \mathcal{R}_{s'}$ }
    {
        $S \gets C(R)$\;
        $r\left(S\right) \gets \max\left( \left \lceil \frac{q\left( S \right)}{Q} \right \rceil, \left \lceil \frac{v\left( S \right)}{V} \right \rceil \right)$\;
        \If{\upshape $R$ disconnected from depot} 
        {
            Add GSEC($S, r\left(S\right)$) \tcp*{$R$ is a subtour without the depot.}
            $cutAdded \gets$ true
        }
        \Else
        {
            \If{\upshape $r\left(S\right) > 1$} 
            {
                Add GSEC($S, r\left(S\right)$) \tcp*{Infeasible volume and/or weight capacity.}
                $cutAdded \gets$ true
            }
            \Else
            {
                \lIf{\upshape $R \in \mathcal{R}$}
                {
                    continue
                }
                \If{$\exists {T} \in \mathcal{C}^{\text{Inf}}: S \supseteq {T}$}
                {
                    Add GSEC($S, 2$) \tcp*{$S$ contains an infeasible customer combination ${T}$.}
                    $cutAdded \gets$ true\;
                    continue
                }
                \lIf{\upshape $R$ is feasible according to heuristic}
                {
                    Add $R$ to $\mathcal{R}$
                }
                \Else
                {
                    $cutAdded \gets$ Check $R$ with CP model \tcp*{See flow chart in Figure \ref{fig:FlowChartExactRouteCheck}.}
                    \lIf{\upshape $cutAdded$ is false} {Add $R$ to $\mathcal{R}$}
                }
            }
        }
    }
    \If{\upshape $cutAdded$}
    {
        \Return false \tcp*{Do not accept $s'$ as new incumbent.}
    }
    
    \Return true \tcp*{Accept $s'$ as new incumbent.}
\caption{IntegerSolutionCallback}
\label{alg:CallbackIntegerSol}
\end{algorithm}

\begin{figure}[h]
    \centering  
    {\tikzstyle{start} = [rectangle, rounded corners, 
minimum width=2cm, 
minimum height=0.75cm,
font=\footnotesize,
text centered,
align=left,
draw=black, 
fill=gray!50]

\tikzstyle{stopPositive} = [rectangle, rounded corners, 
minimum width=2cm, 
minimum height=0.75cm,
font=\footnotesize,
text centered, 
draw=black, 
fill=blue!40]

\tikzstyle{stopNegative} = [rectangle, rounded corners, 
minimum width=2cm, 
minimum height=0.75cm,
font=\footnotesize\linespread{0.9}\selectfont,
text centered, 
draw=black, 
fill=red!50]

\tikzstyle{io} = [trapezium, 
trapezium stretches=true, 
trapezium left angle=70, 
trapezium right angle=110, 
minimum width=3cm, 
minimum height=1cm, text centered, 
draw=black, fill=blue!30]

\tikzstyle{process} = [rectangle,
minimum width=1cm, 
minimum height=0.5cm,
inner sep=1pt,
text centered, 
draw=black, 
fill=orange!30]

\tikzstyle{decision} = [diamond, 
font=\footnotesize\linespread{0.9}\selectfont,
inner sep=1pt,
text centered, 
draw=black,
fill=green!30]
\tikzstyle{arrow} = [thick,->,>=stealth]

\begin{tikzpicture}[node distance=2cm]

\node (Start) [start] {New route $R$ \\ ${S}={C}(R)$};

\node (CPLimited) [decision, right of=Start, align=left,xshift=0.75cm] {1. CP(R)?\\\textit{AllConstr}};
\node [process, above of=CPLimited, anchor=south west, yshift=-1.5cm] {T\textsuperscript{Limited}};
\node (AcceptLimited) [stopPositive, above of=CPLimited, yshift=0.25cm] {Accept route};

\node (TwoPC) [decision, right of=CPLimited, align=left,xshift=0.75cm] {2. CP(R)?\\\textit{LIFO\textsuperscript{NoSeq}}};
\node [process, above of=TwoPC, anchor=south west, yshift=-1.5cm] {T\textsuperscript{Limited}};
\node (TwoPCResult) [stopNegative, below of=TwoPC, yshift=-0.25cm,align=left] {Add 2P for ${S}$};
\node (RedInfSet) [start, below of=TwoPCResult, yshift=0.25cm,align=left] {Determine $\hat{S}$ \\ Add $\hat{S}$ to $\mathcal{C}^{\text{inf}}$};

\node (IPC) [decision, right of=TwoPC, align=left,xshift=0.75cm] {3. CP(R)?\\\textit{NoSupport}};
\node [process, above of=IPC, anchor=south west, yshift=-1.5cm] {T\textsuperscript{Limited}};
\node (RedInfPath) [start, below of=IPC, yshift=-0.25cm,align=left] {Determine $\hat{P}$};
\node (IPCResult) [stopNegative, below of=RedInfPath, yshift=0.25cm,align=left] {Add RT for $\hat{P}$};

\node (CPResult) [decision, right of=IPC, align=left,xshift=0.75cm] {1. CP(R)?};

\node (CP) [decision, right of=CPResult, align=left,xshift=0.75cm] {4. CP(R)?\\\textit{AllConstr}};
\node [process, above of=CP,anchor=south west, yshift=-1.5cm] {T\textsuperscript{Unlimited}};
\node (Accept) [stopPositive, above of=CP, yshift=0.25cm] {Accept route};

\node (RevPath) [decision, below of=CPResult, align=left,yshift=-2.5cm] {5. CP($\bar{R}$)?\\\textit{AllConstr}};
\node [process, above of=RevPath, anchor=south west, yshift=-1.5cm] {T\textsuperscript{Limited}};
\node (RevPathResult) [stopNegative, right of=RevPath,align=left, xshift=1cm] {Add UITP for $P$\\+ TT for $\bar{P}$};

\node (ITPC) [stopNegative, below of=CP, align=left,xshift=0.5cm, yshift=-0.5cm] {Add TT for $P$};

\draw [arrow] (Start) -- (CPLimited);

\draw [arrow] (CPLimited) -- node[font=\footnotesize,anchor=west,align=left] {feas.} (AcceptLimited);
\draw [arrow] (CPLimited) -- node[font=\footnotesize\linespread{0.9}\selectfont,below,align=left] {inf./\\ unk.} (TwoPC);

\draw [arrow] (TwoPC) -- node[font=\footnotesize,anchor=west] {inf.} (TwoPCResult);
\draw [arrow] (TwoPCResult) -- (RedInfSet);
\draw [arrow] (TwoPC) -- node[font=\footnotesize\linespread{0.9}\selectfont,anchor=north, align=left] {feas./\\ unk.} (IPC);

\draw [arrow] (IPC) -- node[font=\footnotesize,anchor=west] {inf.} (RedInfPath);
\draw [arrow] (RedInfPath) -- (IPCResult);
\draw [arrow] (IPC) -- node[font=\footnotesize\linespread{0.9}\selectfont,anchor=north,align=left] {feas./\\ unk.} (CPResult);

\draw [arrow] (CPResult) -- node[font=\footnotesize,anchor=north,align=left] {unk.} (CP);
\draw [arrow] (CP) -- node[font=\footnotesize,anchor=west,align=left] {feas.} (Accept);
\draw [arrow] (CP.south) -- node[font=\footnotesize\linespread{0.9}\selectfont,anchor=west,align=left] {inf.} ++ (0,-0.25) -- ++ (0,0) -| (RevPath.north);

\draw [arrow] (CPResult) -- node[font=\footnotesize\linespread{0.9}\selectfont,anchor=east,align=left] {inf.} (RevPath);

\draw [arrow] (RevPath) -- node[font=\footnotesize,anchor=north,pos=0.4] {inf.} (RevPathResult);
\draw [arrow] (CPResult.south) -- ++ (0,-0.5) -- ++ (0,0) |- (ITPC.west);


\begin{scope}[on background layer]
    \fill [draw=black, dashed, rounded corners, fill=gray!10, opacity=0.5] ($(AcceptLimited.west |- AcceptLimited.north) + (-0.2,0.1)$) -- ($(RevPath.east |-  AcceptLimited.north) + (0,0.1)$) -- ($(RevPath.east |- RevPath.north) + (0,0)$) -- ($(RevPathResult.east |- RevPath.north) + (0.1,0)$) -- ($(RevPathResult.east |- RevPath.south) + (0.1,-0.1)$) -- ($(AcceptLimited.west |- RevPath.south) + (-0.2,-0.1)$) -- cycle;
   -- cycle;
   \node[label={above right:Lifting}] at (CPLimited.west |- RevPath.south) {};
\end{scope}
\end{tikzpicture}}
    \caption{Flow chart of the exact feasibility check of a route with a CP model considering the \allConstraints{} variant.}  \label{fig:FlowChartExactRouteCheck}
\end{figure}
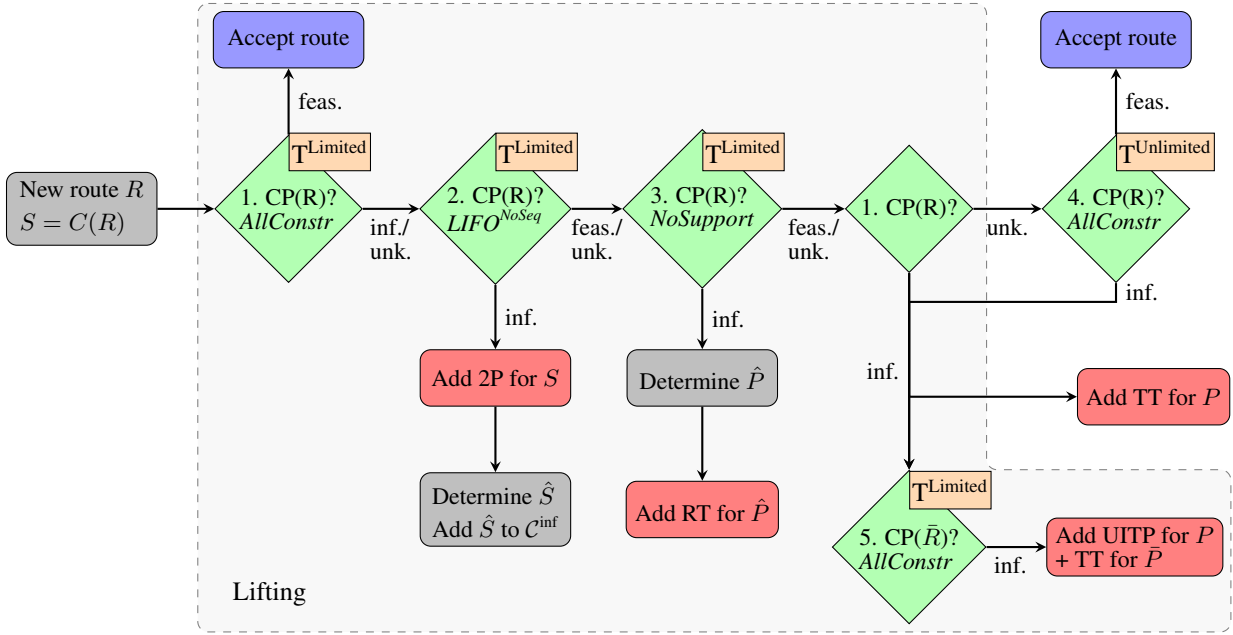

The flow chart in Figure~\ref{fig:FlowChartExactRouteCheck} summarizes the exact feasibility check for the \allConstraints{} variant. The adapted routines for the other loading variants are shown in Appendix A.

The straightforward approach to performing the feasibility check would be to simply call the CP model with the \allConstraints{} variant, specify no time limit, and either accept $R$ or add a TT inequality. However, lifted constraints can be separated beforehand. Prior to separation, we utilize the CP model with a brief time limit to promptly verify a feasible solution, thereby avoiding the need to unnecessarily solve relaxed loading variants. If the result is infeasible or unknown, the problem is relaxed in the second call of the CP model, which considers the \lifoNoSeq{} relaxation. If infeasibility can be proven within the limited run time, we add a 2P inequality with $S=C(R)$. For faster checks in subsequent integer solution callbacks, we reduce $S$ by iteratively removing the customer with the smallest volume and rechecking the remaining set with a time limit until infeasibility cannot be proven. The last reduced infeasible subset $\hat{S}$ is added to set $\mathcal{C}^{\text{Inf}}$. It should be noted that $\hat{S}$ can also be called a minimal infeasible subset if the last checked set is feasible. 

If no 2P inequality has been added, we relax the support constraint in the third call to separate an RT inequality. In the case of proven infeasibility, we explore the possibility of reducing the considered path to a reduced infeasible subpath. Customers are iteratively removed, first from the start of the path and then from the end, until the reduced path is feasible or the time limit is reached. We add an RT inequality for the last reduced infeasible subpath $\hat{P}$. Like infeasible subsets, the reduced infeasible path is minimal if the last checked path is feasible. 

If no inequality is added and if the result of the first call of the CP model is unknown, a fourth call---this time without a time limit---is made to ultimately determine the feasibility status of $R$. If $R$ is feasible, it is accepted, and the check is terminated; otherwise, we add a TT inequality with $P = P^0$. Finally, the reversed route $\bar{R}$ is checked: if $\bar{R}$ is infeasible, a UITP inequality and a TT inequality for the reversed path $\bar{P}$ are added.

\subsection{Separation of cuts in fractional solutions} \label{sec:FractionalSolutionCallback}
We use the separation algorithms introduced in \cite{Lysgaard2004} to identify violations of the CVRP inequalities mentioned in Section~\ref{sec:FractionalCuts}. For the GSECs, we use weight and volume as demand for separating new cuts, while for all other cuts we only use the weight, as this is more restrictive in the test instances. In addition, we identify violated ATSP inequalities with the separation algorithms introduced by \cite{Fischetti1997}. For the first $n^{\text{all}}$ nodes in the branch-and-bound tree, we invoke the fractional solution callback at all nodes and call the separation procedures of all inequalities. After that, we start the fractional solution callback only after every $n^{\text{sep}}$-th node, and we add the inequalities sequentially. This means that we check for violated cuts in a given order and stop when a cut is found for the current type. This strategy reduces the run time of the separation algorithms and avoids adding too many dense rows to the model, which can slow down the solution of the linear relaxation at each brand-and-bound node.

\subsection{Upper bound generation}\label{sec:UBProcedure}
Better UBs accelerate the algorithm, as they help to prune nodes. We employ a constructive heuristic to provide an initial UB and iteratively call a set-partitioning-based heuristic (SP heuristic) during the algorithm each time at least $n^\text{SP}$ new routes have been found. As the SP heuristic utilizes set $\mathcal{R}$, we try to generate more routes by evaluating each path that needs to be checked for feasibility, with the packing heuristic respecting the complete set of constraints of the problem variant instead of the currently needed relaxation. As a result, the path can be used as a new route if it is feasible. Finally, we use a simple intra-tour local search with the two-opt neighborhood to improve routes; this procedure is invoked each time a new feasible route has been found.
\subsubsection{Initial solution} \label{sec:ConstructiveHeuristic}
First, we use the modified savings algorithm introduced by \cite{Zhang2015} to generate an initial set of routes and add them to set $\mathcal{R}$. Thereafter, the local search phase tries to improve each $R \in \mathcal{R}$, and we apply the SP heuristic.
\subsubsection{Set-partitioning-based heuristic} \label{sec:SPHeuristic}
 This heuristic is based on a set partitioning (SP) formulation and tries to find new routes by combining feasible routes stored in set $\mathcal{R}$. First, we define the SP model: for each route $R \in \mathcal{R}$, the binary parameter $a_{iR}$ denotes whether customer $i \in C$ is visited on route $R$, parameter $c_R$ represents its costs, and the binary decision variable $y_R$ indicates whether a route $R$ is used ($y_R =1$) or not ($y_R = 0$). The underlying SP model is described as follows \citep{Balinski1964}:

\begin{align}
\left(SP\right) \quad \min &\sum_{R \in \mathcal{R}} c_{R} y_{R} \label{eq:SP:Objective} \\
\text{s.t.} & \sum_{R \in \mathcal{R}} a_{iR} y_{R} = 1 \quad \forall~i \in C \label{eq:SP:AllExactOnceVisited} \\
 & \sum_{R \in \mathcal{R}} y_{R} \leq K \label{eq:SP:MaxFleet} \\
 & y_R \in \left\{0, 1\right\} \quad \forall ~ R \in \mathcal{R} \label{eq:SP:Binary}
\end{align}

The objective function \eqref{eq:SP:Objective} minimizes the total routing costs. Equalities \eqref{eq:SP:AllExactOnceVisited} ensure that each customer is visited exactly once. Inequality \eqref{eq:SP:MaxFleet} limits the number of routes, and constraints \eqref{eq:SP:Binary} define the decision variables. Usually, in the context of vehicle routing problems, constraints \eqref{eq:SP:AllExactOnceVisited} can be converted to inequalities, leading to a set covering (SC) formulation: customers who are visited multiple times can be removed from a route, with the route remaining feasible. However, this is not permitted in the 3L-CVRP with active support constraints. Nevertheless, we use the SC model's linear relaxation ($LSC$) to generate new routes. As shown in Algorithm~\ref{alg:SPHeuristic}, we determine all customers who are visited multiple times in the $LSC$ and remove them from the corresponding routes. We try to improve $R'$ using the local search procedure, and we add it to $\mathcal{R}$ if it is feasible and not already included. Finally, we solve $SP$ with an updated set $\mathcal{R}$, and we return the new incumbent solution $s'$ if it improves the UB.


\begin{algorithm}
\DontPrintSemicolon
  
  \KwInput{set of feasible routes $\mathcal{R}$, incumbent solution $s$}
  \KwOutput{solution $s'$}
  \KwData{solution of $LSC$ $\bar{s}$}
    
    $\bar{s} \gets$ Solve $LSC$ with $\mathcal{R}$\;
    \lIf{\upshape $f\left( \bar{s} \right) \geq f(s)$}
    {
        \Return $\emptyset$
    }
    \ForEach{\upshape customer node $i \in C$}
    {
        Determine routes $\mathcal{R}_i \in \bar{s}$ in which $i$ is visited\;
        \If{\upshape $|\mathcal{R}_i| > 1$}
        {
            \ForEach{\upshape $R \in \mathcal{R}_i$}
            {
                $R' \gets$ $R$ with $i$ removed\;
                Improve $R'$ using two-opt\;
                \lIf{\upshape $R' \notin \mathcal{R}$ and feasible}
                {
                    Add $R'$ to $\mathcal{R}$
                }
            }
        }
        
    }
    
    $s' \gets$ Solve $SP$ with $\mathcal{R}$\;
    \lIf{\upshape $f\left( s' \right) \geq f(s)$}
	{
		\Return $\emptyset$
	}
    \Return $s'$
\caption{SPHeuristic}
\label{alg:SPHeuristic}
\end{algorithm}


\section{Computational studies}\label{sec:ComputationalStudies}


Our computational study includes the first 19 benchmark instances of \cite{Gendreau2006} for all five loading variants, resulting in a total of 95 different instances. The instances consist of up to 50 customers, each requesting one to three items, with the largest instance containing 99 items in total.

The presented framework is implemented in C++ using Gurobi 10.0 and its associated callback routines for the branch-and-cut algorithm along with the CP-SAT solver of Google OR-Tools 9.8 for the CP model. To separate violated CVRP inequalities, we apply the open-source CVRPSEP library, which implements the algorithms presented in \cite{Lysgaard2004}. In addition, we adapt the open-source code of \cite{Broek2021} for the separation algorithms of $D_k^+$, $D_k^-$, and CAT-cuts. Parameters for the fractional solution callback are set to $n^\text{all} = 200$ and $n^\text{sep} = 100$.  We limit the run time for calling the CP model in the lifting procedures to one second in each separation, except for separating 2P inequalities, where the time limit is four seconds. The SP heuristic is called each time at least $n^{\text{SP}}=20$ new routes have been added.

The study is performed on an AMD EPYC 7513 with 2.6 GHz clock speed and 32 GB RAM. The maximum number of threads is set to eight, and an overall time limit of eight hours (28,800 seconds) is defined. Run times throughout the experiments are reported in seconds unless otherwise specified. All solutions, a simple web app to visualize solutions and solver statistics as well as the source code of the described branch-and-cut algorithm are publicly available in \cite{tamkef-2024}.


Our presentation of the results is structured into two parts. The first two sections consider exclusively the \allConstraints{} variant, while the last two sections consider all five loading variants. In Section~\ref{sec:Effectiveness}, we evaluate the effectiveness of extensions added to the branch-and-cut algorithm. Subsequently, in Section~\ref{sec:HardInstances}, we try to pinpoint challenges in solving the instances. In Section~\ref{sec:ComparisonBKS}, we compare the results of our branch-and-cut algorithm with the current best-known heuristic solutions from the literature. Finally, we conduct a comparative analysis of the results obtained from the different loading variants in Section~\ref{sec:LoadingVariantComparison}.

\subsection{Evaluation of the effectiveness of enhancements to the branch-and-cut-algorithm} \label{sec:Effectiveness}
We examine two configurations of the algorithm. The first configuration, denoted as \basicBC{}, serves as a benchmark and represents a minimal basic configuration. Essentially, it comprises the algorithm using the two-index vehicle flow formulation combined with ITP inequalities; thus, route feasibility checks are performed using exclusively the CP model, and memory management, model strengthening, and the UB procedures are omitted. The second configuration, denoted as \completeBC{}, incorporates all introduced enhancements. The valid inequalities presented in \ref{sec:valid_inequalities} and the preprocessing procedure described in Section~\ref{sec:Preprocessing} are used in both configurations so that the arc set $A$ is always identical. However, all checks in the preprocessing procedure in configuration \basicBC{} are only performed with the CP model. The results are shown in Table~\ref{tab:ComparisonSimpleVsCompleteBC}. It should be noted that the total run time $t$ includes the run times for preprocessing, the start solution procedure, and the actual branch-and-cut part. However, the latter has by far the largest share.

\begin{table}[htpb]
    \centering
    \begin{threeparttable}
    \caption{Comparison of \basicBC{} and \completeBC{}} \label{tab:ComparisonSimpleVsCompleteBC} 
    {\begin{tabular}{l@{\hskip 8pt}r@{\hskip 8pt}r@{\hskip 8pt}r@{\hskip 8pt}r@{\hskip 8pt}rr@{\hskip 8pt}r@{\hskip 8pt}r@{\hskip 8pt}r@{\hskip 8pt}r@{\hskip 8pt}r@{\hskip 8pt}r}
    \toprule
    & \multicolumn{5}{c}{B\&C\textsuperscript{B}} & \multicolumn{7}{c}{B\&C\textsuperscript{C}}  \\
    \cmidrule(rl{3pt}){2-6} \cmidrule(rl{3pt}){7-13} \\[-6pt]
    Instance &
    Obj & Gap & $t$ & $N$ & $t_{\mathrm{Int}}$  &
    Obj & Gap & $t$ & $N$ & $t_{\mathrm{Int}}$ & $t_{\mathrm{Frac}}$ & $t_{\mathrm{SP}}$\\
    \midrule
    E016-03m & 301.66 & 0.0 & 57 & 55\,750 & 38 & 301.66 & 0.0 & \textbf{20} & 6\,261 & 15 & 1 & 1 \\
E016-05m & 334.96 & 0.0 & 9 & 4\,627 & 2 & 334.96 & 0.0 & \textbf{2} & 294 & 0 & 2 & 0 \\
E021-04m & 385.53 & 0.0 & 2\,042 & 790\,777 & 255 & 385.53 & 0.0 & \textbf{257} & 83\,140 & 195 & 5 & 2 \\
E021-06m & 430.88 & 0.0 & 31 & 20\,597 & 16 & 430.88 & 0.0 & \textbf{6} & 332 & 1 & 4 & 0 \\
E022-04g & 427.56 & 0.0 & 2\,903 & 1\,114\,326 & 725 & 427.56 & 0.0 & \textbf{680} & 123\,947 & 581 & 6 & 8 \\
E022-06m & 498.16 & 0.0 & 54 & 50\,581 & 28 & 498.16 & 0.0 & \textbf{15} & 2\,829 & 9 & 4 & 0 \\
E023-03g & 757.88 & 5.3 & TL & 2\,331\,349 & 1\,492 & 757.88 & \textbf{0.0} & 1\,684 & 762\,444 & 939 & 8 & 2 \\
E023-05s & 800.98 & 11.8 & TL & 1\,304\,217 & 2\,637 & 798.65 & \textbf{0.0} & 18\,151 & 2\,948\,306 & 900 & 37 & 1 \\
E026-08m & 630.13 & 0.0 & 455 & 534\,927 & 31 & 630.13 & 0.0 & \textbf{36} & 13\,316 & 9 & 9 & 0 \\
E030-03g & - & - & TL & 1\,145\,008 & 8\,835 & 769.32 & \textbf{16.0} & TL & 1\,559\,599 & 8\,529 & 63 & 148 \\
E030-04s & - & - & TL & 792\,133 & 6\,585 & 728.32 & \textbf{15.2} & TL & 2\,239\,712 & 2\,346 & 76 & 7 \\
E031-09h & 610.23 & 4.4 & TL & 4\,095\,469 & 169 & 610.23 & \textbf{0.0} & 1\,071 & 966\,306 & 63 & 39 & 1 \\
E033-03n & - & - & TL & 645\,616 & 6\,579 & 2\,619.78 & \textbf{12.4} & TL & 1\,121\,856 & 17\,003 & 60 & 67 \\
E033-04g & - & - & TL & 616\,005 & 4\,583 & 1\,320.84 & \textbf{23.0} & TL & 794\,572 & 15\,655 & 185 & 202 \\
E033-05s & - & - & TL & 514\,153 & 5\,112 & 1\,265.73 & \textbf{20.9} & TL & 778\,744 & 16\,819 & 141 & 1\,071 \\
E036-11h & 703.35 & 6.2 & TL & 2\,779\,378 & 59 & 698.61 & \textbf{0.0} & 336 & 244\,155 & 2 & 32 & 0 \\
E041-14h & 872.79 & 8.5 & TL & 1\,627\,793 & 150 & 866.40 & \textbf{3.2} & TL & 2\,655\,765 & 16 & 558 & 0 \\
E045-04f & - & - & TL & 646\,364 & 2\,894 & 1\,204.73 & \textbf{17.9} & TL & 363\,266 & 11\,437 & 313 & 128 \\
E051-05e & - & - & TL & 646\,275 & 4\,597 & 721.12 & \textbf{18.0} & TL & 567\,682 & 10\,765 & 574 & 89 \\

    \bottomrule
\end{tabular}}
    \begin{tablenotes}[para,flushleft]
        \item[] Solution quality is compared by obj, gap (\%), and solution process metrics, including the total run time ($t$; TL = time limit) and explored branch-and-bound nodes ($N$). Additionally, the breakdown of run time specifies route validation ($t_{\mathrm{Int}}$), fractional cut separation ($t_{\mathrm{Frac}}$), and SP heuristic times ($t_{\mathrm{SP}}$). Bold indicates the minimum gap per instance. If both variants prove optimality, bold indicates minimum run time.
    \end{tablenotes}
    \end{threeparttable}
\end{table}

Both configurations can optimally solve all instances with up to 22 nodes plus instance \texttt{E026-08m} and achieve the same optimal values. However, the significant differences in run times and in the number of explored branch-and-bound nodes highlight the advantages of the various extensions for these smaller instances: specifically,
due to the enhancements, \completeBC{} can solve four additional instances to optimality. For larger instances, while it faces challenges in proving optimality, \basicBC{} struggles to find feasible solutions for most instances with more than 26 nodes, underscoring the importance of the UB procedures. In general, only a short time is spent in the SP heuristic and in separating inequalities in fractional solutions; in contrast, most of the run time is devoted to the integer solution callback and to exploring the branch-and-bound tree.

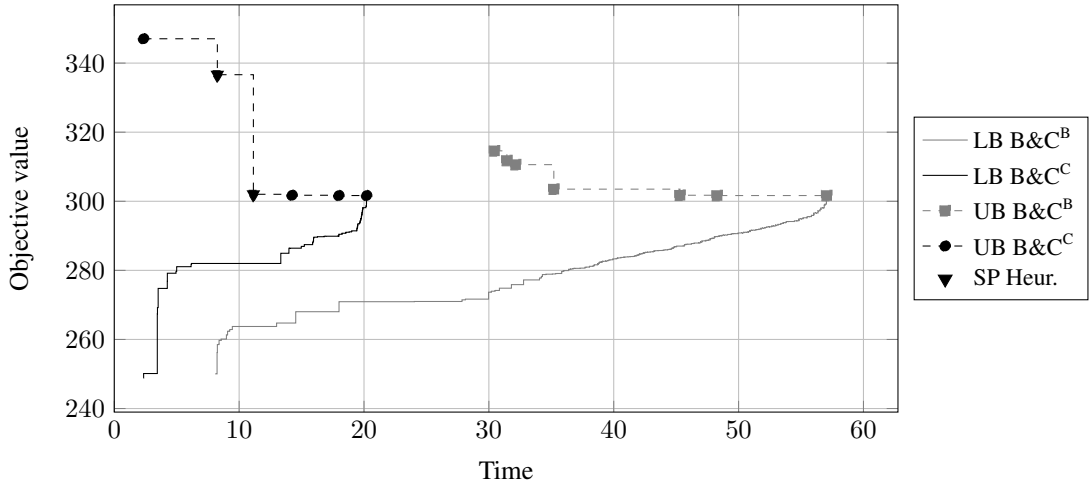
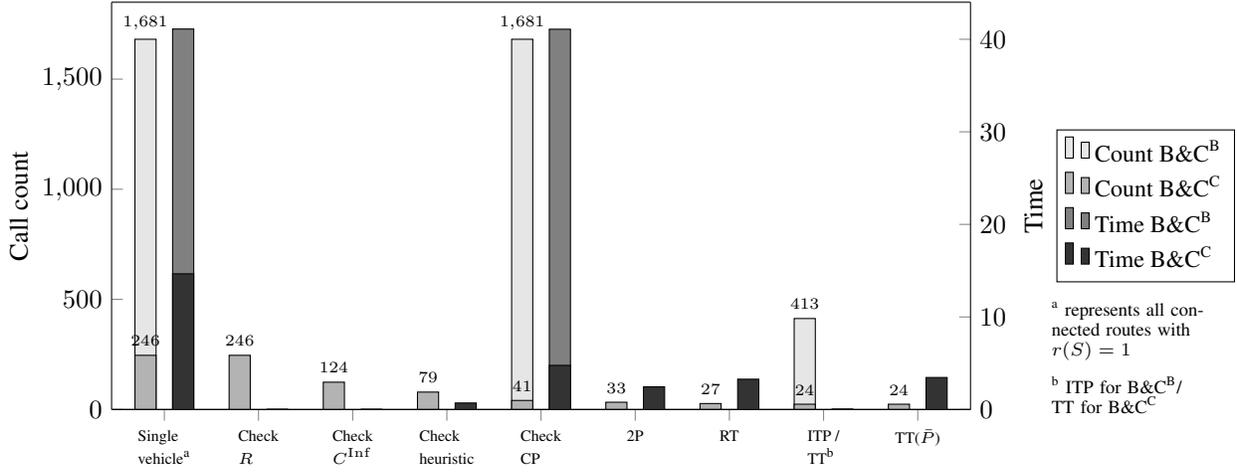
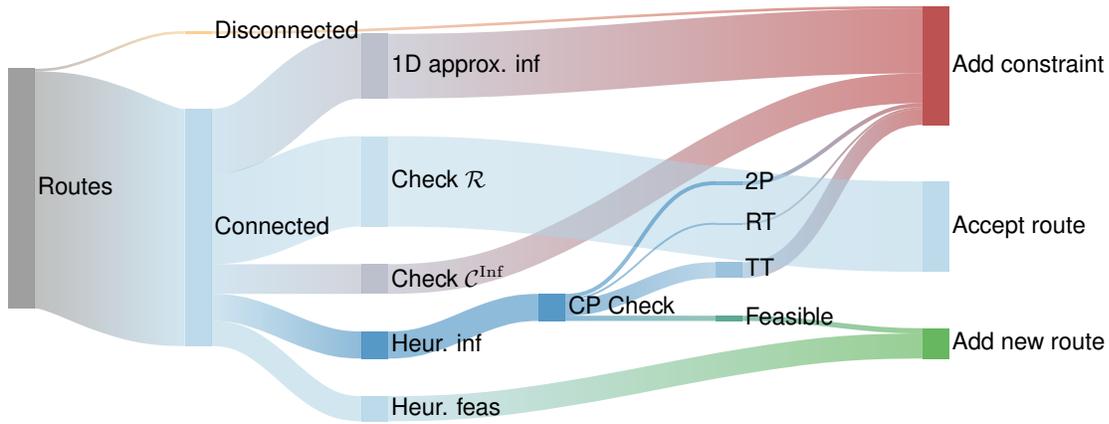
\begin{figure}
    \centering  
    \begin{subfigure}{\textwidth}        
        \begin{tikzpicture}
    \begin{axis}[
        xlabel={Time},
        height = 7cm,
        width = 12cm,
        xmin=0,
        ylabel={Objective value},
        legend style={at={(1.02,0.5)},anchor=west},
        legend style={font=\footnotesize,align=left},
        legend cell align={left},
        legend entries = {LB B\&C\textsuperscript{B}, LB B\&C\textsuperscript{C}, UB B\&C\textsuperscript{B}, UB B\&C\textsuperscript{C}, SP Heur.},
        grid=both
      ]
    \addlegendimage{no markers,black!50}
    \addlegendimage{no markers,black}
    \addlegendimage{mark=square*, black!50, dashed}
    \addlegendimage{mark=*,black, dashed}
    \addlegendimage{only marks, mark=triangle*, mark options={rotate=180, mark size=3pt}}      
    
      \addplot+[const plot mark left, mark=none, color=black!50] table [x=Time, y=LB, col sep=comma] {Figures/ProgressLBBasic_Total.csv};
    
      \addplot+[const plot mark left, mark=none, color=black] table [x=Time, y=LB, col sep=comma] {Figures/ProgressLBComplete_Total.csv};

      \addplot+[const plot mark left, mark=square*, draw= black!50, dashed, mark options={fill=black!50,draw=black!50}] table [x=Time, y=ObjVal, col sep=comma]{Figures/ProgressUBBasic_Total.csv};
    
      \addplot+[scatter, only marks, 
        point meta=explicit symbolic,
        scatter/classes={True={mark=triangle*,mark size=3pt, rotate=180}, False={mark=*,mark size = 2pt}},const plot mark left, dashed, color=black] table [x=Time, y=ObjVal, col sep=comma,meta=SPHeur]{Figures/ProgressUBComplete_Total.csv};
    \end{axis}
\end{tikzpicture}
        \caption[size=scriptsize, labelfont=bf]{Progress of upper bound (UB) and lower bound (LB) of the objective value after preprocessing. For \completeBC{}, it is also displayed when the SP heuristic has found a new solution.}
        \label{fig:SolutionProcess}
    \end{subfigure}    
    \vspace{0.1cm}
    \begin{subfigure}{\textwidth}        
        \begin{tikzpicture}
    \pgfplotstableread[col sep=comma]{Figures/BarPlotSimpleVsComplete.csv}\Compare
    \begin{axis}[
            ybar,
            bar width=8,
            axis y line*=left,
            axis x line*=bottom,
            width = 13cm,
            height = 7cm,
            ylabel = {Call count},
            enlarge y limits = {upper=0.1},
            enlarge x limits = 0.07,
            xtick={0,1,...,8},
            xmin=0, xmax=8,
            xticklabels={Single\\vehicle\textsuperscript{a}, Check\\$R$, Check\\$C^{\mathrm{Inf}}$, Check\\heuristic, Check\\CP, 2P, RT,ITP / \\TT\textsuperscript{b},TT($\bar{P}$)},
            xticklabel style={font=\tiny,align=left},
            yticklabel = {\pgfmathparse{\tick}\pgfmathprintnumber{\pgfmathresult}},
            nodes near coords = {\pgfmathfloatifflags{\pgfplotspointmeta}{0}{}{\pgfmathprintnumber{\pgfplotspointmeta}}}, 
            every node near coord/.append style={font=\tiny, /pgf/number format/.cd,fixed, zerofill, precision=0},
        ]
            \addplot[draw,color=black, bar shift=-0pt,fill=black!10] table[x={CountIndexSimple}, y={CountSimple}] {\Compare};\label{plot:CountSimple}
            \addplot[draw,color=black,fill=black!30, bar shift=-0pt] table[x={CountIndexComplete}, y={CountComplete}] {\Compare};\label{plot:CountComplete}
        \end{axis}
        \begin{axis}[
            ybar,
            axis y line*=right,
            xtick=\empty,
            width = 13cm,
            height = 7cm,
            bar width= 8,
            yticklabel pos=right, ylabel near ticks,
            ylabel = Time,
            enlarge x limits = 0.07,
            ymin=0, ymax=44,
            xmin=0, xmax=8,
            legend style={font=\footnotesize,align=left},
            legend style={at={(1.1,0.5)},anchor=west,name=legendCompare},
            legend cell align={left},
            ]
        \addlegendimage{/pgfplots/refstyle=plot:CountSimple}\addlegendentry{Count \basicBC}
        \addlegendimage{/pgfplots/refstyle=plot:CountComplete}\addlegendentry{Count \completeBC}
        \addplot[draw, bar shift=0pt,fill=black!50] table[x={TimeIndexSimple}, y={TimeSimple}] {\Compare};
        \addlegendentry{Time \basicBC};
        \addplot[draw, bar shift=0pt,fill=black!80] table[x={TimeIndexComplete}, y={TimeComplete}] {\Compare};
        \addlegendentry{Time \completeBC};
        \coordinate (left) at (rel axis cs:0,0);
        \coordinate (right) at (rel axis cs:1,0);
    \end{axis}
    \node[anchor=north, inner sep=0pt, text width=2.5cm, 
    align=left,font=\scriptsize\linespread{0.95}\selectfont] at ($(legendCompare.south) + (0,-0.25)$) {
    \textsuperscript{a} represents all connected routes with $r(S)=1$\\[5pt]
    \textsuperscript{b} ITP for \basicBC / \\ TT for \completeBC};
\end{tikzpicture}
		\caption{Comparison of the call count (left axis) and the time spent (right axis) of the algorithm components for checking the feasibility of the routes.}
		\label{fig:BarPlot}
    \end{subfigure}    
    \vspace{0.1cm}    
    \begin{subfigure}{\textwidth}	
	    \begin{tikzpicture}
  \begin{sankeydiagram}
    \sffamily
    \sankeyset{
      ratio=1pt/4,
      outin steps=3,
      draw/.style={draw=none,line width=0pt},
      color/.style={fill/.style={fill=#1,fill opacity=.75}},
      shade/.style 2 args={fill/.style={left color=#1,
          right color=#2,fill opacity=.5}},
      @define HTML color/.code args={#1/#2}{\definecolor{#1}{HTML}{#2}},
      @define HTML color/.list={
        cyan/a6cee3,lime/b2df8a,red/a81919,orange/fdbf6f,
        violet/cab2d6,yellow/ffff99,blue/1f78b4,green/33a02c
      },
      @let country color/.code args={#1/#2}{\colorlet{#1}[rgb]{#2}},
      @let country color/.list={
        p0/gray,p1/orange,p2/orange,p3/cyan,p4/cyan!80!red,p5/cyan!80,
        p6/cyan!80!red,p19/blue,p7/cyan,p8/blue,p9/blue!50!green,
        p11/blue,p12/blue!80,p13/blue!60,p14/red,p15/cyan,p16/green
      },
    }
    \def\vdist{0.5cm}
    \def\hwidth{10pt}
    \def\hdist{2cm}        

    \sankeynode{name=p0,quantity=363}
    \sankeyadvance[color=p0]{p0}{\hwidth}
    \sankeyfork{p0}{4/p0-to-p2,359/p0-to-p3} 

    \sankeynode{name=p3,quantity=359,
      at={[xshift=\hdist,yshift=-\vdist]p0.right},anchor=right}
    \sankeynode{name=p2,quantity=4,
      at={[yshift=2*\vdist]p3.left},anchor=right}

    \sankeyfork{p2}{4/p2-from-p0}
    \sankeyfork{p3}{359/p3-from-p0}    

    \foreach \myNode in {p2,p3}{
      \sankeyadvance[color=\myNode]{\myNode}{\hwidth}
    }

    \sankeyfork{p3}{98/p3-to-p4,137/p3-to-p5,45/p3-to-p6,41/p3-to-p19,38/p3-to-p7}

    \sankeynode{name=p7,quantity=38,
      at={[xshift=\hdist,yshift=-2*\vdist]p3.right},anchor=right}
    \sankeynode{name=p19,quantity=41,
      at={[yshift=\vdist]p7.left},anchor=right}
    \sankeynode{name=p6,quantity=45,
      at={[yshift=\vdist]p19.left},anchor=right}
    \sankeynode{name=p5,quantity=137,
      at={[yshift=\vdist]p6.left},anchor=right}
    \sankeynode{name=p4,quantity=98,
      at={[yshift=\vdist]p5.left},anchor=right}

    \sankeyfork{p4}{98/p4-from-p3}
    \sankeyfork{p5}{137/p5-from-p3}
    \sankeyfork{p6}{45/p6-from-p3}
    \sankeyfork{p7}{38/p7-from-p3}
    \sankeyfork{p19}{41/p19-from-p3}

    \foreach \myNode in {p4,p5,p6,p7,p19}{
      \sankeyadvance[color=\myNode]{\myNode}{\hwidth}
    }

    
    \sankeynode{name=p8,quantity=41,
      at={[xshift=\hdist,yshift=\vdist]p19.right},anchor=right}

    \sankeyoutin[shade={p19}{p8}]{p19}{p8}
    
    \sankeyadvance[color=p8]{p8}{\hwidth}

    \sankeyfork{p8}{6/p8-to-p11,3/p8-to-p12,24/p8-to-p13,8/p8-to-p9}

    \sankeynode{name=p9,quantity=8,
      at={[xshift=\hdist]p8.right},anchor=right}
    \sankeynode{name=p13,quantity=24,
      at={[yshift=\vdist]p9.left},anchor=right}
    \sankeynode{name=p12,quantity=3,
      at={[yshift=\vdist]p13.left},anchor=right}
    \sankeynode{name=p11,quantity=6,
      at={[yshift=\vdist]p12.left},anchor=right}

    \sankeyfork{p9}{8/p9-from-p8}
    \sankeyfork{p11}{6/p11-from-p8}
    \sankeyfork{p12}{3/p12-from-p8}
    \sankeyfork{p13}{24/p13-from-p8}

    \foreach \myNode in {p9,p11,p12,p13}{
      \sankeyadvance[color=\myNode]{\myNode}{\hwidth}
    }

    \sankeynode{name=p16,quantity=46,
      at={[xshift=1.2*\hdist,yshift=-\vdist]p9.right},anchor=right}
    \sankeynode{name=p15,quantity=137,
      at={[yshift=1.5*\vdist]p16.left},anchor=right}
    \sankeynode{name=p14,quantity=180,
      at={[yshift=0.5*\vdist,yshift=\vdist]p15.left},anchor=right}

    \sankeyfork{p14}{4/p14-from-p2,98/p14-from-p4,45/p14-from-p6,6/p14-from-p11,3/p14-from-p12,24/p14-from-p13}
    \sankeyfork{p15}{137/p15-from-p5}
    \sankeyfork{p16}{8/p16-from-p9,38/p16-from-p7}

    \foreach \myNode in {p14,p15,p16}{
      \sankeyadvance[color=\myNode]{\myNode}{\hwidth}
    }

    \foreach \startnode/\endnode in {
    p2/p14, p4/p14, p6/p14, p11/p14, p12/p14, p13/p14,
    p5/p15,
    p9/p16,p7/p16}
    {
        \sankeyoutin[shade={\startnode}{\endnode}]{\startnode}{\endnode-from-\startnode}
    }
    
    \foreach \startnode/\nodes in {
      p0/{p2,p3},    p3/{p4,p5,p6,p7,p19},
      p8/{p9,p11,p12,p13},
      }
    {
      \foreach \endnode in \nodes {
        \sankeyoutin[shade={\startnode}{\endnode}]
        {\startnode-to-\endnode}{\endnode-from-\startnode}
      }
    }

    \foreach \myNode/\nodename in {p0/Routes
    ,p2/Disconnected, p3/Connected, p4/1D approx. inf, p5/Check $\mathcal{R}$, p6/Check $\mathcal{C}^{\mathrm{Inf}}$, p7/Heur. feas,p19/Heur. inf,p8/CP Check,p9/Feasible,p11/2P,p12/RT,p13/TT
    ,p14/Add constraint,p15/Accept route,p16/Add new route
      }
    {
      \node[anchor=west,inner sep=.1em,font=\footnotesize]
      at (\myNode) {\nodename\vphantom{Ag}}; 
    }

  \end{sankeydiagram}
\end{tikzpicture} 
	    \caption{Sankey diagram representing the call count to the feasibility checks of the routes in \completeBC{}}
	    \label{fig:Sankey}
    \end{subfigure}   
    \caption{Detailed evaluation of the \basicBC{} and \completeBC{} at instance \texttt{E016-03m} considering \allConstraints{}.}\label{fig:DetailedEvaluation}
\end{figure}

Figure~\ref{fig:SolutionProcess} provides more detailed insight into the different solution processes of both algorithm configurations using the instance \texttt{E016-03m}, showing the progression of UB and LB over time. It should be noted that we only show the time for the branch-and-cut part; therefore, the curves do not start at zero on the time axis. Comparing the LB curves reveals that \completeBC{} initially spends more time in the root node, where it is most likely to separate violated inequalities for fractional and integer solutions; however, due to the added inequalities, it subsequently improves the LB faster. This results in an LB of 274.8 after processing the root node for \completeBC{} compared to an LB of 260.1 for \basicBC{}. With regard to the UB, \completeBC{} can achieve a near-optimal solution in less than 12 seconds using an initial solution and the SP heuristic. In contrast, \basicBC{} requires over 30 seconds to obtain a feasible solution. This aligns with the findings in Table~\ref{tab:ComparisonSimpleVsCompleteBC} and highlights the difficulty that \basicBC{} faces in discovering new solutions without an external UB procedure. Overall, \completeBC{} converges much faster than \basicBC{}, proving optimality nearly three times faster in this instance.

Figure~\ref{fig:BarPlot} presents a breakdown of the individual algorithm components for feasibility checks of routes in integer solutions, detailing call counts and total time spent. We only show connected routes where the one-dimensional approximation is feasible ($r(S)=1$), i.e., routes where the loading would have to be evaluated. It is evident that route validation accounts for most of the solution process as it constitutes 41s out of the 57s of \basicBC{} and 15s out of the 20s for \completeBC{}. With \basicBC{}, all 1681 routes must be checked with the CP solver, resulting in 413 ITP inequalities; in contrast, \completeBC{} enhances efficiency through effective memory management and heuristic checking, which require almost zero additional time, meaning that only approximately one-sixth of all routes (41 out of 246) must undergo the time-consuming CP check. Here, 33 infeasible routes are identified, nine of which are successfully lifted to either 2P or RT inequalities. It should be noted that calls are sequential, and a route is checked multiple times if lifting proves unsuccessful.
Although the algorithm is time-bound due to lifting separation, the benefits of stronger constraints offset this limitation.

The Sankey diagram for \completeBC{} in Figure~\ref{fig:Sankey} (with the height of each bar corresponding to the number of calls) represents a flow of decisions that begins on the left with all the routes to be evaluated and ends on the right with the outcome of the check. Notably, only a tiny fraction of the examined routes are found to be disconnected from the depot and can be promptly excluded; for the remaining routes, a large part is already known to be feasible, can be declared as infeasible through the one-dimensional approximation, or is found feasible with the loading heuristic and added to $\mathcal{R}$. This also applies if the route is determined as feasible with the CP solver. In cases of infeasibility, the diagram illustrates the exclusion of routes using 2P or RT inequalities instead of TT inequalities. 

Overall, the results show that the algorithmic enhancements are highly effective and that \completeBC{} is clearly superior to \basicBC{}. Therefore, in the following experiments with the four remaining loading constraints, we will use only the configuration \completeBC{} and will refer to \completeBC{} as 'the B\&C algorithm'.

\subsection{Approximability as a determinant for hard-to-solve instances} \label{sec:HardInstances} 
Table~\ref{tab:ComparisonSimpleVsCompleteBC} highlights the greater difficulty of solving larger instances. However, solvability varies within the larger instances, prompting the question: What makes instances hard to solve? A feasible route must respect the actual loading constraints and the maximum weight. The latter is a one-dimensional parameter and is therefore easier to check, and if it is violated, the added GSEC can considerably restrict the solution space, as it applies independently of the customers' sequence. To explore the impact of both feasibility criteria, we categorize the instances into a group of heavy items ($\setOfHeavyItems$) and a group of light items ($\setOfLightItems$), using an average weight utilization threshold of 70\% in the routes of the best solutions found. The grouping of instances is detailed in Table~\ref{tab:ComputationalResults-AllConstraintsComplexity}.

\begin{table}[h]
    \centering
    \begin{threeparttable}        
    \caption{Results of gap (\%), total run time (s), and route structure for the \allConstraints{} variant.} \label{tab:ComputationalResults-AllConstraintsComplexity}
    {\begin{tabular}{L{1.2cm}lrrrrrrrr} 
     \toprule
             & & & & \multicolumn{6}{c}{Routes} \\
             \cmidrule(rl{3pt}){5-10}  \\[-9pt]
            Set & Instance & Gap & $t$ & \#{K} & $K^{\mathrm{min}}$ & $\overline{v}$ & $\overline{w}$ & $|C|$ & $\setOfItems$ \\
            \midrule
            \multirow{10}{=}{$\setOfHeavyItems$ $\overline{w} \geq 0.7$} & E016-03m & 0.0 & 20 & 4 & 3 & .54 & .72 & 3.75 & 8.00 \\
             & E016-05m & 0.0 & 2 & 5 & 5 & .33 & .94 & 3.00 & 5.20 \\
             & E021-04m & 0.0 & 257 & 4 & 4 & .56 & .97 & 5.00 & 9.25 \\
             & E021-06m & 0.0 & 6 & 6 & 6 & .37 & .95 & 3.33 & 6.00 \\
             & E022-04g & 0.0 & 680 & 5 & 4 & .59 & .75 & 4.20 & 9.00 \\
             & E022-06m & 0.0 & 15 & 6 & 6 & .37 & .94 & 3.50 & 6.67 \\
             & E026-08m & 0.0 & 36 & 8 & 8 & .45 & .96 & 3.12 & 6.25 \\
             & E031-09h & 0.0 & 1\,071 & 9 & 9 & .44 & .96 & 3.33 & 7.00 \\
             & E036-11h & 0.0 & 336 & 11 & 11 & .34 & .95 & 3.18 & 5.73 \\
             & E041-14h & 3.2 & TL & 14 & 14 & .34 & .95 & 2.86 & 5.64 \\ \midrule
             & Mean & 0.3 & 3\,123 & 7 & 7 & .43 & .91 & 3.53 & 6.87 \\ \midrule
             \multirow{9}{=}{$\setOfLightItems$ $\overline{w} < 0.7$} & E023-03g & 0.0 & 1\,684 & 5 & 4 & .57 & .45 & 4.40 & 9.20 \\
             & E023-05s & 0.0 & 18\,151 & 6 & 4 & .47 & .38 & 3.67 & 7.17 \\
             & E030-03g & 16.0 & TL & 6 & 5 & .67 & .47 & 4.83 & 10.33 \\
             & E030-04s & 15.2 & TL & 7 & 4 & .55 & .40 & 4.14 & 8.29 \\
             & E033-03n & 12.4 & TL & 6 & 4 & .63 & .43 & 5.33 & 10.17 \\
             & E033-04g & 23.0 & TL & 7 & 5 & .65 & .52 & 4.57 & 10.29 \\
             & E033-05s & 20.9 & TL & 6 & 5 & .72 & .61 & 5.33 & 11.33 \\
             & E045-04f & 17.9 & TL & 10 & 6 & .59 & .36 & 4.40 & 9.40 \\
             & E051-05e & 18.0 & TL & 10 & 7 & .61 & .49 & 5.00 & 9.90 \\ \midrule
             & Mean & 13.7 & 24\,621 & 7 & 5 & .61 & .46 & 4.63 & 9.56 \\ \bottomrule
\end{tabular}

 }
    \begin{tablenotes}[para,flushleft]
        \item[] The routes are described with the number of vehicles \#{K}, the lower bound for required vehicles $K^{\mathrm{min}}$, the average utilization of volume $\overline{v}$ and weight $\overline{w}$, as well as the average number of visited customers $|C|$ and loaded items $\setOfItems$ per route.        
    \end{tablenotes}
    \end{threeparttable}
\end{table}

With a single exception, instances in $\setOfHeavyItems$ can be optimally solved in less than 1200 seconds. The number of routes in the best solution \#{K} predominantly aligns with the lower bound $K^{\mathrm{min}}$, and vehicles exhibit an average weight utilization of over 90\%. In contrast, achieving optimality is notably more challenging for instances in $\setOfLightItems$, for several reasons. First, it is observed that \#{K} consistently exceeds $K^{\mathrm{min}}$, meaning that the one-dimensional approximation deviates further from actual feasible solutions with respect to the loading constraints; as a result, the one-dimensional approximation in the integer solution feasibility check is less effective and the expensive CP check must be called more often, with only weaker inequalities are added. Furthermore, lighter items allow for incorporation of more customers with more items on a route, intensifying the complexity of the loading problems. Finally, cuts for fractional solutions only consider a one-dimensional metric and are less likely to be violated, leading to a less positive impact on lower bounds.  In general, we can state that instances become more challenging to solve as the loading exerts a more substantial influence while the weight's impact weakens.

\subsection{Comparison to best-known heuristic results} \label{sec:ComparisonBKS}

To the best of our knowledge, no exact approach considers all five variants introduced in \cite{Gendreau2006}. Therefore, we compare the results of our algorithm to existing heuristic solutions. We take the best solution for each instance and variant from the literature; if several references provide the same objective, the first published work is referenced.

\begin{landscape}
\begin{table}
    \renewcommand{\arraystretch}{0.8}
    \centering
    \caption{Comparison to best-known solutions for all loading variants.}\label{tab:BKS_Comparison}
    \begin{threeparttable}
    \begin{filecontents*}{Chapters/AllVariants-TransformedFlat.csv}
Instance,AllConstraints-VehicleCount,AllConstraints-Objective,AllConstraints-GapPCT,AllConstraints-RPDDeltaBKS,AllConstraints-AvgRouteWeight,AllConstraints-AvgRouteVolume,AllConstraints-CiteMarkBKS,AllConstraints-ADVehicleCount,AllConstraints-RPDVehicleCount,AllConstraints-RPDVolumeWeightWeightUtil,AllConstraints-RPDVolumeWeightVolUtil,AllConstraints-RPDVolumeWeight,NoFragility-VehicleCount,NoFragility-Objective,NoFragility-GapPCT,NoFragility-ADVehicleCount,NoFragility-RPDVehicleCount,NoFragility-RPDVolumeWeightWeightUtil,NoFragility-RPDVolumeWeightVolUtil,NoFragility-RPDVolumeWeight,NoFragility-RPDDeltaBKS,NoFragility-AvgRouteWeight,NoFragility-AvgRouteVolume,NoLifo-ADVehicleCount,NoLifo-RPDVehicleCount,NoLifo-RPDVolumeWeightWeightUtil,NoLifo-RPDVolumeWeightVolUtil,NoLifo-RPDVolumeWeight,NoFragility-CiteMarkBKS,NoLifo-VehicleCount,NoLifo-Objective,NoSupport-ADVehicleCount,NoSupport-RPDVehicleCount,NoSupport-RPDVolumeWeightWeightUtil,NoSupport-RPDVolumeWeightVolUtil,NoSupport-RPDVolumeWeight,NoLifo-GapPCT,NoLifo-RPDDeltaBKS,NoLifo-AvgRouteWeight,LoadingOnly-ADVehicleCount,LoadingOnly-RPDVehicleCount,LoadingOnly-RPDVolumeWeightWeightUtil,LoadingOnly-RPDVolumeWeightVolUtil,LoadingOnly-RPDVolumeWeight,NoLifo-AvgRouteVolume,NoLifo-CiteMarkBKS,NoSupport-VehicleCount,VolumeWeightApproximation-ADVehicleCount,VolumeWeightApproximation-RPDVehicleCount,VolumeWeightApproximation-RPDVolumeWeightWeightUtil,VolumeWeightApproximation-RPDVolumeWeightVolUtil,VolumeWeightApproximation-RPDVolumeWeight,NoSupport-Objective,NoSupport-GapPCT,NoSupport-RPDDeltaBKS,NoSupport-AvgRouteWeight,NoSupport-AvgRouteVolume,NoSupport-CiteMarkBKS,LoadingOnly-VehicleCount,LoadingOnly-Objective,LoadingOnly-GapPCT,LoadingOnly-RPDDeltaBKS,LoadingOnly-AvgRouteWeight,LoadingOnly-AvgRouteVolume,LoadingOnly-CiteMarkBKS,VolumeWeightApproximation-VehicleCount,VolumeWeightApproximation-Objective,VolumeWeightApproximation-GapPCT,VolumeWeightApproximation-RPDDeltaBKS,VolumeWeightApproximation-AvgRouteWeight,VolumeWeightApproximation-AvgRouteVolume,VolumeWeightApproximation-CiteMarkBKS,adjustedInstance
E016-03m,4.0,301.6582384371,0.0,-0.03,0.7166666667,0.5354222222,\textsuperscript{a},1.0,0.3333333333333333,-25.000000000000007,-25.000000003501917,8.127068716116879,4.0,301.6582384371,0.0,1.0,0.3333333333333333,-25.000000000000007,-25.000000003501917,8.127068716116879,-0.03,0.7166666667,0.5354222222,1.0,0.3333333333333333,-25.000000000000007,-25.000000003501917,6.69066022905316,\textsuperscript{d},4.0,297.6508750726,1.0,0.3333333333333333,-25.000000000000007,-25.000000003501917,6.69066022905316,0.0,0.0,0.7166666667,1.0,0.3333333333333333,-25.000000000000007,-25.000000003501917,6.69066022905316,0.5354222222,\textsuperscript{d},4.0,,,,,,297.6508750726,0.0,0.0,0.7166666667,0.5354222222,\textsuperscript{d},4.0,297.6508750726,0.0,0.0,0.7166666667,0.5354222222,\textsuperscript{e},3.0,278.9849406064,0.0,0.0,0.9555555556,0.7138962963000001,\textsuperscript{j},E016-03m\textsuperscript{$\setOfHeavyItems$}
E016-05m,5.0,334.9638862376,0.0,0.0,0.9381818182,0.3322,\textsuperscript{b},0.0,0.0,0.0,0.0,0.0,5.0,334.9638862376,0.0,0.0,0.0,0.0,0.0,0.0,0.0,0.9381818182,0.3322,0.0,0.0,0.0,0.0,0.0,\textsuperscript{e},5.0,334.9638862376,0.0,0.0,0.0,0.0,0.0,0.0,0.0,0.9381818182,0.0,0.0,0.0,0.0,0.0,0.3322,\textsuperscript{e},5.0,,,,,,334.9638862376,0.0,0.0,0.9381818182,0.3322,\textsuperscript{e},5.0,334.9638862376,0.0,0.0,0.9381818182,0.3322,\textsuperscript{e},5.0,334.9638862376,0.0,0.0,0.9381818182,0.3322,\textsuperscript{j},E016-05m\textsuperscript{$\setOfHeavyItems$}
E021-04m,4.0,385.531612043,0.0,0.0,0.9676470588,0.5601833333,\textsuperscript{c},0.0,0.0,0.0,0.0,7.569521693719902,4.0,373.0101194861,0.0,0.0,0.0,0.0,0.0,4.075823840774506,0.0,0.9676470588,0.5601833333,0.0,0.0,0.0,0.0,1.0795714830415148,\textsuperscript{c},4.0,362.2714828967,0.0,0.0,0.0,0.0,1.0795714830415148,0.0,0.0,0.9676470588,0.0,0.0,0.0,0.0,1.0795714830415148,0.5601833333,\textsuperscript{g},4.0,,,,,,362.2714828967,0.0,0.0,0.9676470588,0.5601833333,\textsuperscript{c},4.0,362.2714828967,0.0,0.0,0.9676470588,0.5601833333,\textsuperscript{e},4.0,358.4022741504,0.0,0.0,0.9676470588,0.5601833333,\textsuperscript{j},E021-04m\textsuperscript{$\setOfHeavyItems$}
E021-06m,6.0,430.8846754242,0.0,-1.44,0.9454022989,0.36855555560000003,\textsuperscript{a},0.0,0.0,0.0,0.0,0.0,6.0,430.8846754242,0.0,0.0,0.0,0.0,0.0,0.0,0.0,0.9454022989,0.36855555560000003,0.0,0.0,0.0,0.0,0.0,\textsuperscript{c},6.0,430.8846754242,0.0,0.0,0.0,0.0,0.0,0.0,0.0,0.9454022989,0.0,0.0,0.0,0.0,0.0,0.36855555560000003,\textsuperscript{e},6.0,,,,,,430.8846754242,0.0,0.0,0.9454022989,0.36855555560000003,\textsuperscript{e},6.0,430.8846754242,0.0,0.0,0.9454022989,0.36855555560000003,\textsuperscript{e},6.0,430.8846754242,0.0,0.0,0.9454022989,0.36855555560000003,\textsuperscript{j},E021-06m\textsuperscript{$\setOfHeavyItems$}
E022-04g,5.0,427.5637664265,0.0,-2.04,0.75,0.58696,\textsuperscript{a},1.0,0.25,-20.0,-20.0,13.931999822276769,5.0,395.6358415971,0.0,1.0,0.25,-20.0,-20.0,5.424234170403673,-8.32,0.75,0.58696,1.0,0.25,-20.0,-20.0,5.424234170403673,\textsuperscript{c},5.0,395.6358415971,0.0,0.0,0.0,0.0,3.8526960481615267,0.0,-2.67,0.75,0.0,0.0,0.0,0.0,1.106190265228123,0.58696,\textsuperscript{c},4.0,,,,,,389.738176677,0.0,-1.49,0.9375,0.7337,\textsuperscript{c},4.0,379.4310956208,0.0,0.0,0.9375,0.7337,\textsuperscript{f},4.0,375.279787148,0.0,0.0,0.9375,0.7337,\textsuperscript{j},E022-04g\textsuperscript{$\setOfHeavyItems$}
E022-06m,6.0,498.1572394012,0.0,0.0,0.9375,0.3722666667,\textsuperscript{d},0.0,0.0,0.0,0.0,0.46576173977106144,6.0,495.8477702001,0.0,0.0,0.0,0.0,0.0,0.0,0.0,0.9375,0.3722666667,0.0,0.0,0.0,0.0,0.0,\textsuperscript{c},6.0,495.8477702001,0.0,0.0,0.0,0.0,0.0,0.0,0.0,0.9375,0.0,0.0,0.0,0.0,0.0,0.3722666667,\textsuperscript{c},6.0,,,,,,495.8477702001,0.0,0.0,0.9375,0.3722666667,\textsuperscript{g},6.0,495.8477702001,0.0,0.0,0.9375,0.3722666667,\textsuperscript{e},6.0,495.8477702001,0.0,0.0,0.9375,0.3722666667,\textsuperscript{j},E022-06m\textsuperscript{$\setOfHeavyItems$}
E023-03g,5.0,757.8756055701,0.0019801281,-1.25,0.4528444444,0.5722888889000001,\textsuperscript{a},1.0,0.25,-20.000000014132887,-19.9999999972042,15.14194488935266,5.0,750.3772954739,0.0,1.0,0.25,-20.000000014132887,-19.9999999972042,14.002747372617067,-1.26,0.4528444444,0.5722888889000001,1.0,0.25,-20.000000014132887,-19.9999999972042,11.288952398513135,\textsuperscript{c},5.0,732.5148300503,1.0,0.25,-20.000000014132887,-19.9999999972042,12.764776339497336,0.0,0.0,0.4528444444,0.0,0.0,0.0,0.0,5.8876253652429416,0.5722888889000001,\textsuperscript{f},5.0,,,,,,742.2288483784,0.0,-0.12,0.4528444444,0.5722888889000001,\textsuperscript{f},4.0,696.9627642922,0.0,-2.5300000000000002,0.5660555556,0.7153611111,\textsuperscript{h},4.0,658.2098350852,0.0,0.0,0.5660555556,0.7153611111,\textsuperscript{j},E023-03g\textsuperscript{$\setOfLightItems$}
E023-05s,6.0,798.647378203,0.0,-0.66,0.3773703704,0.47164444440000003,\textsuperscript{a},2.0,0.5,-33.33333333333333,-33.333333342756625,20.19459668348126,5.0,779.6612491802,0.0,1.0,0.25,-20.000000014132887,-20.00000000848097,17.337227858693183,-0.33,0.4528444444,0.5659733333,1.0,0.25,-20.000000014132887,-20.00000000848097,9.962330092265985,\textsuperscript{c},5.0,730.6578586103,1.0,0.25,-20.000000014132887,-20.00000000848097,10.637191513352253,0.0,-0.61,0.4528444444,1.0,0.25,-20.000000014132887,-20.00000000848097,9.962330092265985,0.5659733333,\textsuperscript{c},5.0,,,,,,735.1420560657,0.0,-5.24,0.4528444444,0.5659733333,\textsuperscript{c},5.0,730.6578586103,0.0,0.0,0.4528444444,0.5659733333,\textsuperscript{f},4.0,664.4619643811,0.0,0.0,0.5660555556,0.7074666667,\textsuperscript{j},E023-05s\textsuperscript{$\setOfLightItems$}
E026-08m,8.0,630.1275729831,0.0,0.0,0.9557291667000001,0.4518472222,\textsuperscript{a},0.0,0.0,0.0,0.0,3.69893727806681,8.0,630.1275729831,0.0,0.0,0.0,0.0,0.0,3.69893727806681,0.0,0.9557291667000001,0.4518472222,0.0,0.0,0.0,0.0,3.69893727806681,\textsuperscript{f},8.0,630.1275729831,0.0,0.0,0.0,0.0,3.69893727806681,0.0,0.0,0.9557291667000001,0.0,0.0,0.0,0.0,3.69893727806681,0.4518472222,\textsuperscript{b},8.0,,,,,,630.1275729831,0.0,0.0,0.9557291667000001,0.4518472222,\textsuperscript{b},8.0,630.1275729831,0.0,0.0,0.9557291667000001,0.4518472222,\textsuperscript{e},8.0,607.6509456345,0.0,0.0,0.9557291667000001,0.4518472222,\textsuperscript{j},E026-08m\textsuperscript{$\setOfHeavyItems$}
E030-03g,6.0,769.3189972358,15.9692356148,-6.22,0.47222222220000004,0.6699592593,\textsuperscript{d},1.0,0.2,-16.666666675490188,-16.66666666044739,28.91609365956328,6.0,759.650527411,13.7172389425,1.0,0.2,-16.666666675490188,-16.66666666044739,27.295931716394016,-0.81,0.47222222220000004,0.6699592593,1.0,0.2,-16.666666675490188,-16.66666666044739,18.356144894242853,\textsuperscript{c},6.0,706.3015029542,1.0,0.2,-16.666666675490188,-16.66666666044739,18.90805377721973,5.5402860184,-4.54,0.47222222220000004,0.0,0.0,0.0,0.0,8.520488204669146,0.6699592593,\textsuperscript{c},6.0,,,,,,709.595071479,4.8772672076,-4.82,0.47222222220000004,0.6699592593,\textsuperscript{c},5.0,647.6062902249,0.0,-5.8100000000000005,0.5666666667,0.8039511111000001,\textsuperscript{i},5.0,596.7594699754,0.0,0.0,0.5666666667,0.8039511111000001,\textsuperscript{j},E030-03g\textsuperscript{$\setOfLightItems$}
E030-04s,7.0,728.3201915273,15.2241316423,-5.2,0.4047619048,0.5537269841,\textsuperscript{a},3.0,0.75,-42.85714284907563,-42.857142858617095,31.613975069412692,7.0,728.0054337266,15.658307556,3.0,0.75,-42.85714284907563,-42.857142858617095,31.557095518610502,-5.19,0.4047619048,0.5537269841,3.0,0.75,-42.85714284907563,-42.857142858617095,29.793230638389094,\textsuperscript{c},7.0,718.2446282596,3.0,0.75,-42.85714284907563,-42.857142858617095,30.28635069899116,14.0661710622,0.0,0.4047619048,2.0,0.5,-33.33333333333333,-33.33333333333333,27.749350506712002,0.5537269841,\textsuperscript{c},7.0,,,,,,720.9734364792,14.1935583364,0.0,0.4047619048,0.5537269841,\textsuperscript{c},6.0,706.9342855078,12.3514833889,-1.45,0.47222222220000004,0.6460148148,\textsuperscript{f},4.0,553.3760310356,0.0,0.0,0.7083333333,0.9690222222,\textsuperscript{j},E030-04s\textsuperscript{$\setOfLightItems$}
E031-09h,9.0,610.2341337694,0.0,0.0,0.9640522876000001,0.4408567901,\textsuperscript{a},0.0,0.0,0.0,0.0,0.03783733954209724,9.0,610.0033247402,0.0038764005,0.0,0.0,0.0,0.0,0.0,-0.04,0.9640522876000001,0.4408567901,0.0,0.0,0.0,0.0,0.0,\textsuperscript{c},9.0,610.0033247402,0.0,0.0,0.0,0.0,0.0,0.0005130089,0.0,0.9640522876000001,0.0,0.0,0.0,0.0,0.0,0.4408567901,\textsuperscript{c},9.0,,,,,,610.0033247402,0.0039738967,0.0,0.9640522876000001,0.4408567901,\textsuperscript{c},9.0,610.0033247402,0.008298523700000001,0.0,0.9640522876000001,0.4408567901,\textsuperscript{f},9.0,610.0033247402,0.0,0.0,0.9640522876000001,0.4408567901,\textsuperscript{j},E031-09h\textsuperscript{$\setOfHeavyItems$}
E033-03n,6.0,2619.7827583108,12.4158560976,-0.99,0.43230263160000004,0.626462963,\textsuperscript{d},2.0,0.5,-33.33333333333333,-33.33333332623883,22.948881729176467,6.0,2545.611827101,9.7355331921,2.0,0.5,-33.33333333333333,-33.33333332623883,19.46796980236596,-2.29,0.43230263160000004,0.626462963,1.0,0.25,-20.00000000308426,-19.999999991486593,8.342981581704779,\textsuperscript{c},5.0,2308.5616651394,1.0,0.25,-20.00000000308426,-19.999999991486593,12.19521447571984,3.9883690501,-5.29,0.5187631579,1.0,0.25,-20.00000000308426,-19.999999991486593,7.591728888990193,0.7517555556000001,\textsuperscript{c},5.0,,,,,,2390.6446672359,3.6951721104,-2.77,0.5187631579,0.7517555556000001,\textsuperscript{c},5.0,2292.5540461694,3.8111637812,-0.58,0.5187631579,0.7517555556000001,\textsuperscript{i},4.0,2130.7902288054,0.0,0.0,0.6484539474000001,0.9396944444,\textsuperscript{j},E033-03n\textsuperscript{$\setOfLightItems$}
E033-04g,7.0,1320.8360534148,23.0417060022,-3.1,0.5244642857,0.6474603175,\textsuperscript{c},2.0,0.4,-28.57142857337419,-28.571428563548494,31.183180508061763,7.0,1304.6004811055,21.9271401411,2.0,0.4,-28.57142857337419,-28.571428563548494,29.570691200705053,-2.81,0.5244642857,0.6474603175,1.0,0.2,-16.666666666666682,-16.66666665931193,22.38677617309061,\textsuperscript{c},6.0,1232.2682359473,1.0,0.2,-16.666666666666682,-16.66666665931193,19.62377357369656,18.6689174597,-4.21,0.611875,1.0,0.2,-16.666666666666682,-16.66666665931193,17.242515841083343,0.7553703704,\textsuperscript{c},6.0,,,,,,1204.4485609339,15.0600221912,-6.95,0.611875,0.7553703704,\textsuperscript{c},6.0,1180.4725370753,13.8378239626,-0.32,0.611875,0.7553703704,\textsuperscript{c},5.0,1006.8638740876,0.0,0.0,0.7342500000000001,0.9064444444,\textsuperscript{j},E033-04g\textsuperscript{$\setOfLightItems$}
E033-05s,6.0,1265.725296606,20.872134067,-5.42,0.611875,0.7247814815,\textsuperscript{c},1.0,0.2,-16.666666666666682,-16.666666666666675,28.40201198256413,6.0,1210.4524174345,16.8245355738,1.0,0.2,-16.666666666666682,-16.666666666666675,22.794832515802625,-8.54,0.611875,0.7247814815,1.0,0.2,-16.666666666666682,-16.666666666666675,18.60050775203218,\textsuperscript{c},6.0,1169.1067806044,1.0,0.2,-16.666666666666682,-16.666666666666675,18.634894420981325,15.2961756911,-1.1,0.611875,1.0,0.2,-16.666666666666682,-16.666666666666675,13.897518184154373,0.7247814815,\textsuperscript{c},6.0,,,,,,1169.4457478533,13.5089040356,-2.7800000000000002,0.611875,0.7247814815,\textsuperscript{c},6.0,1122.7469707087,11.7396352178,-1.65,0.611875,0.7247814815,\textsuperscript{f},5.0,985.7519185742,0.0,0.0,0.7342500000000001,0.8697377778000001,\textsuperscript{j},E033-05s\textsuperscript{$\setOfLightItems$}
E036-11h,11.0,698.6054332259,0.0,0.0,0.9484396201,0.3389030303,\textsuperscript{b},0.0,0.0,0.0,0.0,0.0,11.0,698.6054332259,0.0,0.0,0.0,0.0,0.0,0.0,0.0,0.9484396201,0.3389030303,0.0,0.0,0.0,0.0,0.0,\textsuperscript{f},11.0,698.6054332259,0.0,0.0,0.0,0.0,0.0,0.0,0.0,0.9484396201,0.0,0.0,0.0,0.0,0.0,0.3389030303,\textsuperscript{b},11.0,,,,,,698.6054332259,0.0,0.0,0.9484396201,0.3389030303,\textsuperscript{b},11.0,698.6054332259,0.0,0.0,0.9484396201,0.3389030303,\textsuperscript{e},11.0,698.6054332259,0.0,0.0,0.9484396201,0.3389030303,\textsuperscript{j},E036-11h\textsuperscript{$\setOfHeavyItems$}
E041-14h,14.0,866.3977009136,3.1768315629,0.0,0.9500000000000001,0.3405206349,\textsuperscript{d},0.0,0.0,0.0,0.0,0.5350086111793289,14.0,866.3977009136,3.1564711289,0.0,0.0,0.0,0.0,0.5350086111793289,0.0,0.9500000000000001,0.3405206349,0.0,0.0,0.0,0.0,0.0,\textsuperscript{b},14.0,861.787065901,0.0,0.0,0.0,0.0,0.1716143677619248,2.7423793136,-0.53,0.9500000000000001,0.0,0.0,0.0,0.0,0.0,0.3405206349,\textsuperscript{b},14.0,,,,,,863.2660163256,2.7233218263000003,0.0,0.9500000000000001,0.3405206349,\textsuperscript{c},14.0,861.787065901,2.5332788293,0.0,0.9500000000000001,0.3405206349,\textsuperscript{b},14.0,861.787065901,2.6654110376,0.0,0.9500000000000001,0.3405206349,\textsuperscript{j},E041-14h\textsuperscript{$\setOfHeavyItems$}
E045-04f,10.0,1204.7332702772,17.9480142417,0.24,0.3592039801,0.59252,\textsuperscript{c},4.0,0.6666666666666666,-40.000000003340716,-39.99999999797475,31.239096504439125,9.0,1128.0252341676,11.0454992122,3.0,0.5,-33.33333334446906,-33.333333326582505,22.88281250197488,-3.06,0.3991155334,0.6583555556,3.0,0.5,-33.33333334446906,-33.333333326582505,18.6912299409039,\textsuperscript{c},9.0,1089.5478360375,2.0,0.3333333333333333,-25.00000000835179,-24.99999999746844,20.67437834774356,10.1637163364,-3.52,0.3991155334,2.0,0.3333333333333333,-25.00000000835179,-24.99999999746844,16.513214721968257,0.6583555556,\textsuperscript{c},8.0,,,,,,1107.752509258,9.4668765813,-2.0,0.44900497510000004,0.74065,\textsuperscript{c},8.0,1069.5543473035,10.2132054546,-0.8200000000000001,0.44900497510000004,0.74065,\textsuperscript{i},6.0,917.968274977,0.0,0.0,0.5986733002,0.9875333333,\textsuperscript{j},E045-04f\textsuperscript{$\setOfLightItems$}
E051-05e,10.0,721.1234587472,18.0049679916,-2.77,0.485625,0.6099488889,\textsuperscript{c},3.0,0.42857142857142855,-30.000000000000014,-30.000000002295273,21.66725749575513,9.0,702.7825333764,16.1324860498,2.0,0.2857142857142857,-22.222222227027043,-22.222222220947067,18.57279417920142,-1.43,0.5395833333,0.6777209877,2.0,0.2857142857142857,-22.222222227027043,-22.222222220947067,12.281985406115972,\textsuperscript{c},9.0,665.4968258316,2.0,0.2857142857142857,-22.222222227027043,-22.222222220947067,13.560351944556306,12.3635185686,-1.05,0.5395833333,1.0,0.14285714285714285,-12.500000000000009,-12.500000005738187,9.420189845509814,0.6777209877,\textsuperscript{c},9.0,,,,,,673.0737213639,11.7621239292,-2.25,0.5395833333,0.6777209877,\textsuperscript{c},8.0,648.5349253551,9.6647887046,-1.15,0.60703125,0.7624361111,\textsuperscript{h},7.0,592.7013344345,0.0,0.0,0.6937500000000001,0.8713555556,\textsuperscript{j},E051-05e\textsuperscript{$\setOfLightItems$}
Mean,7.105263157894737,808.9730141449369,6.666045123589473,-1.5200000000000002,0.6944362503263158,0.5156057201526317,,1.105263157894737,0.2357142857142857,-16.127819549934053,-16.127819548381645,13.982798616972598,6.947368421052632,791.9106085379841,5.694794115626315,0.9473684210526315,0.20626566416040096,-14.665831247247468,-14.66583124380185,11.860167120152942,-1.795263157894737,0.7033490903578947,0.5276023803263159,0.8421052631578947,0.1825814536340852,-13.337510445302385,-13.337510441749806,8.768291686201247,,6.842105263157895,761.6041100901632,0.7368421052631579,0.16065162907268168,-11.846282374980424,-11.846282371796436,9.093603394623315,4.359476132052632,-1.2378947368421052,0.7125002082842105,0.5263157894736842,0.11629072681704261,-8.903508773275561,-8.903508771894106,6.808437942420299,0.5398762037736843,,6.7368421052631575,,,,,,766.6665175173842,3.962695795510526,-1.4957894736842106,0.7249943894263157,0.5519306482157896,,6.526315789473684,747.2419582920736,3.3768251506684206,-0.7531578947368421,0.7430240625736841,0.5758289128526316,,6.0,692.5943702433842,0.1402847914526316,0.0,0.8061313770894736,0.6538682208421054,,Mean
\end{filecontents*}
\pgfplotstabletypeset[
    	col sep=comma,
	    set thousands separator={$\,$},
    	columns={adjustedInstance, AllConstraints-Objective,AllConstraints-GapPCT,AllConstraints-RPDDeltaBKS,AllConstraints-CiteMarkBKS, 
        NoFragility-Objective,NoFragility-GapPCT,NoFragility-RPDDeltaBKS,NoFragility-CiteMarkBKS,
        NoLifo-Objective,NoLifo-GapPCT,NoLifo-RPDDeltaBKS,NoLifo-CiteMarkBKS,
        NoSupport-Objective,NoSupport-GapPCT,NoSupport-RPDDeltaBKS,NoSupport-CiteMarkBKS,
        LoadingOnly-Objective,LoadingOnly-GapPCT,LoadingOnly-RPDDeltaBKS,LoadingOnly-CiteMarkBKS
        },
    	display columns/adjustedInstance/.style={string type},
    	display columns/AllConstraints-CiteMarkBKS/.style={string type},
    	display columns/NoFragility-CiteMarkBKS/.style={string type},
    	display columns/NoLifo-CiteMarkBKS/.style={string type},
    	display columns/NoSupport-CiteMarkBKS/.style={string type},
    	display columns/LoadingOnly-CiteMarkBKS/.style={string type},  
    	column type=r, 
    	precision=2, fixed, fixed zerofill,
        columns/adjustedInstance/.style={string type, column type={l}},
        columns/AllConstraints-Objective/.style={column type={@{\hskip 4pt}r}}, 
        columns/AllConstraints-GapPCT/.append style={precision=1,
        },
        columns/AllConstraints-RPDDeltaBKS/.append style={column type={r@{}}}, 
        columns/AllConstraints-CiteMarkBKS/.append style={string type, column type={@{}l}},
        columns/NoFragility-GapPCT/.append style={precision=1,
        },
        columns/NoFragility-RPDDeltaBKS/.append style={column type={r@{}}},
        columns/NoFragility-CiteMarkBKS/.append style={string type, column type={@{}l}},
        columns/NoLifo-GapPCT/.append style={precision=1,
        },
        columns/NoLifo-RPDDeltaBKS/.append style={column type={r@{}}},
        columns/NoLifo-CiteMarkBKS/.append style={string type, column type={@{}l}},
        columns/NoSupport-GapPCT/.append style={precision=1,
        },
        columns/NoSupport-RPDDeltaBKS/.append style={column type={r@{}}},
        columns/NoSupport-CiteMarkBKS/.append style={string type, column type={@{}l}},
        columns/LoadingOnly-GapPCT/.append style={precision=1,
        },
        columns/LoadingOnly-RPDDeltaBKS/.append style={column type={r@{}}},
        columns/LoadingOnly-CiteMarkBKS/.append style={string type, column type={@{}l}},
        every head row/.append style={
            output empty row,
            before row={\toprule
            & \multicolumn{4}{c}{\allConstraints{}} & \multicolumn{4}{c}{\noFragility{}} & \multicolumn{4}{c}{\noLifo{}} & \multicolumn{4}{c}{\noSupport{}} & \multicolumn{4}{c}{\loadingOnly{}} \\
            \cmidrule(rl{3pt}){2-5} \cmidrule(rl{3pt}){6-9} \cmidrule(rl{3pt}){10-13} \cmidrule(rl{3pt}){14-17} \cmidrule(rl{3pt}){18-21} \\[-6pt]
            Instance & Obj & Gap & $\Delta_{\mathrm{BKS}}$ & 
            & Obj & Gap & \llap{$\Delta_{\mathrm{BKS}}$} &  
            & Obj & Gap & \llap{$\Delta_{\mathrm{BKS}}$} &  
            & Obj & Gap & \llap{$\Delta_{\mathrm{BKS}}$} &  
            & Obj & Gap & \llap{$\Delta_{\mathrm{BKS}}$} & \\
            }, 
            after row={\midrule}},
        every last row/.append style={after row=\bottomrule, before row=\midrule},
    ]{Chapters/AllVariants-TransformedFlat.csv}
    \begin{tablenotes}[flushleft]
        \item[] For each variant, the objective, gap (\%) and the relative percentage deviation to the best-known solution $\Delta_{\mathrm{BKS}}$ (\%) are presented. Previously best-known solutions were first found by:       
        \setlength{\columnsep}{0.8cm} 
        \setlength{\multicolsep}{0cm}
        \begin{multicols}{3}
            \item[a] \cite{Wang2010}
            \item[b] \cite{Tarantilis2009}
            \item[c] \cite{Zhang2015}
            \item[d] \cite{Bortfeldt2012}
            \item[e] \cite{Gendreau2006}
            \item[f] \cite{Zhu2012}
            \item[g] \cite{Fuellerer2010}
            \item[h] \cite{Wei2014}
            \item[i] \cite{Lacomme2013}
        \end{multicols}
        \item[] Some references are excluded from the analysis because validation is not possible as only summarized results are published. In addition, comparability cannot be established for the following reasons. \cite{Escobar2016} assume that weight equals volume. \cite{Escobar2015} allow more rotations as defined in the benchmarks, while \cite{Hokama2016} do not consider rotation nor fragility. \cite{Miao2012} report incomplete solutions. In addition, all references are excluded that report at least one objective value that is less than the lower bounds of our method \citep{Ceschia2013, Ruan2013, Wei2014, Tao2015, Mahvash2017, Koch2020}. A standardized solution format and the publishing of detailed solutions simplify comparability and validation \citep{Krebs2023a}. It should be noted that items can hover in our solutions if support constraints are disabled as the supported area can be zero. This is prohibited in \citep{Krebs2023a}. Thus, our solutions for variants \noSupport{} and \loadingOnly{} might be infeasible using their solution validator. While with the \loadingOnly{} variant all floating items could be lowered enough to touch an underlying item or the container floor to generate a feasible solution respecting their definition, this is not viable with the \noSupport{} variant due to the fragility constraint.
    \end{tablenotes}
    \end{threeparttable}
\end{table}
\end{landscape}

The B\&C algorithm produces proven optimal solutions for eleven instances across all loading variants and additionally solves instance \texttt{E030-03g} to optimality for the \loadingOnly{} variant (see Table~\ref{tab:BKS_Comparison}). We can prove the optimality of 39 best-known solutions, find the same solutions but not prove optimality in six instances, generate 49 new best solutions, and produce worse results in only one instance. Notably, existing heuristics demonstrate mostly very good results for small instances, instances with heavy items, and for the \noSupport{}, \noLifo{}, and \loadingOnly{} variants. However, for instances with light items and with more complex constraints in the loading problem, it becomes evident that finding good solutions is challenging: here, existing heuristics prove less effective, indicating a need for algorithmic improvement, especially for the loading problem. As a result, current best solutions can frequently be improved by more than 2\% and at maximum by more than 8\%.

\subsection{Impact of constraints in the loading subproblem} \label{sec:LoadingVariantComparison} 

Table~\ref{tab:ComputationalResults_Brief} explores the impact of loading constraints on solutions and their solvability. The benchmark for comparing the different variants is provided by solving each instance as CVRP with limited volume and weight.

\begin{table}[htpb]
    \centering
    \begin{threeparttable}
    \caption{Optimal objective values of the CVRP with limited volume and weight compared to all other loading variants.} \label{tab:ComputationalResults_Brief}
    {\footnotesize \pgfplotstableset{
    row style/.style 2 args={
        every row #1 column 3/.style={#2},
        every row #1 column 6/.style={#2},
        every row #1 column 9/.style={#2},
        every row #1 column 12/.style={#2},
        every row #1 column 15/.style={#2},
        every row #1 column 18/.style={#2},
    }
}

\pgfplotstabletypeset[
    	col sep=comma,
	    set thousands separator={$\,$},
    	header=has colnames,
    	columns={adjustedInstance,
    	VolumeWeightApproximation-Objective, VolumeWeightApproximation-GapPCT,VolumeWeightApproximation-VehicleCount,
    	LoadingOnly-RPDVolumeWeight,LoadingOnly-GapPCT,LoadingOnly-ADVehicleCount, NoSupport-RPDVolumeWeight,NoSupport-GapPCT, NoSupport-ADVehicleCount, NoLifo-RPDVolumeWeight,NoLifo-GapPCT,NoLifo-ADVehicleCount, NoFragility-RPDVolumeWeight,NoFragility-GapPCT,NoFragility-ADVehicleCount, AllConstraints-RPDVolumeWeight,AllConstraints-GapPCT,AllConstraints-ADVehicleCount},
    	display columns/Instance/.style={string type},
    	column type=r, 
    	precision=2, fixed, fixed zerofill,
        columns/adjustedInstance/.style={string type, column type={l@{\hskip 4pt}}},
        columns/VolumeWeightApproximation-VehicleCount/.style={column type={@{\hskip 6pt}r},precision=0}, 
        columns/AllConstraints-ADVehicleCount/.style={column type={@{\hskip 8pt}r},precision=0},
        columns/NoFragility-ADVehicleCount/.style={column type={@{\hskip 8pt}r},precision=0},
        columns/NoLifo-ADVehicleCount/.style={column type={@{\hskip 8pt}r},precision=0},
        columns/NoSupport-ADVehicleCount/.style={column type={@{\hskip 8pt}r},precision=0},
        columns/LoadingOnly-ADVehicleCount/.style={column type={@{\hskip 8pt}r},precision=0},
        columns/VolumeWeightApproximation-GapPCT/.style={column type={@{\hskip 6pt}r},precision=1,
        }, 
        columns/AllConstraints-RPDVolumeWeight/.style={column type={@{\hskip 8pt}r},precision=1,
        },
        columns/NoFragility-RPDVolumeWeight/.style={column type={@{\hskip 8pt}r},precision=1,
        },
        columns/NoLifo-RPDVolumeWeight/.style={column type={@{\hskip 8pt}r},precision=1,
        },
        columns/NoSupport-RPDVolumeWeight/.style={column type={@{\hskip 8pt}r},precision=1,
        },
        columns/LoadingOnly-RPDVolumeWeight/.style={column type={@{\hskip 8pt}r},precision=1,
        },
        columns/AllConstraints-GapPCT/.style={column type={@{\hskip 6pt}r@{\hskip 2pt}},precision=1,
        },
        columns/NoFragility-GapPCT/.style={column type={@{\hskip 6pt}r},precision=1,
        },
        columns/NoLifo-GapPCT/.style={column type={@{\hskip 6pt}r},precision=1,
        },
        columns/NoSupport-GapPCT/.style={column type={@{\hskip 6pt}r},precision=1,
        },
        columns/LoadingOnly-GapPCT/.style={column type={@{\hskip 6pt}r},precision=1,
        },
        row style = {19}{precision=1, skip 0.},
        every head row/.append style={
            output empty row,
            before row={\toprule
            & \multicolumn{3}{c}{CVRP} & \multicolumn{3}{c}{\loadingOnly{}} & \multicolumn{3}{c}{\noSupport{}} & \multicolumn{3}{c}{\noLifo} & \multicolumn{3}{c}{\noFragility} & \multicolumn{3}{c}{\allConstraints} \\
            \cmidrule(rl{3pt}){2-4} \cmidrule(rl{3pt}){5-7} \cmidrule(rl{3pt}){8-10} \cmidrule(rl{3pt}){11-13} \cmidrule(rl{3pt}){14-16} \cmidrule(rl{3pt}){17-19} \\[-6pt]
            Instance & Obj & Gap & \llap{\#K} 
            & $\Delta_\mathrm{Obj}$ & Gap & \llap{$\Delta_\mathrm{K}$} 
            & $\Delta_\mathrm{Obj}$ & Gap & \llap{$\Delta_\mathrm{K}$}
            & $\Delta_\mathrm{Obj}$ & Gap & \llap{$\Delta_\mathrm{K}$}
            & $\Delta_\mathrm{Obj}$ & Gap & \llap{$\Delta_\mathrm{K}$}
            & $\Delta_\mathrm{Obj}$ & Gap & \llap{$\Delta_\mathrm{K}$}\\
            }, 
            after row={\midrule}},
        every last row/.append style={after row=\bottomrule,before row=\midrule},
    ]{Chapters/AllVariants-TransformedFlat.csv} }
    \begin{tablenotes}[para,flushleft]
        \item[] Comparison is made with relative percentage deviation to the CVRP objective $\Delta_\mathrm{Obj}$ (\%), gap (\%), and absolute deviation of the number of routes $\Delta_\mathrm{K}$.    
    \end{tablenotes}    
    \end{threeparttable}
\end{table}

All instances except instance \texttt{E041-14h} are optimally solved as a CVRP. When different loading variants are compared, a clear pattern emerges: the more the problem deviates from the pure CVRP by introducing loading restrictions, the higher the costs, the more routes are needed, and the more challenging it becomes. The \loadingOnly{} variant exhibits relatively good solvability with an average gap value of 3.4\%, and its solutions have, on average 6.8\% higher costs than those of the CVRP. Of the three loading variants between the extremes of \loadingOnly{} and \allConstraints{}, the \noSupport{} variant exhibits the smallest gap. If the incremental feasibility property does not need to be considered, the solution process becomes easier, as stronger constraints can be added.
However, there are instances---specifically \texttt{E016-05m} and \texttt{E036-11h}---where the CVRP solution is valid for all other loading variants; in such cases, the cuts for fractional solutions are very effective, enabling proof of optimality even for medium-sized instances.
As expected, the \allConstraints{} variant deviates the most from the CVRP, especially for instances with light items. Due to the strongly restricted loading problem, up to four additional vehicles are sometimes required, leading to the highest gap values and more than 30\% higher costs in some instances. The results also highlight that approximating the real loading problem with one-dimensional metrics can very likely generate seemingly appealing solutions in terms of cost that are nevertheless ultimately useless because they are infeasible. Therefore, it is important to assess in advance when approximation is sufficient and when the loading problem needs to be considered in more detail.



\section{Conclusion}\label{sec:Conclusion}




In this paper, we have studied the three-dimensional loading capacitated vehicle routing problem (3L-CVRP), which integrates container loading and vehicle routing. We describe how under certain conditions, previously infeasible loadings can become feasible by adding more items: we call this the \textit{incremental infeasibility} property. To address the latter, we propose novel inequalities and lifting procedures to prohibit infeasible paths, and we integrate them into a branch-and-cut algorithm. We apply several enhancements to speed up the algorithm, including memory management to track previously-evaluated routes and infeasible combinations of customers, a heuristic for the loading problem, cuts for fractional solutions, and an upper bound procedure. In our comprehensive computational study on small to medium-sized benchmark instances with five different loading variants, we show through comparison to a basic algorithm variant that the proposed enhancements are very effective. In total, we are able to solve 56 instances to optimality for the first time, to prove the optimality of 39 current best-known heuristic solutions, and to generate 49 new best solutions. Furthermore, we show that the loading problem has a greater effect in instances with light items, which are more difficult to solve. 

The infeasible path concept can be relatively straightforwardly extended to consider 3L-CVRP variants, such as time windows or pickups and delivery.
With support constraints, feasible routes can become infeasible by removing items and infeasible routes can become feasible by adding items; accordingly, it might be useful to consider the \textit{incremental feasibility} property when designing heuristics.
With regard to branch-and-price-and-cut, which is currently the best exact approach for most other vehicle routing problems, it entails multiple challenges. First, the profitable tour problem with three-dimensional loading constraints (3L-PTP) arises as the pricing problem; the 3L-PTP is not well-studied, and solving it requires solving numerous loading problems with many items and high volume utilization, which are particularly challenging (as demonstrated in this paper). In addition, the \textit{incremental feasibility} property might present a challenge for efficient dominance checks in labeling algorithms. Furthermore, cuts based on one-dimensional capacities like volume or weight can be inefficient; therefore, improved lower bounds for three-dimensional loading problems or better approximations to obtain quicker feasibility checks are worth investigating.




\bibliographystyle{apalike} %
\bibliography{3L-CVRP}  

\begin{thebibliography}{}

\bibitem[Ascheuer et~al., 2000]{Ascheuer2000}
Ascheuer, N., Fischetti, M., and Gr{\"o}tschel, M. (2000).
\newblock A polyhedral study of the asymmetric traveling salesman problem with
  time windows.
\newblock {\em Networks: An International Journal}, 36(2):69--79.

\bibitem[Balas, 1989]{Balas1989}
Balas, E. (1989).
\newblock The asymmetric assignment problem and some new facets of the
  traveling salesman polytope on a directed graph.
\newblock {\em {SIAM} Journal on Discrete Mathematics}, 2(4):425--451.

\bibitem[Balinski and Quandt, 1964]{Balinski1964}
Balinski, M.~L. and Quandt, R.~E. (1964).
\newblock On an integer program for a delivery problem.
\newblock {\em Operations Research}, 12(2):300--304.

\bibitem[Beasley, 1985]{Beasley1985a}
Beasley, J. (1985).
\newblock Algorithms for unconstrained two-dimensional guillotine cutting.
\newblock {\em Journal of the Operational Research Society}, 36(4):297--306.

\bibitem[Bortfeldt, 2012]{Bortfeldt2012}
Bortfeldt, A. (2012).
\newblock A hybrid algorithm for the capacitated vehicle routing problem with
  three-dimensional loading constraints.
\newblock {\em Computers \& Operations Research}, 39(9):2248--2257.

\bibitem[Bortfeldt and Homberger, 2013]{Bortfeldt2013}
Bortfeldt, A. and Homberger, J. (2013).
\newblock Packing first, routing second—a heuristic for the vehicle routing
  and loading problem.
\newblock {\em Computers \& Operations Research}, 40(3):873--885.

\bibitem[Bortfeldt and Wäscher, 2013]{Bortfeldt2013a}
Bortfeldt, A. and Wäscher, G. (2013).
\newblock Constraints in container loading {\textendash} a state-of-the-art
  review.
\newblock {\em European Journal of Operational Research}, 229(1):1--20.

\bibitem[Bortfeldt and Yi, 2020]{Bortfeldt2020}
Bortfeldt, A. and Yi, J. (2020).
\newblock The split delivery vehicle routing problem with three-dimensional
  loading constraints.
\newblock {\em European Journal of Operational Research}, 282(2):545--558.

\bibitem[Boschetti et~al., 2002]{Boschetti2002}
Boschetti, M.~A., Mingozzi, A., and Hadjiconstantinou, E. (2002).
\newblock New upper bounds for the two-dimensional orthogonal non-guillotine
  cutting stock problem.
\newblock {\em IMA Journal of Management Mathematics}, 13(2):95--119.

\bibitem[Ceschia et~al., 2013]{Ceschia2013}
Ceschia, S., Schaerf, A., and St{\"u}tzle, T. (2013).
\newblock Local search techniques for a routing-packing problem.
\newblock {\em Computers \& industrial engineering}, 66(4):1138--1149.

\bibitem[Christofides and Whitlock, 1977]{Christofides1977}
Christofides, N. and Whitlock, C. (1977).
\newblock An algorithm for two-dimensional cutting problems.
\newblock {\em Operations Research}, 25(1):30--44.

\bibitem[C{\^o}t{\'e} et~al., 2017]{Cote2017}
C{\^o}t{\'e}, J.-F., Guastaroba, G., and Speranza, M.~G. (2017).
\newblock The value of integrating loading and routing.
\newblock {\em European Journal of Operational Research}, 257(1):89--105.

\bibitem[C{\^o}t{\'e} and Iori, 2018]{Cote2018}
C{\^o}t{\'e}, J.-F. and Iori, M. (2018).
\newblock The meet-in-the-middle principle for cutting and packing problems.
\newblock {\em INFORMS Journal on Computing}, 30(4):646--661.

\bibitem[Crainic et~al., 2008]{Crainic2008}
Crainic, T.~G., Perboli, G., and Tadei, R. (2008).
\newblock Extreme point-based heuristics for three-dimensional bin packing.
\newblock {\em Informs Journal on computing}, 20(3):368--384.

\bibitem[Escobar et~al., 2015]{Escobar2015}
Escobar, L.~M., Martinez, D.~A., Escobar, J.~W., Linfati, R., and Mauricio,
  G.~E. (2015).
\newblock A hybrid metaheuristic approach for the capacitated vehicle routing
  problem with container loading constraints.
\newblock In {\em 2015 International Conference on Industrial Engineering and
  Systems Management (IESM)}, pages 1374--1382.

\bibitem[Escobar-Falc{\'o}n et~al., 2016]{Escobar2016}
Escobar-Falc{\'o}n, L.~M., {\'A}lvarez-Mart{\'\i}nez, D., Granada-Echeverri,
  M., Escobar, J.~W., and Romero-L{\'a}zaro, R.~A. (2016).
\newblock A matheuristic algorithm for the three-dimensional loading
  capacitated vehicle routing problem ({3L-CVRP}).
\newblock {\em Revista Facultad de Ingenier{\'\i}a Universidad de Antioquia},
  (78):09--20.

\bibitem[Fischetti and Toth, 1997]{Fischetti1997}
Fischetti, M. and Toth, P. (1997).
\newblock A polyhedral approach to the asymmetric traveling salesman problem.
\newblock {\em Management Science}, 43(11):1520--1536.

\bibitem[Fuellerer et~al., 2010]{Fuellerer2010}
Fuellerer, G., Doerner, K.~F., Hartl, R.~F., and Iori, M. (2010).
\newblock Metaheuristics for vehicle routing problems with three-dimensional
  loading constraints.
\newblock {\em European Journal of Operational Research}, 201(3):751--759.

\bibitem[Gendreau et~al., 2006]{Gendreau2006}
Gendreau, M., Iori, M., Laporte, G., and Martello, S. (2006).
\newblock A tabu search algorithm for a routing and container loading problem.
\newblock {\em Transportation Science}, 40(3):342--350.

\bibitem[Gr{\"o}tschel et~al., 1985]{Grotschel1985}
Gr{\"o}tschel, M., Padberg, M.~W., et~al. (1985).
\newblock Polyhedral theory.
\newblock {\em The traveling salesman problem}, pages 251--305.

\bibitem[Hansen et~al., 2017]{Hansen2017}
Hansen, P., Mladenovi{\'c}, N., Todosijevi{\'c}, R., and Hanafi, S. (2017).
\newblock Variable neighborhood search: basics and variants.
\newblock {\em EURO Journal on Computational Optimization}, 5(3):423--454.

\bibitem[Herz, 1972]{Herz1972}
Herz, J. (1972).
\newblock Recursive computational procedure for two-dimensional stock cutting.
\newblock {\em IBM Journal of Research and Development}, 16(5):462--469.

\bibitem[Hokama et~al., 2016]{Hokama2016}
Hokama, P., Miyazawa, F.~K., and Xavier, E.~C. (2016).
\newblock A branch-and-cut approach for the vehicle routing problem with
  loading constraints.
\newblock {\em Expert Systems with Applications}, 47:1--13.

\bibitem[Junqueira et~al., 2013]{Junqueira2013}
Junqueira, L., Oliveira, J.~F., Carravilla, M.~A., and Morabito, R. (2013).
\newblock An optimization model for the vehicle routing problem with practical
  three-dimensional loading constraints.
\newblock {\em International Transactions in Operational Research},
  20(5):645--666.

\bibitem[Koch et~al., 2020]{Koch2020}
Koch, H., Schl{\"o}gell, M., and Bortfeldt, A. (2020).
\newblock A hybrid algorithm for the vehicle routing problem with
  three-dimensional loading constraints and mixed backhauls.
\newblock {\em Journal of Scheduling}, 23(1):71--93.

\bibitem[Kohl et~al., 1999]{Kohl1999}
Kohl, N., Desrosiers, J., Madsen, O.~B., Solomon, M.~M., and Soumis, F. (1999).
\newblock 2-path cuts for the vehicle routing problem with time windows.
\newblock {\em Transportation Science}, 33(1):101--116.

\bibitem[Krebs and Ehmke, 2023]{Krebs2023a}
Krebs, C. and Ehmke, J.~F. (2023).
\newblock Solution validator and visualizer for (combined) vehicle routing and
  container loading problems.
\newblock {\em Annals of Operations Research}.

\bibitem[Krebs et~al., 2021]{Krebs2021}
Krebs, C., Ehmke, J.~F., and Koch, H. (2021).
\newblock Advanced loading constraints for {3D} vehicle routing problems.
\newblock {\em OR Spectrum}, 43(4):835--875.

\bibitem[Kurpel et~al., 2020]{Kurpel2020}
Kurpel, D.~V., Scarpin, C.~T., Junior, J. E.~P., Schenekemberg, C.~M., and
  Coelho, L.~C. (2020).
\newblock The exact solutions of several types of container loading problems.
\newblock {\em European Journal of Operational Research}, 284(1):87--107.

\bibitem[Kü{\c{c}}ük and Yildiz, 2022]{Kuecuek2022}
Kü{\c{c}}ük, M. and Yildiz, S.~T. (2022).
\newblock Constraint programming-based solution approaches for
  three-dimensional loading capacitated vehicle routing problems.
\newblock {\em Computers \& Industrial Engineering}, 171:108505.

\bibitem[Lacomme et~al., 2013]{Lacomme2013}
Lacomme, P., Toussaint, H., and Duhamel, C. (2013).
\newblock A {GRASP}$\times${ELS} for the vehicle routing problem with basic
  three-dimensional loading constraints.
\newblock {\em Engineering Applications of Artificial Intelligence},
  26(8):1795--1810.

\bibitem[Lysgaard et~al., 2004]{Lysgaard2004}
Lysgaard, J., Letchford, A.~N., and Eglese, R.~W. (2004).
\newblock A new branch-and-cut algorithm for the capacitated vehicle routing
  problem.
\newblock {\em Mathematical Programming}, 100(2):423--445.

\bibitem[Mahvash et~al., 2017]{Mahvash2017}
Mahvash, B., Awasthi, A., and Chauhan, S. (2017).
\newblock A column generation based heuristic for the capacitated vehicle
  routing problem with three-dimensional loading constraints.
\newblock {\em International Journal of Production Research}, 55(6):1730--1747.

\bibitem[M{\"a}nnel and Bortfeldt, 2016]{Mannel2016}
M{\"a}nnel, D. and Bortfeldt, A. (2016).
\newblock A hybrid algorithm for the vehicle routing problem with pickup and
  delivery and three-dimensional loading constraints.
\newblock {\em European Journal of Operational Research}, 254(3):840--858.

\bibitem[Martello et~al., 2000]{Martello2000}
Martello, S., Pisinger, D., and Vigo, D. (2000).
\newblock The three-dimensional bin packing problem.
\newblock {\em Operations Research}, 48(2):256--267.

\bibitem[Miao et~al., 2012]{Miao2012}
Miao, L., Ruan, Q., Woghiren, K., and Ruo, Q. (2012).
\newblock A hybrid genetic algorithm for the vehicle routing problem with
  three-dimensional loading constraints.
\newblock {\em RAIRO-Operations Research}, 46(1):63--82.

\bibitem[Moura and Oliveira, 2009]{Moura2009}
Moura, A. and Oliveira, J.~F. (2009).
\newblock An integrated approach to the vehicle routing and container loading
  problems.
\newblock {\em {OR} Spectrum}, 31(4):775--800.

\bibitem[Pollaris et~al., 2015]{Pollaris2015}
Pollaris, H., Braekers, K., Caris, A., Janssens, G.~K., and Limbourg, S.
  (2015).
\newblock Vehicle routing problems with loading constraints: state-of-the-art
  and future directions.
\newblock {\em OR Spectrum}, 37(2):297--330.

\bibitem[Rajaei et~al., 2022]{Rajaei2022}
Rajaei, M., Moslehi, G., and Reisi-Nafchi, M. (2022).
\newblock The split heterogeneous vehicle routing problem with
  three-dimensional loading constraints on a large scale.
\newblock {\em European Journal of Operational Research}, 299(2):706--721.

\bibitem[Reil et~al., 2018]{Reil2018}
Reil, S., Bortfeldt, A., and M{\"o}nch, L. (2018).
\newblock Heuristics for vehicle routing problems with backhauls, time windows,
  and {3D} loading constraints.
\newblock {\em European Journal of Operational Research}, 266(3):877--894.

\bibitem[Ropke et~al., 2007]{ropkes-2007}
Ropke, S., Cordeau, J.-F., and Laporte, G. (2007).
\newblock Models and branch-and-cut algorithms for pickup and delivery problems
  with time windows.
\newblock {\em Networks}, 49(4):258--272.

\bibitem[Ruan et~al., 2013]{Ruan2013}
Ruan, Q., Zhang, Z., Miao, L., and Shen, H. (2013).
\newblock A hybrid approach for the vehicle routing problem with
  three-dimensional loading constraints.
\newblock {\em Computers \& Operations Research}, 40(6):1579--1589.

\bibitem[Tamke et~al., 2024]{tamkef-2024}
Tamke, F., Linß, F., and Kuttner, L. (2024).
\newblock {3L-CVRP}.

\bibitem[Tao and Wang, 2015]{Tao2015}
Tao, Y. and Wang, F. (2015).
\newblock An effective tabu search approach with improved loading algorithms
  for the {3L-CVRP}.
\newblock {\em Computers \& Operations Research}, 55:127--140.

\bibitem[Tarantilis et~al., 2009]{Tarantilis2009}
Tarantilis, C.~D., Zachariadis, E.~E., and Kiranoudis, C.~T. (2009).
\newblock A hybrid metaheuristic algorithm for the integrated vehicle routing
  and three-dimensional container-loading problem.
\newblock {\em IEEE Transactions on Intelligent Transportation Systems},
  10(2):255--271.

\bibitem[uit~het Broek et~al., 2021]{Broek2021}
uit~het Broek, M. A.~J., Schrotenboer, A.~H., Jargalsaikhan, B., Roodbergen,
  K.~J., and Coelho, L.~C. (2021).
\newblock Asymmetric multidepot vehicle routing problems: Valid inequalities
  and a branch-and-cut algorithm.
\newblock {\em Operations Research}, 69(2):380--409.

\bibitem[Vega-Mej{\'\i}a et~al., 2019]{Vega2019}
Vega-Mej{\'\i}a, C.~A., Montoya-Torres, J.~R., and Islam, S.~M. (2019).
\newblock A nonlinear optimization model for the balanced vehicle routing
  problem with loading constraints.
\newblock {\em International Transactions in Operational Research},
  26(3):794--835.

\bibitem[Wang et~al., 2010]{Wang2010}
Wang, L., Guo, S., Chen, S., Zhu, W., and Lim, A. (2010).
\newblock Two natural heuristics for {3D} packing with practical loading
  constraints.
\newblock In {\em Pacific Rim International Conference on Artificial
  Intelligence}, pages 256--267. Springer.

\bibitem[Wei et~al., 2014]{Wei2014}
Wei, L., Zhang, Z., and Lim, A. (2014).
\newblock An adaptive variable neighborhood search for a heterogeneous fleet
  vehicle routing problem with three-dimensional loading constraints.
\newblock {\em IEEE Computational Intelligence Magazine}, 9(4):18--30.

\bibitem[Zhang et~al., 2022a]{Zhang2022a}
Zhang, X., Chen, L., Gendreau, M., and Langevin, A. (2022a).
\newblock A branch-and-cut algorithm for the vehicle routing problem with
  two-dimensional loading constraints.
\newblock {\em European Journal of Operational Research}, 302(1):259--269.

\bibitem[Zhang et~al., 2022b]{Zhang2022}
Zhang, X., Chen, L., Gendreau, M., and Langevin, A. (2022b).
\newblock A branch-and-price-and-cut algorithm for the vehicle routing problem
  with two-dimensional loading constraints.
\newblock {\em Transportation Science}, 56(6):1618--1635.

\bibitem[Zhang et~al., 2015]{Zhang2015}
Zhang, Z., Wei, L., and Lim, A. (2015).
\newblock An evolutionary local search for the capacitated vehicle routing
  problem minimizing fuel consumption under three-dimensional loading
  constraints.
\newblock {\em Transportation Research Part B: Methodological}, 82:20--35.

\bibitem[Zhu et~al., 2012]{Zhu2012}
Zhu, W., Qin, H., Lim, A., and Wang, L. (2012).
\newblock A two-stage tabu search algorithm with enhanced packing heuristics
  for the {3L-CVRP} and {M3L-CVRP}.
\newblock {\em Computers \& Operations Research}, 39(9):2178--2195.

\end{thebibliography}

\clearpage
\section*{Appendix A: Flow charts for the exact feasibility check of different loading variants}\label{sec:ExactFeasibilityCheck}

The exact feasibility check for the loading is adapted to the different loading variants. In particular, the considered loading variants in the lifting steps are adapted and unnecessary lifting steps are skipped. The following flow charts illustrate the feasibility checks.

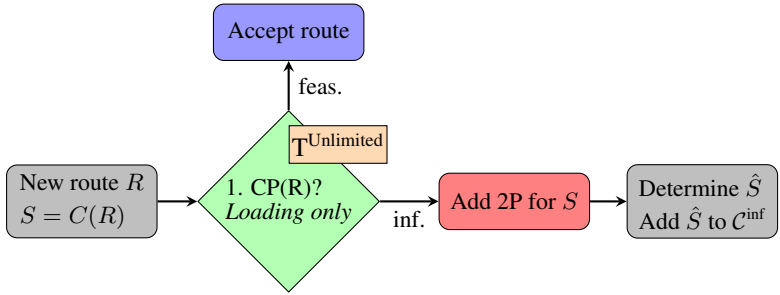
\begin{figure}[h]
    \tikzstyle{start} = [rectangle, rounded corners, 
minimum width=2cm, 
minimum height=0.75cm,
font=\footnotesize,
text centered,
align=left,
draw=black, 
fill=gray!50]

\tikzstyle{stopPositive} = [rectangle, rounded corners, 
minimum width=2cm, 
minimum height=0.75cm,
font=\footnotesize,
text centered, 
draw=black, 
fill=blue!40]

\tikzstyle{stopNegative} = [rectangle, rounded corners, 
minimum width=2cm, 
minimum height=0.75cm,
font=\footnotesize\linespread{0.9}\selectfont,
text centered, 
draw=black, 
fill=red!50]

\tikzstyle{io} = [trapezium, 
trapezium stretches=true, 
trapezium left angle=70, 
trapezium right angle=110, 
minimum width=3cm, 
minimum height=1cm, text centered, 
draw=black, fill=blue!30]

\tikzstyle{process} = [rectangle,
minimum width=1cm, 
minimum height=0.5cm,
inner sep=1pt,
text centered, 
draw=black, 
fill=orange!30]

\tikzstyle{decision} = [diamond, 
font=\footnotesize\linespread{0.9}\selectfont,
inner sep=1pt,
text centered, 
draw=black,
fill=green!30]
\tikzstyle{arrow} = [thick,->,>=stealth]
\tikzstyle{noArrow} = [thick,>=stealth]

\begin{tikzpicture}[node distance=2cm]

\node (Start) [start] {New route $R$ \\ ${S}={C}(R)$};







\node (CP) [decision, right of=Start, align=left,xshift=0.75cm] {1. CP(R)?\\\textit{Loading only}};
\node [process, above of=CP,anchor=south west, yshift=-1.5cm] {T\textsuperscript{Unlimited}};
\node (Accept) [stopPositive, above of=CP, yshift=0.25cm] {Accept route};


\node (ITPC) [stopNegative, right of=CP, align=left,xshift=1cm] {Add 2P for $S$};
\node (RedInfSet) [start, right of=ITPC, align=left,xshift=0.5cm] {Determine $\hat{S}$ \\ Add $\hat{S}$ to $\mathcal{C}^{\text{inf}}$};

\draw [arrow] (Start) -- (CP);




\draw [arrow] (CP) -- node[font=\footnotesize,anchor=west,align=left] {feas.} (Accept);
\draw [arrow] (CP) -- node[font=\footnotesize\linespread{0.9}\selectfont,anchor=north,align=left] {inf.} (ITPC);

\draw [arrow] (ITPC) -- (RedInfSet);






\end{tikzpicture}
    \caption{Flow chart of the exact feasibility check of a route with a CP model considering the \textit{loading only} variant. \label{fig:FlowChartExactRouteCheck_NoSupport}}    
\end{figure}

\begin{figure}[h]    
    \tikzstyle{start} = [rectangle, rounded corners, 
minimum width=2cm, 
minimum height=0.75cm,
font=\footnotesize,
text centered,
align=left,
draw=black, 
fill=gray!50]

\tikzstyle{stopPositive} = [rectangle, rounded corners, 
minimum width=2cm, 
minimum height=0.75cm,
font=\footnotesize,
text centered, 
draw=black, 
fill=blue!40]

\tikzstyle{stopNegative} = [rectangle, rounded corners, 
minimum width=2cm, 
minimum height=0.75cm,
font=\footnotesize\linespread{0.9}\selectfont,
text centered, 
draw=black, 
fill=red!50]

\tikzstyle{io} = [trapezium, 
trapezium stretches=true, 
trapezium left angle=70, 
trapezium right angle=110, 
minimum width=3cm, 
minimum height=1cm, text centered, 
draw=black, fill=blue!30]

\tikzstyle{process} = [rectangle,
minimum width=1cm, 
minimum height=0.5cm,
inner sep=1pt,
text centered, 
draw=black, 
fill=orange!30]

\tikzstyle{decision} = [diamond, 
font=\footnotesize\linespread{0.9}\selectfont,
inner sep=1pt,
text centered, 
draw=black,
fill=green!30]
\tikzstyle{arrow} = [thick,->,>=stealth]
\tikzstyle{noArrow} = [thick,>=stealth]

\begin{tikzpicture}[node distance=2cm]

\node (Start) [start] {New route $R$ \\ ${S}={C}(R)$};

\node (CPLimited) [decision, right of=Start, align=left,xshift=0.75cm] {1. CP(R)?\\\textit{NoLIFO}};
\node [process, above of=CPLimited, anchor=south west, yshift=-1.5cm] {T\textsuperscript{Limited}};
\node (AcceptLimited) [stopPositive, above of=CPLimited, yshift=0.25cm] {Accept route};

\node (TwoPC) [decision, right of=CPLimited, align=left,xshift=0.75cm] {2. CP(R)?\\\textit{LoadingOnly}};
\node [process, above of=TwoPC, anchor=south west, yshift=-1.5cm] {T\textsuperscript{Limited}};
\node (TwoPCResult) [stopNegative, below of=TwoPC, align=left] {Add 2P for ${S}$};
\node (RedInfSet) [start, below of=TwoPCResult, yshift=0.75cm,align=left] {Determine $\hat{S}$ \\ Add $\hat{S}$ to $\mathcal{C}^{\text{inf}}$};



\node (CPResult) [decision, right of=TwoPC, align=left,xshift=0.75cm] {1. CP(R)?};

\node (CP) [decision, right of=CPResult, align=left,xshift=0.75cm] {3. CP(R)?\\\textit{NoLIFO}};
\node [process, above of=CP,anchor=south west, yshift=-1.5cm] {T\textsuperscript{Unlimited}};
\node (Accept) [stopPositive, above of=CP, yshift=0.25cm] {Accept route};


\node (ITPC) [stopNegative, below of=CP, align=left,xshift=0.5cm, yshift=-0.5cm] {Add 2PT for $S$};

\draw [arrow] (Start) -- (CPLimited);

\draw [arrow] (CPLimited) -- node[font=\footnotesize,anchor=west,align=left] {feas.} (AcceptLimited);
\draw [arrow] (CPLimited) -- node[font=\footnotesize\linespread{0.9}\selectfont,below,align=left] {inf./\\ unk.} (TwoPC);

\draw [arrow] (TwoPC) -- node[font=\footnotesize,anchor=west] {inf.} (TwoPCResult);
\draw [arrow] (TwoPCResult) -- (RedInfSet);
\draw [arrow] (TwoPC) -- node[font=\footnotesize\linespread{0.9}\selectfont,anchor=north, align=left] {feas./\\ unk.} (CPResult);


\draw [arrow] (CPResult) -- node[font=\footnotesize,anchor=north,align=left] {unk.} (CP);
\draw [arrow] (CP) -- node[font=\footnotesize,anchor=west,align=left] {feas.} (Accept);
\draw [noArrow] (CP.south) -- node[font=\footnotesize\linespread{0.9}\selectfont,anchor=west,align=left] {inf.} ++ (0,-0.25) -- ++ (0,0) -| (CPResult.south);



\draw [arrow] (CPResult.south) -- ++ (0,-0.5) -- node[font=\footnotesize\linespread{0.9}\selectfont,anchor=east,align=left] {inf.} ++ (0,0) |- (ITPC.west) ;



\begin{scope}[on background layer]
    \fill [draw=black, dashed, rounded corners, fill=gray!10, opacity=0.5] ($(AcceptLimited.west |- AcceptLimited.north) + (-0.2,0.1)$) 
    -- ($(CPResult.east |-  AcceptLimited.north) + (0,0.1)$) 
    -- ($(CPResult.east |- RedInfSet.north) + (0.1,0)$) 
    -- ($(RedInfSet.east |- RedInfSet.north) + (0.2,0)$) 
    -- ($(RedInfSet.east |- RedInfSet.south) + (0.2,-0.1)$) 
    -- ($(AcceptLimited.west |- RedInfSet.south) + (-0.2,-0.1)$) -- cycle;   
   \node[label={above right:Lifting}] at (CPLimited.west |- RedInfSet.south) {};
\end{scope}
\end{tikzpicture}
    \caption{Flow chart of the exact feasibility check of a route with a CP model considering the \textit{no LIFO} variant. \label{fig:FlowChartExactRouteCheck_NoSupport}}    
\end{figure}
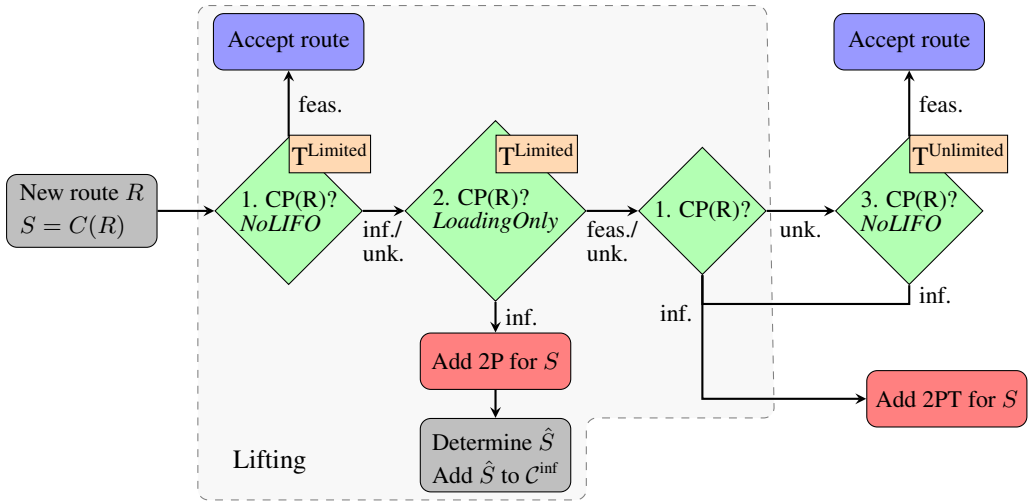

\begin{figure}[h]    
    {\tikzstyle{start} = [rectangle, rounded corners, 
minimum width=2cm, 
minimum height=0.75cm,
font=\footnotesize,
text centered,
align=left,
draw=black, 
fill=gray!50]

\tikzstyle{stopPositive} = [rectangle, rounded corners, 
minimum width=2cm, 
minimum height=0.75cm,
font=\footnotesize,
text centered, 
draw=black, 
fill=blue!40]

\tikzstyle{stopNegative} = [rectangle, rounded corners, 
minimum width=2cm, 
minimum height=0.75cm,
font=\footnotesize\linespread{0.9}\selectfont,
text centered, 
draw=black, 
fill=red!50]

\tikzstyle{io} = [trapezium, 
trapezium stretches=true, 
trapezium left angle=70, 
trapezium right angle=110, 
minimum width=3cm, 
minimum height=1cm, text centered, 
draw=black, fill=blue!30]

\tikzstyle{process} = [rectangle,
minimum width=1cm, 
minimum height=0.5cm,
inner sep=1pt,
text centered, 
draw=black, 
fill=orange!30]

\tikzstyle{decision} = [diamond, 
font=\footnotesize\linespread{0.9}\selectfont,
inner sep=1pt,
text centered, 
draw=black,
fill=green!30]
\tikzstyle{arrow} = [thick,->,>=stealth]

\begin{tikzpicture}[node distance=2cm]

\node (Start) [start] {New route $R$ \\ ${S}={C}(R)$};

\node (CPLimited) [decision, right of=Start, align=left,xshift=0.75cm] {1. CP(R)?\\\textit{NoSupport}};
\node [process, above of=CPLimited, anchor=south west, yshift=-1.5cm] {T\textsuperscript{Limited}};
\node (AcceptLimited) [stopPositive, above of=CPLimited, yshift=0.25cm] {Accept route};

\node (TwoPC) [decision, right of=CPLimited, align=left,xshift=0.75cm] {2. CP(R)?\\\textit{LIFO\textsuperscript{NoSeq}}};
\node [process, above of=TwoPC, anchor=south west, yshift=-1.5cm] {T\textsuperscript{Limited}};
\node (TwoPCResult) [stopNegative, below of=TwoPC, yshift=-0.25cm,align=left] {Add 2P for ${S}$};
\node (RedInfSet) [start, below of=TwoPCResult, yshift=0.25cm,align=left] {Determine $\hat{S}$ \\ Add $\hat{S}$ to $\mathcal{C}^{\text{inf}}$};




\node (CPResult) [decision, right of=TwoPC, align=left,xshift=0.75cm] {1. CP(R)?};

\node (CP) [decision, right of=CPResult, align=left,xshift=0.75cm] {3. CP(R)?\\\textit{NoSupport}};
\node [process, above of=CP,anchor=south west, yshift=-1.5cm] {T\textsuperscript{Unlimited}};
\node (Accept) [stopPositive, above of=CP, yshift=0.25cm] {Accept route};

\node (RevPath) [decision, below of=CPResult, align=left,yshift=-2.5cm] {4. CP($\bar{R}$)?\\\textit{NoSupport}};
\node [process, above of=RevPath, anchor=south west, yshift=-1.5cm] {T\textsuperscript{Limited}};
\node (RevPathResult) [stopNegative, right of=RevPath,align=left, xshift=1cm] {Add UIP for $P$\\+ RT for $\bar{P}$};

\node (ITPC) [stopNegative, below of=CP, align=left,xshift=1cm, yshift=-0.5cm] {Add RT for $\hat{P}$};
\node (RedInfPath) [start, left of=ITPC, xshift=-0.3cm, align=left] {Determine $\hat{P}$};

\draw [arrow] (Start) -- (CPLimited);

\draw [arrow] (CPLimited) -- node[font=\footnotesize,anchor=west,align=left] {feas.} (AcceptLimited);
\draw [arrow] (CPLimited) -- node[font=\footnotesize\linespread{0.9}\selectfont,below,align=left] {inf./\\ unk.} (TwoPC);

\draw [arrow] (TwoPC) -- node[font=\footnotesize,anchor=west] {inf.} (TwoPCResult);
\draw [arrow] (TwoPCResult) -- (RedInfSet);
\draw [arrow] (TwoPC) -- node[font=\footnotesize\linespread{0.9}\selectfont,anchor=north, align=left] {feas./\\ unk.} (CPResult);


\draw [arrow] (CPResult) -- node[font=\footnotesize,anchor=north,align=left] {unk.} (CP);
\draw [arrow] (CP) -- node[font=\footnotesize,anchor=west,align=left] {feas.} (Accept);
\draw [arrow] (CP.south) -- node[font=\footnotesize\linespread{0.9}\selectfont,anchor=west,align=left] {inf.} ++ (0,-0.25) -- ++ (0,0) -| (RevPath.north);

\draw [arrow] (CPResult) -- node[font=\footnotesize\linespread{0.9}\selectfont,anchor=east,align=left] {inf.} (RevPath);

\draw [arrow] (RevPath) -- node[font=\footnotesize,anchor=north,pos=0.4] {inf.} (RevPathResult);
\draw [arrow] (CPResult.south) -- ++ (0,-0.5) -- ++ (0,0) |- (RedInfPath.west);
\draw [arrow] (RedInfPath) -- (ITPC);


\begin{scope}[on background layer]
    \fill [draw=black, dashed, rounded corners, fill=gray!10, opacity=0.5] ($(AcceptLimited.west |- AcceptLimited.north) + (-0.2,0.1)$) -- ($(RevPath.east |-  AcceptLimited.north) + (0,0.1)$) -- ($(RevPath.east |- RevPath.north) + (0,0)$) -- ($(RevPathResult.east |- RevPath.north) + (0.1,0)$) -- ($(RevPathResult.east |- RevPath.south) + (0.1,-0.1)$) -- ($(AcceptLimited.west |- RevPath.south) + (-0.2,-0.1)$) -- cycle;
   -- cycle;
   \node[label={above right:Lifting}] at (CPLimited.west |- RevPath.south) {};
\end{scope}
\end{tikzpicture}}
    \caption{Flow chart of the exact feasibility check of a route with a CP model considering the \textit{no support} variant. \label{fig:FlowChartExactRouteCheck_NoSupport}}    
\end{figure}
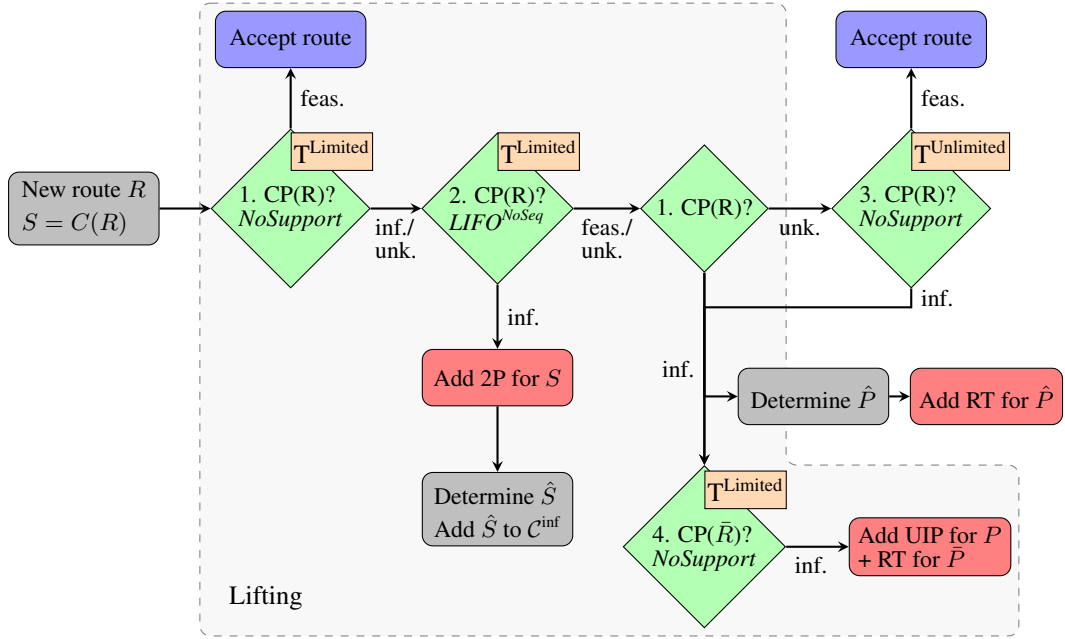

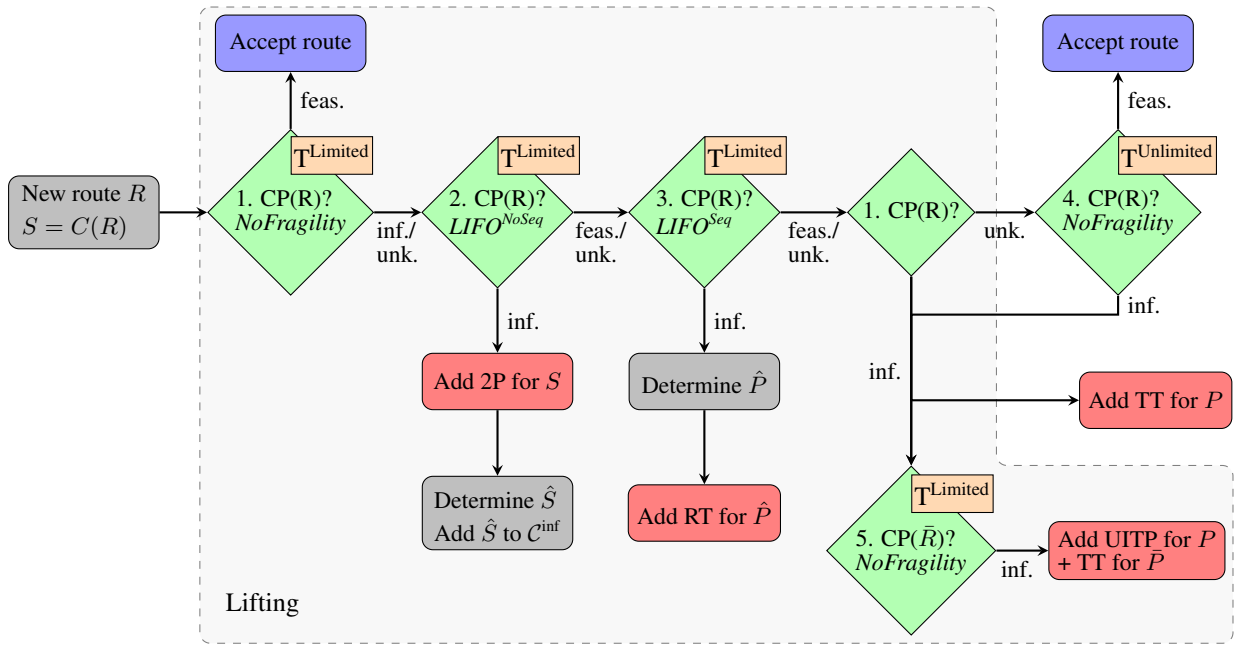
\begin{figure}[h]
    \tikzstyle{start} = [rectangle, rounded corners, 
minimum width=2cm, 
minimum height=0.75cm,
font=\footnotesize,
text centered,
align=left,
draw=black, 
fill=gray!50]

\tikzstyle{stopPositive} = [rectangle, rounded corners, 
minimum width=2cm, 
minimum height=0.75cm,
font=\footnotesize,
text centered, 
draw=black, 
fill=blue!40]

\tikzstyle{stopNegative} = [rectangle, rounded corners, 
minimum width=2cm, 
minimum height=0.75cm,
font=\footnotesize\linespread{0.9}\selectfont,
text centered, 
draw=black, 
fill=red!50]

\tikzstyle{io} = [trapezium, 
trapezium stretches=true, 
trapezium left angle=70, 
trapezium right angle=110, 
minimum width=3cm, 
minimum height=1cm, text centered, 
draw=black, fill=blue!30]

\tikzstyle{process} = [rectangle,
minimum width=1cm, 
minimum height=0.5cm,
inner sep=1pt,
text centered, 
draw=black, 
fill=orange!30]

\tikzstyle{decision} = [diamond, 
font=\footnotesize\linespread{0.9}\selectfont,
inner sep=1pt,
text centered, 
draw=black,
fill=green!30]
\tikzstyle{arrow} = [thick,->,>=stealth]

\begin{tikzpicture}[node distance=2cm]

\node (Start) [start] {New route $R$ \\ ${S}={C}(R)$};

\node (CPLimited) [decision, right of=Start, align=left,xshift=0.75cm] {1. CP(R)?\\\textit{NoFragility}};
\node [process, above of=CPLimited, anchor=south west, yshift=-1.5cm] {T\textsuperscript{Limited}};
\node (AcceptLimited) [stopPositive, above of=CPLimited, yshift=0.25cm] {Accept route};

\node (TwoPC) [decision, right of=CPLimited, align=left,xshift=0.75cm] {2. CP(R)?\\\textit{LIFO\textsuperscript{NoSeq}}};
\node [process, above of=TwoPC, anchor=south west, yshift=-1.5cm] {T\textsuperscript{Limited}};
\node (TwoPCResult) [stopNegative, below of=TwoPC, yshift=-0.25cm,align=left] {Add 2P for ${S}$};
\node (RedInfSet) [start, below of=TwoPCResult, yshift=0.25cm,align=left] {Determine $\hat{S}$ \\ Add $\hat{S}$ to $\mathcal{C}^{\text{inf}}$};

\node (IPC) [decision, right of=TwoPC, align=left,xshift=0.75cm] {3. CP(R)?\\\textit{LIFO\textsuperscript{Seq}}};
\node [process, above of=IPC, anchor=south west, yshift=-1.5cm] {T\textsuperscript{Limited}};
\node (RedInfPath) [start, below of=IPC, yshift=-0.25cm,align=left] {Determine $\hat{P}$};
\node (IPCResult) [stopNegative, below of=RedInfPath, yshift=0.25cm,align=left] {Add RT for $\hat{P}$};

\node (CPResult) [decision, right of=IPC, align=left,xshift=0.75cm] {1. CP(R)?};

\node (CP) [decision, right of=CPResult, align=left,xshift=0.75cm] {4. CP(R)?\\\textit{NoFragility}};
\node [process, above of=CP,anchor=south west, yshift=-1.5cm] {T\textsuperscript{Unlimited}};
\node (Accept) [stopPositive, above of=CP, yshift=0.25cm] {Accept route};

\node (RevPath) [decision, below of=CPResult, align=left,yshift=-2.5cm] {5. CP($\bar{R}$)?\\\textit{NoFragility}};
\node [process, above of=RevPath, anchor=south west, yshift=-1.5cm] {T\textsuperscript{Limited}};
\node (RevPathResult) [stopNegative, right of=RevPath,align=left, xshift=1cm] {Add UITP for $P$\\+ TT for $\bar{P}$};

\node (ITPC) [stopNegative, below of=CP, align=left,xshift=0.5cm, yshift=-0.5cm] {Add TT for $P$};

\draw [arrow] (Start) -- (CPLimited);

\draw [arrow] (CPLimited) -- node[font=\footnotesize,anchor=west,align=left] {feas.} (AcceptLimited);
\draw [arrow] (CPLimited) -- node[font=\footnotesize\linespread{0.9}\selectfont,below,align=left] {inf./\\ unk.} (TwoPC);

\draw [arrow] (TwoPC) -- node[font=\footnotesize,anchor=west] {inf.} (TwoPCResult);
\draw [arrow] (TwoPCResult) -- (RedInfSet);
\draw [arrow] (TwoPC) -- node[font=\footnotesize\linespread{0.9}\selectfont,anchor=north, align=left] {feas./\\ unk.} (IPC);

\draw [arrow] (IPC) -- node[font=\footnotesize,anchor=west] {inf.} (RedInfPath);
\draw [arrow] (RedInfPath) -- (IPCResult);
\draw [arrow] (IPC) -- node[font=\footnotesize\linespread{0.9}\selectfont,anchor=north,align=left] {feas./\\ unk.} (CPResult);

\draw [arrow] (CPResult) -- node[font=\footnotesize,anchor=north,align=left] {unk.} (CP);
\draw [arrow] (CP) -- node[font=\footnotesize,anchor=west,align=left] {feas.} (Accept);
\draw [arrow] (CP.south) -- node[font=\footnotesize\linespread{0.9}\selectfont,anchor=west,align=left] {inf.} ++ (0,-0.25) -- ++ (0,0) -| (RevPath.north);

\draw [arrow] (CPResult) -- node[font=\footnotesize\linespread{0.9}\selectfont,anchor=east,align=left] {inf.} (RevPath);

\draw [arrow] (RevPath) -- node[font=\footnotesize,anchor=north,pos=0.4] {inf.} (RevPathResult);
\draw [arrow] (CPResult.south) -- ++ (0,-0.5) -- ++ (0,0) |- (ITPC.west);


\begin{scope}[on background layer]
    \fill [draw=black, dashed, rounded corners, fill=gray!10, opacity=0.5] ($(AcceptLimited.west |- AcceptLimited.north) + (-0.2,0.1)$) -- ($(RevPath.east |-  AcceptLimited.north) + (0,0.1)$) -- ($(RevPath.east |- RevPath.north) + (0,0)$) -- ($(RevPathResult.east |- RevPath.north) + (0.1,0)$) -- ($(RevPathResult.east |- RevPath.south) + (0.1,-0.1)$) -- ($(AcceptLimited.west |- RevPath.south) + (-0.2,-0.1)$) -- cycle;
   -- cycle;
   \node[label={above right:Lifting}] at (CPLimited.west |- RevPath.south) {};
\end{scope}
\end{tikzpicture}
    \caption{Flow chart of the exact feasibility check of a route with a CP model considering the \textit{no fragility} variant. The relaxation \textit{LIFO\textsuperscript{Seq}} results from the lifting of the infeasible tail path inequalities by relaxing the support constraints \label{fig:FlowChartExactRouteCheck_NoSupport}}
\end{figure}

\end{document}